\documentclass[12pt]{amsart}      
\usepackage{amssymb}
\usepackage{eucal}
\usepackage{graphicx}
\usepackage{amsmath}
\usepackage{amscd}
\usepackage[all]{xy}           
\usepackage{longtable, lscape}
\usepackage{rotating}

\usepackage{amsfonts,latexsym}
\usepackage{xspace}
\usepackage{epsfig}
\usepackage{float}
\usepackage{color}
\usepackage{fancybox}
\usepackage{colordvi}
\usepackage{multicol}

\newtheorem{theorem}{Theorem}[section]
\newtheorem{mprop}[theorem]{Proposition}
\newtheorem{mdefn}[theorem]{Definition}
\newtheorem{lemma}[theorem]{Lemma}
\newtheorem{coro}[theorem]{Corollary}
\newtheorem{thm-def}[theorem]{Theorem-Definition}
\newtheorem{def-prop}[theorem]{Definition-Proposition}
\newtheorem{prop-def}[theorem]{Proposition-Definition}
\newtheorem{coro-def}[theorem]{Corollary-Definition}

\newcommand{\nc}{\newcommand}
\nc{\tred}[1]{\textcolor{red}{#1}}
\nc{\tblue}[1]{\textcolor{blue}{#1}}
\nc{\tgreen}[1]{\textcolor{green}{#1}}
\nc{\tpurple}[1]{\textcolor{purple}{#1}}
\nc{\btred}[1]{\textcolor{red}{\bf #1}}
\nc{\btblue}[1]{\textcolor{blue}{\bf #1}}
\nc{\btgreen}[1]{\textcolor{green}{\bf #1}}
\nc{\btpurple}[1]{\textcolor{purple}{\bf #1}}

\renewcommand{\Bbb}{\mathbb}
\renewcommand{\frak}{\mathfrak}
\newcommand{\efootnote}[1]{}
\renewcommand{\textbf}[1]{}
\newcommand{\delete}[1]{}
\nc{\dfootnote}[1]{{}}          
\nc{\ffootnote}[1]{\dfootnote{#1}}
\nc{\mfootnote}[1]{\footnote{#1}} 
\nc{\ofootnote}[1]{\footnote{\tiny Older version: #1}} 

\delete{
\nc{\mlabel}[1]{\label{#1}  
{\hfill \hspace{1cm}{\bf{{\ }\hfill(#1)}}}}
\nc{\mcite}[1]{\cite{#1}{{\bf{{\ }(#1)}}}}  
\nc{\mref}[1]{\ref{#1}{{\bf{{\ }(#1)}}}}  
\nc{\mbibitem}[1]{\bibitem[\bf #1]{#1}} 
}

\nc{\mlabel}[1]{\label{#1}}  
\nc{\mcite}[1]{\cite{#1}}  
\nc{\mref}[1]{\ref{#1}}  
\nc{\mbibitem}[1]{\bibitem{#1}} 

\nc{\sbar}{, }
\delete{
{\, {\scriptstyle{|\hspace{-.08cm}|\hspace{-.08cm}|
\hspace{-.08cm}|\hspace{-.08cm}|\hspace{-.08cm}|
\hspace{-.08cm}|\hspace{-.08cm}|\hspace{-.08cm}|\hspace{-.08cm}|
\hspace{-.08cm}|\hspace{-.08cm}|\hspace{-.08cm}|
\hspace{-.08cm}|\hspace{-.08cm}|\hspace{-.08cm}|\hspace{-.08cm}|
\hspace{-.08cm}|\hspace{-.08cm}|\hspace{-.08cm}|
\hspace{-.08cm}|\hspace{-.08cm}|\hspace{-.08cm}|\hspace{-.08cm}|
\hspace{-.08cm}|\hspace{-.08cm}|\hspace{-.08cm}|
\hspace{-.08cm}|\hspace{-.08cm}|}\, }}
}

\nc{\wvec}[2]{{\scriptsize{ [
    \begin{array}{c} #1 \\ #2 \end{array}   ]}}}
\nc{\lp}{\big ( }
\nc{\llp}{\Big (}
\nc{\Llp}{\left (}
\nc{\rp}{\big ) }
\nc{\rrp}{\Big )}
\nc{\Rrp}{\right )}
\nc{\lb}{\big < }
\nc{\llb}{\!\Big \langle }
\nc{\Llb}{\! \left <}
\nc{\rb}{\big >  }
\nc{\rrb}{\Big \rangle \!}
\nc{\Rb}{\Big \rangle\! }
\nc{\length}{{\rm leng}}
\nc{\cop}{{\rm cop}}

\nc{\bin}[2]{ (_{\stackrel{\scs{#1}}{\scs{#2}}})}  
\nc{\binc}[2]{ \big (\! \begin{array}{c} \scs{#1}\\
    \scs{#2} \end{array}\! \big )}  
\nc{\bincc}[2]{  \left ( {\scs{#1} \atop
    \vspace{-1cm}\scs{#2}} \right )}  
\nc{\bs}{\bar{S}}
\nc{\cosum}{\sqsubset}
\nc{\la}{\longrightarrow}
\nc{\rar}{\rightarrow}
\nc{\dar}{\downarrow}
\nc{\dap}[1]{\downarrow \rlap{$\scriptstyle{#1}$}}
\nc{\uap}[1]{\uparrow \rlap{$\scriptstyle{#1}$}}
\nc{\defeq}{\stackrel{\rm def}{=}}
\nc{\disp}[1]{\displaystyle{#1}}
\nc{\dotcup}{\ \displaystyle{\bigcup^\bullet}\ }
\nc{\gzeta}{\bar{\zeta}}
\nc{\hcm}{\ \hat{,}\ }
\nc{\hts}{\hat{\otimes}}
\nc{\barot}{{\otimes}}
\nc{\free}[1]{\bar{#1}}
\nc{\uni}[1]{\tilde{#1}}          
\nc{\hcirc}{\hat{\circ}}
\nc{\lleft}{[}
\nc{\lright}{]}
\nc{\curlyl}{\left \{ \begin{array}{c} {} \\ {} \end{array}
    \right .  \!\!\!\!\!\!\!}
\nc{\curlyr}{ \!\!\!\!\!\!\!
    \left . \begin{array}{c} {} \\ {} \end{array}
    \right \} }
\nc{\longmid}{\left | \begin{array}{c} {} \\ {} \end{array}
    \right . \!\!\!\!\!\!\!}
\nc{\ora}[1]{\stackrel{#1}{\rar}}
\nc{\ola}[1]{\stackrel{#1}{\la}}
\nc{\ot}{\otimes}
\nc{\mot}{{{\sbar}}}
\nc{\otm}{\mot}
\nc{\scs}[1]{\scriptstyle{#1}}
\nc{\subv}{{^{\star}}}
\nc{\cov}{{^{\sharp}}}
\nc{\mrm}[1]{{\rm #1}}
\nc{\dirlim}{\displaystyle{\lim_{\longrightarrow}}\,}
\nc{\invlim}{\displaystyle{\lim_{\longleftarrow}}\,}
\nc{\proofbegin}{\noindent{\bf Proof: }}
\nc{\proofend}{$\quad \square$ \vspace{0.3cm}}
\nc{\sha}{{\mbox{\cyr X}}}  
\nc{\shap}{{\mbox{\cyrs X}}} 
\nc{\shpr}{\diamond}    
\nc{\shplus}{\shpr^+}
\nc{\shprc}{\shpr_c}    
\nc{\msh}{\ast}
\nc{\vep}{\varepsilon}
\nc{\labs}{\mid\!}
\nc{\rabs}{\!\mid}
\nc{\FG}{\mrm{FG}}
\nc{\fp}{\tilde{P}} \nc{\rchar}{\mrm{char}} \nc{\Fil}{\mrm{Fil}}
\nc{\gmzvs}{gMZV\xspace}
\nc{\gmzv}{gMZV\xspace}
\nc{\mzv}{MZV\xspace}
\nc{\mzvs}{MZVs\xspace}
\nc{\MZV}{\mrm{MZV}}
\nc{\Hom}{\mrm{Hom}} \nc{\id}{\mrm{id}} \nc{\im}{\mrm{im}}
\nc{\incl}{\mrm{incl}} \nc{\map}{\mrm{Map}} \nc{\mchar}{\rm char}
\nc{\nz}{\rm NZ} \nc{\supp}{\mathrm Supp}

\nc{\Alg}{\mathbf{Alg}}
\nc{\Bax}{\mathbf{Bax}}
\nc{\bff}{\mathbf f}
\nc{\bfk}{{\bf k}}
\nc{\bfone}{{\bf 1}}
\nc{\bfx}{\mathbf x}
\nc{\bfy}{\mathbf y}
\nc{\base}[1]{\bfone^{\otimes ({#1}+1)}} 
\nc{\Cat}{\mathbf{Cat}}

\nc{\detail}{\marginpar{\bf More detail}
    \noindent{\bf Need more detail!}
    \smallskip}
\nc{\Int}{\mathbf{Int}}
\nc{\Mon}{\mathbf{Mon}}
\nc{\remark}{\noindent{\bf Remark: }}
\nc{\remarks}{\noindent{\bf Remarks: }}
\nc{\Rings}{\mathbf{Rings}}
\nc{\Sets}{\mathbf{Sets}}

\nc{\BA}{{\Bbb A}} \nc{\CC}{{\Bbb C}} \nc{\DD}{{\Bbb D}}
\nc{\EE}{{\Bbb E}} \nc{\FF}{{\Bbb F}} \nc{\GG}{{\Bbb G}}
\nc{\HH}{{\Bbb H}} \nc{\LL}{{\Bbb L}} \nc{\NN}{{\Bbb N}}
\nc{\KK}{{\Bbb K}} \nc{\QQ}{{\Bbb Q}} \nc{\RR}{{\Bbb R}}
\nc{\TT}{{\Bbb T}} \nc{\VV}{{\Bbb V}} \nc{\ZZ}{{\Bbb Z}}

\nc{\cala}{{\mathcal A}} \nc{\calc}{{\mathcal C}}
\nc{\cald}{{\mathcal D}} \nc{\cale}{{\mathcal E}}
\nc{\calf}{{\mathcal F}} \nc{\calg}{{\mathcal G}}
\nc{\calh}{{\mathcal H}} \nc{\cali}{{\mathcal I}}
\nc{\call}{{\mathcal L}} \nc{\calm}{{\mathcal M}}
\nc{\caln}{{\mathcal N}} \nc{\calo}{{\mathcal O}}
\nc{\calp}{{\mathcal P}} \nc{\calr}{{\mathcal R}}
\nc{\cals}{{\mathcal S}}
\nc{\calt}{{\mathcal T}} \nc{\calw}{{\mathcal W}}
\nc{\calk}{{\mathcal K}} \nc{\calx}{{\mathcal X}}
\nc{\CA}{\mathcal{A}}

\nc{\fraka}{{\frak a}}
\nc{\frakA}{{\frak A}}
\nc{\frakb}{{\frak b}}
\nc{\frakB}{{\frak B}}
\nc{\frakH}{{\frak H}}
\nc{\frakM}{{\frak M}}
\nc{\bfrakM}{\overline{\frakM}}
\nc{\frakm}{{\frak m}}
\nc{\frakP}{{\frak P}}
\nc{\frakN}{{\mathfrak N}}
\nc{\frakp}{{\frak p}}
\nc{\frakS}{{\frak S}}

\font\cyr=wncyr10
\font\cyrs=wncyr7
\nc{\redt}[1]{\textcolor{red}{#1}} \nc{\li}[1]{\textcolor{red}{Li:
#1}} \nc{\zb}[1]{\textcolor{blue}{Bin: #1}}
\nc{\zhb}[1]{\textcolor{red}{Bin: #1}}
\begin{document}

\title[Renormalization of multiple zeta values]{Renormalization of multiple zeta values}
%
\author{Li Guo}
\address{Department of Mathematics and Computer Science,
         Rutgers University,
         Newark, NJ 07102}
\email{liguo@newark.rutgers.edu}
\author{Bin Zhang}
\address{Max Planck Institute for Mathematics,
Vivatsgasse 7, D-53111 Bonn, Germany; Yangtze Center of Mathematics,
Sichuan University, Chengdu, 610064, P. R. China }
\email{binzhang@mpim-bonn.mpg.de}

\maketitle

\begin{abstract}
Multiple zeta values (MZVs) in the usual sense are the special values of multiple variable zeta functions at positive integers. Their
extensive studies are important in both mathematics and
physics with broad connections and applications. In contrast, very
little is known about the special values of multiple zeta functions
at non-positive integers since the values are usually singular. We define and study
multiple zeta functions at integer values by adapting methods of
renormalization from quantum field theory, and following the Hopf algebra approach of Connes and Kreimer. This
definition of renormalized MZVs agrees with the convergent MZVs and extends the work of Ihara-Kaneko-Zagier on renormalization of MZVs with
positive arguments. We further show that the important quasi-shuffle (stuffle) relation for usual MZVs remains true for the renormalized MZVs.
\end{abstract}

\delete{
\begin{keyword}
multiple zeta values, algebraic Birkhoff decomposition, renormalization, quasi-shuffle, stuffle


\PACS 11M41, 16W30, 81T15

\end{keyword}

\end{frontmatter}
}

\tableofcontents


\section{Introduction}
Multiple zeta values (MZVs), as we know in the current literature, are defined to be the values of the multi-variable analytic function, called the {\bf multiple zeta function},
\begin{equation} \zeta(s_1,\cdots, s_k)=\sum_{n_1>\cdots>n_k>0}
    \frac{1}{n_1^{s_1}\cdots n_k^{s_k}}
\mlabel{eq:mzv}
\end{equation}
at positive integers $s_1,\cdots,s_k$ with $s_1>1$. With the earliest
study of MZVs went
back to Euler when $k=2$, their systematic study started in early
1990s with the works of Hoffman~\mcite{Ho0} and Zagier~\mcite{Za}.
Since then MZVs and their generalizations have been studied extensively by numerous authors from different point of views with connections to arithmetic geometry,
mathematical physics, quantum groups and knot
theory~\mcite{3BL,B-K,Ca1,Go3,G-M,Ho3,Kr,Te}.

In comparison, little is known about special values of multiple zeta functions at integers that are not all positive.
%
Through the recent work of Zhao~\mcite{Zh2} and Akiyama-Egami-Tanigawa~\mcite{AET} (see also~\mcite{Ma2}),
we know that $\zeta(s_1,\cdots,s_k)$ can be meromorphically continued to $\CC^k$
with singularities on the subvarieties
\begin{equation}
s_1=1; \
s_1+s_2=2,1,0,-2,-4, \cdots; {\rm\ and\ } \mlabel{eq:pole}
\sum_{i=1}^j
s_{i} \in \ZZ_{\leq j}\ (3\leq j\leq k).
\end{equation}
Thus $\zeta(s_1,\cdots,s_k)$ is undefined at most points with non-positive arguments.
%

In~\mcite{AET,A-K}, several definitions were proposed for the non-positive MZVs, that is, the values of $\zeta(s_1,\cdots,s_k)$ when $s_i$ are all non-positive. Some of them are
$$\lim_{r_1\to s_1}\cdots\lim_{r_k\to s_k} \zeta(r_1,\cdots,r_k),
\lim_{r_k\to s_k}\cdots\lim_{r_1\to s_1} \zeta(r_1,\cdots,r_k),
\lim_{r\to 0} \zeta(s_1+r,\cdots,s_k+r).$$
As expected they give different values.
Some good properties of the variously defined non-positive MZVs were obtained in the these papers.
But they fell short of the analogous properties of the usual MZVs, especially the double shuffle relations.


In this paper, we adapt a renormalization procedure (dimensional
regularization plus minimal subtraction) in quantum field theory
(QFT) to define the values of multiple zeta functions
$\zeta(s_1,\cdots,s_k)$ at $(s_1,\cdots,s_k)$ when $s_i$, $1\leq i\leq k$,
are all non-positive or all positive, that we expect to further extend
to when $s_i$ are {\em arbitrary integers}.
For our purpose, the dimensional regularization of
Feynman integrals is replaced by a regularization (or
deformation) of infinite series that has occurred in the study of
Todd classes for toric varieties~\mcite{B-Z}.

The renormalization procedure of QFT was put in the framework
of Hopf algebra and Rota-Baxter algebra by the recent works of Connes
and Kreimer~\mcite{C-K1,C-K2}, continued
in~\mcite{C-M,EGK2,EGK3}, and thus made possible for applications beyond QFT.
A fundamental result in this framework is the Algebraic Birkhoff Decomposition (Theorem~\mref{thm:diffBirk}). It states that for a given triple $(\calh, R, \phi)$ consisting of \begin{itemize}
\item a connected filtered Hopf algebra $\calh$,
\item a commutative Rota-Baxter algebra $R$ on which the Rota-Baxter operator $P:R\to R$ is idempotent, and
\item an algebra homomorphism $\phi: \calh\to R$,
\end{itemize}
there are unique algebra homomorphisms
$
\phi_-: \calh\to \CC+P(R)$ and
$\phi_+: \calh\to \CC+(\id-P)(R)$
such that
\begin{equation}
\phi=\phi_-^{\star (-1)}\star \phi_+.
\mlabel{eq:abd}
\end{equation}
Here $\star$ is the convolution product and $\phi_+$ is called the {\bf renormalization} of $\phi$. This algebraic setup is reviewed in Section~\mref{sec:setup} together with a discussion of quasi-shuffle algebras.

To apply this setup to the renormalization in QFT, one takes
\begin{itemize}
\item
$\calh=\calh_\FG$ to be the Connes-Kreimer Hopf algebra of Feynman diagrams, parameterizing regularized Feynman integrals,
\item
$R=\CC[[\vep,\vep^{-1}]$ to be the Rota-Baxter algebra of Laurent series, and
\item
$\phi$ to be the the regularized Feynman rule that assigns a Feynman diagram to the Laurent series expansion of the corresponding regularized Feynman integral.
\end{itemize}
Then the renormalized values of a Feynman integral is given by $\phi_+(\Gamma)$, where $\Gamma$
is the corresponding Feynman graph, when $\vep$ approaches zero.
For further details see~\mcite{C-K1,C-K2,C-M,EGK2,EGK3,Ma,F-G}.

To apply this setup to our study of renormalized MZVs, we similarly define
\begin{itemize}
\item
$\calh$ to be the quasi-shuffle Hopf algebra parameterizing regularized MZVs,
\item
$R$ to be the Rota-Baxter algebra $\CC[T][[\vep,\vep^{-1}]$ of log Laurent series, and
\item
$\phi$ to be the algebra homomorphism sending a symbol in $\calh$ to the Laurent series expansion of the corresponding regularized MZV.
\end{itemize}
Once these are obtained in Section~\mref{sec:gmzv}, the
Algebraic Birkhoff Decomposition in Eq.~(\mref{eq:abd}) applies to
give the renormalization $\phi_{+}$ from which the renormalized
MZVs can be derived when $\vep$ goes to zero, as in
the QFT case. However, there is an important difference from the QFT renormalization: in order to equip the regularized MZVs and the corresponding Hopf algebra with a suitable algebra structure that reflects the quasi-shuffle relation of the regularized MZVs, an extra parameter vector $\vec{r}$ has to be introduced in the regularized sums in addition to $\vep$. Thus the
renormalized MZVs at $\vec{s}$ from the Algebraic Birkhoff
Decomposition depend on $\vec{r}$, resulting in the {\bf
renormalized directional MZVs} $\zeta(\wvec{\vec{s}}{\vec{r}})$
in Definition~\mref{de:dmzv}.

This dependency on $\vec{r}$ is removed in the following
Section~\mref{sec:exp} in a consistent manner, giving the
{\bf renormalized MZVs} $\gzeta(\vec{s})$ in Definition~\mref{de:rmzv}.
Our main result Theorem~\mref{thm:main} shows that the renormalized MZVs satisfy the quasi-shuffle (or stuffle) relation, and  include as special cases the
MZVs defined either by convergence, by
analytic continuation, or by regularization in the sense of
Ihara-Kaneko-Zagier~\mcite{IKZ}. Parts of the proof are postponed to Section~\mref{sec:gzeta}
and Section~\mref{sss:nonpoqs}.
Here is the hierarchy of
MZVs introduced in this paper:
\begin{eqnarray*}
& \mbox{formal MZVs}\ \zeta(\vec{s}) \to \mbox{directional regularized MZVs}\ Z(\wvec{\vec{s}}{\vec{r}};\vep) &\\
& \to \mbox{renormalized directional MZVs}\ \zeta(\wvec{\vec{s}}{\vec{r}}) \to \mbox{renormalized MZVs}\ \gzeta(\vec{s}) &
\end{eqnarray*}
The concepts of regularization and renormalization have
already been introduced to the study of MZVs by Ihara-Kaneko-Zagier~\mcite{IKZ} to take care of the divergency of the MZVs $\zeta(s_1,\cdots,s_k)$ with $s_1=1$. As a part of their process, the natural algebra homomorphism from the the quasi-shuffle algebra for convergent MZVs to the algebra of convergent MZVs is extended to an algebra homomorphism from a larger quasi-shuffle algebra to an extension of the algebra of convergent MZVs. Thus they obtained their extended MZVs as an {\bf algebraic continuation} (we thank Robert Sczech for suggesting this term) in the sense that their extended MZVs preserves the quasi-shuffle relation.
{}From this point of view, we obtain our renormalized MZVs as an algebraic continuation that goes beyond theirs to cover the MZVs with all non-positive arguments. In a weak sense it covers arbitrary arguments (Definition~\mref{de:dmzv}).
More recently, Manchon and Paycha~\mcite{M-P1,M-P2} have considered renormalization of MZVs from the point of view of Chen integrals and Chen sums of symbols using a similar renormalization approach in
the spirit of Connes and Kreimer. The two approaches should be related though the exact link is still not clear.

This paper should lead to further studies of MZVs with arbitrary arguments.
First we can consider questions related to the renormalization procedure,
such as the renormalization of MZVs with arbitrary arguments, the
dependence of renormalized MZVs on the regularization and
renormalization. We would also like to study the extension of the double shuffle
relation to renormalized MZVs, and the possible connection to
rational associators in the sense of Drinfel'd~\mcite{Dr} and DMR
in the sense of Racinet~\mcite{Ra}.
Possible arithmetic properties of these renormalized
MZVs, such as the Kummer type congruences, are also interesting to investigate.
Some of these directions will be pursued in future works.

\medskip

\noindent {\bf Acknowledgements.} Both authors thank the Max-Planck
Institute for Mathematics in Bonn for the stimulating environment
where this project was started, and thank Matilde Marcolli for
encouragement. We also thank Dominique Manchon and Sylvie Paycha for communications on
their papers~\mcite{M-P1,M-P2} and for their comments on our paper.
The first named author thanks NSF for
support, and is indebted to Kurusch Ebrahimi-Fard and Dirk Kreimer
for their collaborations on QFT that inspired the
renormalization approach in this paper. Thanks also go to Herbert Gangl, Robert Sczech and
Jianqiang Zhao for discussions and comments.

\section{The algebraic setup}
\mlabel{sec:setup}
We describe the general setup for our later applications to renormalization of MZVs.
In the following an algebra means a $\bfk$-algebra where $\bfk$ is a unitary commutative ring that we usually take to be $\CC$. Denote the unit of $\bfk$ by $\bfone$.

\subsection{The algebraic Birkhoff decomposition}
We review the algebraic framework of Connes and Kreimer for renormalization of perturbative quantum field theory.

A {\bf connected filtered Hopf algebra} is a Hopf
algebra $(H, \Delta)$ with $\bfk$-submodules $H^{(n)},\ n\geq 0$ of $H$ such
that
\begin{eqnarray*}
&H^{(n)}\subseteq H^{(n+1)}, \quad
\cup_{n\geq 0} H^{(n)} = H, \quad
H^{(p)} H^{(q)}\subseteq H^{(p+q)}, &\\
&
\Delta(H^{(n)}) \subseteq \sum_{p+q=n} H^{(p)}\otimes H^{(q)}, \quad
H^{(0)}=\bfk {\rm \ (connectedness)}.&
\end{eqnarray*}
%
%
Let $\lambda\in \bfk$. A {\bf Rota--Baxter algebra} of weight $\lambda$ is a pair $(R,P)$ where $R$ is a unitary $\bfk$-algebra and $P:R\to R$ is a linear operator such that
\begin{equation}
P(x)P(y)=P(xP(y))+P(P(x)y)+\lambda P(xy),
\mlabel{eq:rbe}
\end{equation}
for any $x,\,y\in R$. Often $\theta=-\lambda$ is used, especially in the physics literature.
It follows from the definition that $P(R)$ and $(-\lambda-P)(R)$ are non-unitary subalgebras
of $R$. So $\bfk+P(R)$ and $\bfk +(-\lambda-P)(R)$ are unitary subalgebras.

\begin{theorem}
Let $H$ be a commutative connected filtered Hopf algebra.
Let $(R,P)$ be a Rota-Baxter algebra of weight $-1$.
Let $\phi: H \to R$ be an algebra homomorphism.
\begin{enumerate}
\item
There are algebra homomorphisms $\phi_-: H \to \bfk+P(R)$ and
$\phi_+: H \to \bfk+(1-P)(R)$ with the decomposition
\begin{equation}
\phi=\phi_-^{\star\, (-1)}\star \phi_+,
\mlabel{eq:phidecom}
\end{equation}
called the {\bf Algebraic Birkhoff Decomposition} of $\phi$.
Here $\star$ is the convolution product and $\phi_-^{\star\, (-1)}$ is the inverse of $\phi_-$ with respect to $\star$.
Further,
\begin{equation} \phi_-(x)=-P\big(\phi(x)+\sum_{(x)} \phi_-(x')\phi(x'')\big)
\mlabel{eq:phi-}
\end{equation}
and
\begin{equation}
 \phi_+(x)=(\id-P)\big(\phi(x)+\sum_{(x)} \phi_-(x')\phi(x'')\big).
\mlabel{eq:phi+}
\end{equation}
Here we have used the notation
$\Delta(x)=x\ot 1 + 1\ot x + \sum_{(x)} x'\ot x''.$
\mlabel{it:decom}
\item
If $P^2=P$, then the decomposition in Eq.~(\mref{eq:phidecom}) is unique.
\mlabel{it:uni}
\end{enumerate}
\mlabel{thm:diffBirk}
\end{theorem}
\proofbegin
For item~(\mref{it:decom}), see~\cite{C-K1} and \cite[Theorem II.5.1]{Ma}.
For item~(\mref{it:uni}), see~\cite[Theorem 3.7]{EGK3} where one can also find a proof of item~(\mref{it:decom}) using Rota-Baxter algebras. \proofend

\subsection{Quasi-shuffle algebras}
Let $M$ be a commutative semigroup. For each integer $k\geq 0$, let $\bfk M^k$ be the free $\bfk$-module with basis $M^k$, with the convention that $M^0=\{\bfone\}$. Let
\begin{equation}
 \calh_M=\bigcup_{k=0}^\infty \bfk\, M^k.
 \mlabel{eq:qprod}
\end{equation}
Following~\mcite{Ho2}, define the {\bf quasi-shuffle product} $\msh$ by
first taking $\bfone$ to be the multiplication identity.
Next for any $m,n\geq 1$ and
$\vec a:=(a_1,\cdots , a_m)\in M^{m}$ and $\vec b:=(b_1, \cdots, b_n)\in M^{n}$, denote
$\vec a\,'=(a_2,\cdots,a_m)$ and $\vec b\,'=(b_2,\cdots,b_n)$. Recursively define
\begin{equation}
\vec{a} \msh \vec{b} =\big(a_1,   \vec{a}\,'\msh \vec b\big )
  + \big(b_1, \vec a \msh \vec{b}\,' \big) \\
  + \big(a_1 b_1,   \vec{a}\,' \msh \vec{b}\,'\big)
\mlabel{eq:qshuf}
\end{equation}
with the convention that $\vec{a}\,'=\bfone$ if $m=1$, $\vec{b}\,'=\bfone$ if $n=1$ and
$(a_1b_1,\vec{a}\,'\msh \vec{b}\,')=(a_1b_1)$ if $m=n=1$.

Quasi-shuffle is also known as harmonic product~\mcite{Ho1} and coincides with the stuffle product~\mcite{3BL,Br}
in the study of MZVs. Variations of the stuffle product have also appeared in ~\mcite{Ca2,Eh}.
See \S~\mref{sss:nonpoqs} for further details.
It is shown~\mcite{E-G1} to be the same as the mixable shuffle product~\mcite{G-K1,G-K2}
which is also called overlapping shuffles~\mcite{Ha} and generalized shuffles~\mcite{Go3}, and can be interpreted in terms of Delannoy paths~\mcite{A-H,Fa,Lo}.

The following theorem is a simple generalization of~\cite[Theorem 2.1, 3.1]{Ho2} where $M$ has the extra condition of being a locally finite set to ensure the grading structure on $\calh_M$.

\begin{theorem}
Let $M$ be a commutative semigroup.
Equip $\calh_M$ with
the submodules $\calh_M^{(n)}=\oplus_{i=0}^n \bfk\, M^i$,
The quasi-shuffle product $\msh$,
the deconcatenation coproduct
\begin{eqnarray}
&&\Delta: \calh_M \to \calh_M \barot \calh_M,\\
 \Delta( a_1, \cdots ,  a_k)&=&1\barot (a_1,  \cdots ,  a_k)
+ \sum_{i=1}^{k-1} (a_1,  \cdots,  a_i)\barot (a_{i+1},  \cdots ,  a_k) \notag\\
&&+ (a_1,  \cdots ,  a_k) \barot 1
\mlabel{eq:coprod}
\end{eqnarray}
and the projection counit
$\vep: \calh_M \to \bfk$
onto the direct summand $\bfk \subseteq \calh_M$.
Then $\calh_M$ is a commutative connected filtered Hopf algebra.
\mlabel{thm:hopf}
\end{theorem}

\proofbegin
By the same proofs as \cite[Theorem 2.1]{Ho2} and \cite[Theorem 3.1]{Ho2}, $\calh_M$ is
a bialgebra.
By the definition of $\msh$ and $\Delta$, $\calh_M$ is connected filtered with the submodules $\calh_M^{(n)}, n\geq 0$. Then $\calh_M$ is automatically a Hopf algebra by ~\cite[Proposition 5.3]{F-G}, for example.
\proofend

We prove the following property for later applications.
\begin{mprop}
For $k\geq 1$, let $\Sigma_k$ be the permutation group on $\{1,\cdots,k\}$.
For $\vec{a}=(a_1,\cdots,a_k)\in M^k \subseteq \calh_M$, define
$\sigma(\vec{a})=(a_{\sigma(1)},\cdots,a_{\sigma(k)})$ and define
$\vec{a}^{(\Sigma_k)}=\sum_{\sigma\in \Sigma_k} \sigma(\vec{a}).$
Then for $a_{k+1}\in M$, we have
\begin{equation}
\vec{a}^{(\Sigma_k)} \msh (a_{k+1}) = (a_1,\cdots,a_k,a_{k+1})^{(\Sigma_{k+1})} + \sum_{i=1}^k (a_1,\cdots,a_ia_{k+1},\cdots,a_{k})^{(\Sigma_k)}
\mlabel{eq:permshuf}
\end{equation}
where in the sum,
$\sigma(a_1,\cdots,a_ia_{k+1},\cdots,a_k)=(a_{\sigma(1)},\cdots,a_{\sigma(i)}a_{k+1},\cdots,a_{\sigma(k)}).$
\mlabel{pp:prep}
\end{mprop}
\proofbegin By the quasi-shuffle relation in
Eq.~(\mref{eq:qshuf}), we have
\begin{eqnarray*}
\lefteqn{\vec{a}^{(\Sigma_k)} \msh (a_{k+1})= \sum_{\sigma\in
\Sigma_k}
\Big((a_{\sigma(1)},\cdots,a_{\sigma(k)},a_{k+1})}\\
&& + \sum_{i=1}^k (a_{\sigma(1)},\cdots,a_{\sigma(i)},a_{k+1},a_{\sigma(i+1)},\cdots,a_{k})
 + \sum_{i=1}^k (a_{\sigma(1)},\cdots,a_{\sigma(i)}a_{k+1},\cdots,a_{\sigma(k)})\Big),
\end{eqnarray*}
hence the proposition. It also follows from the Partition Identity of Hoffman~\mcite{Ho0} whose proof
only needs the quasi-shuffle relation~\mcite{Br1},
or the Bohnenblust-Spitzer formula for Rota-Baxter algebras~\mcite{E-G5,Ro3}.
\proofend

\section{Renormalized directional multiple zeta values}
\mlabel{sec:gmzv}
We now introduce directional regularized MZVs and the corresponding Hopf algebra. We then show that the directional regularized MZVs have Laurent series expansion with log coefficients, giving an algebra homomorphism from the Hopf algebra to Laurent series. This allows us to apply the
algebraic Birkhoff decomposition in Theorem~\mref{thm:hopf} to obtain
renormalized directional MZVs.
\subsection{The Hopf algebra of directional regularized multiple zeta values}
We consider the commutative semigroup
\begin{equation}
\frakM= \{{{\wvec{s}{r}}}\ \big|\ (s,r)\in \ZZ \times \RR_{>0}\}
\mlabel{eq:mbase}
\end{equation}
with the multiplication
$ {\wvec{s}{r}} {\wvec{s'}{r'}}={\wvec{s+s'}{r+r'}}.$
%
By Theorem~\mref{thm:hopf},
$$\calh_{\frakM}:=\sum_{k\geq 0} \CC\, \frakM^k,$$
with the quasi-shuffle product $\msh$ and the deconcatenation coproduct $\Delta$, is a
connected filtered Hopf algebra.
The same is true with the sub-semigroup
$$\frakM^-=\{{\wvec{s}{r}}\ |\ (s,r)\in \ZZ_{\leq 0} \times \RR_{>0}\}.$$
For $w_i=\wvec{s_i}{r_i}\in \frakM,\ i=1,\cdots,k$, we use the notations
$$ \vec{w}=(w_1,\cdots,w_k)
=\wvec{s_1,\cdots,s_n}{r_1,\cdots,r_k}=\wvec{\vec{s}}{\vec{r}},\
{\rm where\ } \vec{s}=(s_1,\cdots,s_k), \vec{r}=(r_1,\cdots,r_k).$$
For
$\vep\in \CC$ with ${\rm
Re}(\vep)<0$, define the {\bf directional regularized MZV}:
\begin{equation}
Z(\wvec{\vec{s}}{\vec{r}};\vep)=\sum_{n_1>\cdots>n_k>0}
\frac{e^{n_1\,r_1\vep} \cdots
    e^{n_k\,r_k\vep}}{n_1^{s_1}\cdots n_k^{s_k}}
\mlabel{eq:reggmzv}
\end{equation}
It converges for any $\wvec{\vec{s}}{\vec{r}}$ and is regarded as the regularization of the {\bf formal MZV}
\begin{equation}
\zeta (\vec{s})= \sum_{n_1>\cdots>n_k>0} \frac{1}{n_1^{s_1} \cdots
    n_k^{s_k}}
\mlabel{eq:formgmzv}
\end{equation}
which converges only when $s_i>0$ and $s_1>1$.
It is related to the multiple polylogarithm
$${\rm Li}_{s_1,\cdots,s_k}(z_1,\cdots,z_k)=\sum_{n_1>\cdots n_k>0}
    \frac{z_1^{n_1} \cdots z_k^{n_k}}{n_1^{s_1}\cdots n_k^{s_k}}$$
by a change of variables $z_i=e^{r_i\vep}, 1\leq i\leq k$.
As is well-known~\mcite{3BL,Go3}, the product of multiple polylogarithms as functions satisfies the quasi-shuffle (stuffle) relation of the nested sums. Therefore the product of regularized MZVs as functions also satisfies the quasi-shuffle relation: if
 $ \wvec{\vec{s}}{\vec{r}}\msh \wvec{\vec{s}\,'}{\vec{r}\,'}
=\sum \wvec{\vec{s}\,''}{\vec{r}\,''}$, then
\begin{equation}
Z(\wvec{\vec{s}}{\vec{r}};\vep)Z(\wvec{\vec{s}\,'}{\vec{r}\,'};\vep)
= Z(\wvec{\vec{s}}{\vec{r}}\msh \wvec{\vec{s}\,'}{\vec{r}\,'};\vep)
:= \sum Z(\wvec{\vec{s}\,''}{\vec{r}\,''};\vep).
\mlabel{eq:qsh1}
\end{equation}
 We thus obtained an algebra homomorphism
\begin{equation}
Z: \calh_\frakM \to \sum_{{\wvec{\vec{s}}{\vec{r}}\in \cup_{n\geq 0} \frakM^n}}\
    \CC\, Z(\wvec{\vec{s}}{\vec{r}};\vep),
    \wvec{\vec{s}}{\vec{r}} \mapsto Z(\wvec{\vec{s}}{\vec{r}};\vep).
\mlabel{eq:regz}
\end{equation}
With this map, $\calh_\frakM$ is a parametrization of the directional regularized MZVs that also reflects their multiplication property.

\subsection{Log Laurent series of directional regularized multiple zeta values}
We first construct Laurent series with log coefficients. We then show that the nested sums from directional regularized MZVs are such log Laurent series.

Let $\CC\{\{\vep,\vep^{-1}\}$ be the algebra of convergent
Laurent series, regarded as a subalgebra of the algebra of (germs
of) complex valued functions meromorphic in a neighborhood of
$\vep=0$. Take $\ln \vep$ to be analytic on $\CC\backslash (-\infty,0]$.

\begin {lemma} $\ln (-\vep) $ is transcendental over
$\CC\{\{\vep,\vep^{-1}\}$.
\mlabel{lem:tran}
\end {lemma}

\proofbegin We give a simple proof for the lack of references.
Assume $\ln (-\vep) $ is algebraic over the field
$\CC\{\{\vep,\vep^{-1}\}$ with the monic minimal polynomial
$$\ln ^n(-\vep)+a_{n-1}(\vep)\ln ^{n-1}(-\vep)+\cdots +a_0(\vep)=0.$$
Differentiating the above equation, we have
$$\sum _{i=0}^n \big(a'_i(\vep )\ln ^i(-\vep)+\frac i{\vep}a_i(\vep)\ln ^{i-1}(-\vep )\big)=0.$$
The highest power term in $\ln (-\vep)$ is $(\frac
n{\vep}+a'_{n-1}(\vep ))\ln ^{n-1}(-\vep)$. Because of the minimality, $n/\vep+a'_{n-1}(\vep )$ has to be 0, which is impossible
for $a_{n-1}(\vep )\in \CC\{\{\vep,\vep^{-1}\}$.
\proofend

\begin{lemma}
$\CC\{\{\vep,\vep^{-1}\} [\ln (-\vep)]$ is closed under the differential operator
$d/d\vep$. It is also closed under the indefinite
integral operator: the antiderivatives of any $f\in \CC\{\{\vep,\vep^{-1}\} [\ln (-\vep)]$ are in
$\CC\{\{\vep,\vep^{-1}\} [\ln (-\vep)]$. \mlabel{lem:diffint}
\end{lemma}
\proofbegin
Let $f(\vep)\in \CC\{\{\vep,\vep^{-1}\} [\ln (-\vep)]$. Then
$$f(\vep)=\sum_{n=0}^M a_n(\vep)\ln^n (- \vep)
=\sum_{n=0}^M \big(\sum_{k\geq N_n} a_{n,k}\vep^k \ln^n(-\vep)\big)$$
with $\sum_{k\geq N_n} a_{n,k}\vep^k \in \CC\{\{\vep,\vep^{-1}\}$.
For each $0\leq n\leq M$, the series for the inside sum converges absolutely and uniformly in a nonempty
open interval of $\{\vep\in \CC\,\big |\, -\infty<\vep <0\}$. Thus the
series can be differentiated and integrated term by term. Thus we
only need to show that  the derivative and anti-derivatives of
$\vep^k \ln ^n (-\vep)$, $k\in \ZZ, n\in \ZZ_{\geq 0}$ are linear
combinations of functions of the same form. This is easy  to check for
derivatives.

For anti-derivatives, we use induction on $n$. It is clear when $n=0$.
The induction step follows from the integration by parts formula
$$ \int \vep^k \ln ^n (-\vep) d\vep = \frac{1}{k+1}\vep^{k+1} \ln ^n (-\vep)
+ \frac{n}{k+1} \int  \vep^{k} \ln ^{n-1} (-\vep) d\vep$$
when $k\neq -1$ and
$ \int \frac{\ln ^n (-\vep)}{\vep}d\vep = \frac{1}{n+1} \ln ^{n+1} (-\vep)+C.$
\proofend

Because of Lemma~\mref{lem:tran}, we have
\begin{equation}
\CC\{\{\vep,\vep^{-1}\} [\ln (-\vep)]\cong
\CC\{\{\vep,\vep^{-1}\} [T]
\hookrightarrow \CC[[\vep,\vep^{-1}][T]
\mlabel{eq:lau2}
\end{equation}
sending $-\ln (-\vep)$ to $T$.
Here $\CC[[\vep,\vep^{-1}][T]$ denotes the polynomial algebra over the formal
Laurent series $\CC[[\vep,\vep^{-1}]$.

An element of $\CC[[\vep,\vep^{-1}] [T]$ is of the form
$\sum_{n=0}^M a_n(\vep)T^n$ with
$$a_n(\vep)=\sum_{k\geq N_n} a_{n,k}\vep^k \in \CC[[\vep,\vep^{-1}], 0\leq n\leq M.$$
Taking $N=\min_{0\leq n\leq M} N_n$ and letting $a_{n,k}=0$ for $N\leq k< N_n$, we have
$$\sum_{n=0}^M a_n(\vep)T^n = \sum_{n=0}^M \big(\sum_{k\geq N} a_{n,k}\vep^k\big) T^n =\sum_{k\geq
N} \big(\sum_{n=0}^M a_{n,k} T^n\big) \vep^k.$$
This gives an element of the {\bf algebra of log Laurent series}
$\CC[T][[\vep,\vep^{-1}]$ with coefficients in
$\CC[T]$.
Combining with Eq.~(\mref{eq:lau2}), we obtain a natural algebra injection
\begin{equation}
u: \CC\{\{\vep,\vep^{-1}\} [\ln (-\vep)]  \to \CC[T][[\vep,\vep^{-1}]
\mlabel{eq:formulau}
\end{equation}
with which we identify $\CC\{\{\vep,\vep^{-1}\} [\ln (-\vep)]$ as a subalgebra of $\CC[T][[\vep,\vep^{-1}]$.


\begin{theorem}
For any $\vec s\in \ZZ ^k$, $\vec r \in \Bbb {Z}_{>0}^k$,
$Z(\wvec{\vec s}{\vec r};\vep)$ is in $\CC \{\{\vep,\vep^{-1}\}[ln(-\vep)]$ and can thus be regarded as
an element in $\CC[T][[\vep,\vep^{-1}]$ by Eq.~(\mref{eq:formulau}).
If $\vec s$ is in $\Bbb {Z}_{\le 0}^k$,
then $Z(\wvec{\vec s}{\vec r};\vep)$ is in $\CC\{\{\vep,\vep^{-1}\}$.
\mlabel{thm:zlaurent}
\end{theorem}
\proofbegin
First notice that, for $r,i\in \ZZ_{>0}$,
\begin{equation}
 \sum _{n\geq i}e^{nr\vep }=\frac {1}{1-e^{r\vep}}\,e^{ir\vep}
 \mlabel{eq:geo}
 \end{equation}
has a Laurent series expansion at $\vep=0$. Since $Z(\wvec{s}{r};\vep)$ is
uniformly convergent on compact subsets in $Re(\vep)<0$, by repeatedly differentiating Eq.~(\mref{eq:geo}), we see that, for $s\in \ZZ_{<0}$,
\begin{eqnarray}
\sum _{n\geq i}n^{-s}e^{nr\vep}
&=& r^s\Big( \sum_{p=0}^{-s} \binc{-s}{p} \big(\frac{1}{1-e^{r\vep}}\big)^{(p)}(e^{ir\vep})^{(-s-p)}\Big) \notag\\
&=&\sum _{p=0}^{-s}\binc{-s}{p}
\big(\frac{1}{1-e^{r\vep}}\big)^{(p)}
\frac{e^{ir\vep}}{r^{p}\,i^{s+p}}
\mlabel{eq:diffsum}
\end{eqnarray}
has a Laurent series expansion at $\vep=0$.

Now we prove by induction on $k\geq 1$.
Let $k=1$. Then $\vec s=s\in \ZZ$. The case when $s\leq 0$ follows
from Eq.~(\mref{eq:diffsum}) with $i=1$. When $s>0$, we note that
$Z'(\wvec{s}{r};\vep)=rZ(\wvec{s-1}{r};\vep)$. By Eq.~(\mref{eq:geo}), $Z(\wvec{0}{r};\vep)$ is in $\CC\{\{\vep,\vep^{-1}\}[\ln(-\vep)]$
which is closed under integration (Lemma~\mref{lem:diffint}). Thus $Z(\wvec{1}{r};\vep)$ is in $\CC\{\{\vep,\vep^{-1}\}[\ln(-\vep)]$ and,
by an induction on $s$, the same holds for $Z(\wvec{s}{r};\vep)$ for any $s>0$.

Assume that the statements hold for $k\geq 1$ and prove for $Z(\wvec{\vec s}{\vec r}; \vep)$ with $\vec s=(s_1,\cdots,s_{k+1})$
in two cases.

\noindent
{\bf Case 1. } Suppose $s_i\leq 0$ for some $1\leq i\leq k+1$.
Then for fixed $n_{i-1}>n_{i+2}+1>0$,
\begin{eqnarray*}
\lefteqn{
\sum_{n_{i-1}>n_i>n_{i+1}>n_{i+2}}
\frac{e^{n_ir_i\vep+n_{i+1}r_{i+1}\vep}}{n_i^{s_i}n_{i+1}^{s_{i+1}}}}\\
&=& \sum_{n_{i-1}>n_i\geq n_{i+1}>n_{i+2}}
\frac{e^{n_ir_i\vep+n_{i+1}r_{i+1}\vep}}{n_i^{s_i}n_{i+1}^{s_{i+1}}}
-\sum_{n_{i-1}> n_i=n_{i+1}>n_{i+2}}
\frac{e^{n_ir_i\vep+n_{i+1}r_{i+1}\vep}}{n_i^{s_i}n_{i+1}^{s_{i+1}}}\\
&=& \sum_{n_{i-1}>n_{i+1}>n_{i+2}}
\frac{e^{n_{i+1}r_{i+1}\vep}}{n_{i+1}^{s_{i+1}}}
\sum_{n_i=n_{i+1}}^{n_{i-1}-1}\frac{e^{n_ir_i\vep}}{n_i^{s_i}}
-\sum_{n_{i-1}> n_{i+1}>n_{i+2}}
\frac{e^{n_{i+1}(r_i+r_{i+1})\vep}}{n_{i+1}^{s_i+s_{i+1}}}.
\end{eqnarray*}
Applying Eq.~(\mref{eq:diffsum}) to the inner sum of the first term, we have
\allowdisplaybreaks{
\begin{eqnarray*}
\sum _{n_i=n_{i+1}}^{n_{i-1}-1}\frac{e^{n_ir_i\vep}}{n_i^{s_i}}&=&
\sum _{n_i=n_{i+1}}^\infty \frac{e^{n_ir_i\vep}}{n_i^{s_i}}
-\sum _{n_i=n_{i-1}}^\infty \frac{e^{n_ir_i\vep}}{n_i^{s_i}}\\
&=&
\sum _{p=0}^{-s_i} \frac{1}{r_i^p} \binc{-s_i}{p}
\Big(\frac{1}{1-e^{r_i\vep}}\Big)^{(p)}
\Big ( \frac{e^{n_{i+1}r_i\vep}}{n_{i+1}^{s_i+p}}
-\frac{e^{n_{i-1}r_i\vep}}{n_{i-1}^{s_i+p}}\Big ).
\end{eqnarray*}
}
Thus
\allowdisplaybreaks
{
\begin{eqnarray*}
\lefteqn{
\sum_{n_{i-1}>n_i>n_{i+1}>n_{i+2}}
\frac{e^{n_ir_i\vep+n_{i+1}r_{i+1}\vep}}{n_i^{s_i}n_{i+1}^{s_{i+1}}}}\\
&=&
\sum _{p=0}^{-s_i} \frac{1}{r_i^p} \binc{-s_i}{p}
\Big(\frac{1}{1-e^{r_i\vep}}\Big)^{(p)}
\Big ( \sum_{n_{i-1}>n_{i+1}>n_{i+2}}
\Big ( \frac{e^{n_{i+1}(r_i+r_{i+1})\vep}}{n_{i+1}^{s_i+s_{i+1}+p}}
-\frac{e^{n_{i+1}r_{i+1}\vep}}{n_{i+1}^{s_{i+1}}}
\frac{e^{n_{i-1}r_{i}\vep}}{n_{i-1}^{s_{i}+p}}
\Big ) \Big)\\
&& -\sum_{n_{i-1}> n_{i+1}>n_{i+2}}
\frac{e^{n_{i+1}(r_i+r_{i+1})\vep}}{n_{i+1}^{s_i+s_{i+1}}}.
\end{eqnarray*}
}
Then we have
\allowdisplaybreaks{
\begin{eqnarray*}
\lefteqn{Z(\wvec{\vec s}{\vec r};\vep)
=\sum_{n_{1}>\cdots
>n_{i-1}} \frac {e^{n_1r_1\vep+\cdots +n_{i-1}r_{i-1}\vep}}{n_1^{s_1}\cdots
n_{i-1}^{s_{i-1}}}
\sum_{n_{i-1}>n_{i}>n_{i+1}>n_{i+2}}\frac{e^{n_ir_i\vep}}{n_i^{s_i}}
\frac{e^{n_{i+1}r_{i+1}\vep}}{n_{i+1}^{s_{i+1}}}}\\
&&
\times \sum_{n_{i+1}>n_{i+2}>\cdots >n_{k+1}>0}
\frac {e^{n_{i+2}r_{i+2}\vep+\cdots +n_{k+1}r_{k+1}\vep}}
{n_{i+2}^{s_{i+2}}\cdots n_{k+1}^{s_{k+1}}}\\
&=&
\sum _{p=0}^{-s_i} \frac{1}{r_i^p} \binc{-s_i}{p}
\left(\frac{1}{1-e^{r_i\vep}}\right)^{(p)}
\Big(Z\big(\wvec{s_1,\cdots,s_{i-1},s_i+s_{i+1}+p,s_{i+2},\cdots,s_{k+1}}
    {r_1,\cdots,r_{i-1},r_i+r_{i+1},r_{i+2},\cdots,r_{k+1}};\vep \big)\\
&&
 -Z\big(\wvec{s_1,\cdots,s_{i-2},s_{i-1}+s_i+p,s_{i+1},\cdots,s_{k+1}}
{r_1,\cdots,r_{i-2},r_{i-1}+r_i,r_{i+1},\cdots,r_{k+1}};\vep\big)\Big)\\
&&-Z\big(\wvec{s_1,\cdots,s_{i-1},s_i+s_{i+1},s_{i+2},\cdots,s_{k+1}}
{r_1,\cdots,r_{i-1},r_i+r_{i+1},r_{i+2},\cdots,r_{k+1}};\vep\big).
\end{eqnarray*}
}
The induction hypothesis applies to each term on the right hand side, completing the induction on $k$ in this case. In particular this completes the induction when $\vec{s}\in \ZZ^{k+1}_{\leq 0}$.

\noindent {\bf Case 2.} Suppose $s_i>0$ for all $1\leq i\leq k+1$.
We use induction on the sum $s:=\sum_{i=1}^{k+1} s_i$. Then $s\geq
k+1$. If $s=k+1$, then $s_i=1$ for $1\leq i\leq k+1$. Note that
\begin{equation}
\frac {d}{d\vep } Z(\wvec{\vec s}{\vec r};\vep)=\sum r_i Z(\wvec{\vec s-\vec e _i}{\vec r};\vep),
\mlabel{eq:zdiff}
\end{equation}
where $\vec e_i, 1\leq i\leq k+1$ is the $i$-th unit vector in
$\ZZ ^{k+1}$. Each term on the right hand
side is in $\CC\{\{\vep,\vep^{-1}\}[\ln(-\vep)]$ by Case 1.
So by Lemma~\mref{lem:diffint}, $Z(\wvec{\vec s}{\vec r};\vep)$ is in $\CC[T][[\vep,\vep^{-1}]$. The
inductive step follows from Eq.~(\mref{eq:zdiff}) and the induction assumption.
\proofend

\subsection{Renormalized directional MZVs}

Combining Eq.~(\mref{eq:regz}), Theorem~\mref{thm:zlaurent}
and Eq.~(\mref{eq:formulau}), we obtain an algebra homomorphism
\begin{equation}
\uni{Z}: \calh_\frakM \to \CC\{\{\vep,\vep^{-1}\}[\log (-\vep)] \ola{u} \CC[T][[\vep,\vep^{-1}],\ \wvec{\vec{s}}{\vec{r}} \mapsto
    u\big(Z(\wvec{\vec{s}}{\vec{r}};\vep)\big).
\mlabel{eq:zmap}
\end{equation}
In the same
way, $\uni{Z}$ restricts to an algebra homomorphism
$$
\uni{Z}: \calh_{\frakM^-} \to \CC[[\vep,\vep^{-1}].
$$
For any commutative $\bfk$-algebra $K$, $K[[\vep,\vep^{-1}]$
is a Rota-Baxter algebra of weight -1 with the Rota-Baxter operator $P$ to be the projection to $\vep^{-1} K [\vep^{-1}]$:
\begin{equation}
 P\big ( \sum_{n\geq N} \alpha_k \vep^k\big )=\sum_{k\leq -1}\alpha_k \vep^k.
 \mlabel{eq:prb}
\end{equation}
This can be directly verified as with the well-known case of $\CC[[\vep,\vep^{-1}]$ in~\mcite{C-K1}.

Thus we can apply the Algebraic Birkhoff Decomposition in Theorem~\mref{thm:diffBirk} and obtain

\begin{coro}
We have
$$ \uni{Z}=\uni{Z}_-^{-1} \star \uni{Z}_+$$
and the map $\uni{Z}_+: \calh_\frakM\to \CC[T][[\vep]]$ is an algebra homomorphism which restricts to an algebra homomorphism
$\uni{Z}_+: \calh_{\frakM^-}\to \CC[[\vep]]$.
\mlabel{co:renormz}
\end{coro}

Because of Corollary~\mref{co:renormz}, the following definition is valid.
\begin{mdefn}
For $\vec{s}=(s_1,\cdots,s_k)\in \ZZ^k$ and $\vec{r}=(r_1,\cdots,r_k)\in \RR_{>0}^k$, define
the {\bf renormalized directional MZV} by
\begin{equation} \zeta\lp \wvec{\vec{s}}{\vec{r}}\rp = \lim_{\vep\to 0} \uni{Z}_+\lp \wvec{\vec{s}}{\vec{r}};\vep\rp .
\mlabel{eq:dmzv}
\end{equation}
Here $\vec{r}$ is called the {\bf direction vector}.
\mlabel{de:dmzv}
\end{mdefn}

As a consequence of Corollary~\mref{co:renormz}, we have
\begin{coro}
The renormalized directional MZVs satisfies the quasi-shuffle relation
\begin{equation}
\zeta\lp \wvec{\vec{s}}{\vec{r}}\rp \zeta\lp \wvec{\vec{s}\,'}{\vec{r}\,'}\rp
= \zeta\lp \wvec{\vec{s}}{\vec{r}} \msh \wvec{\vec{s}\,'}{\vec{r}\,'} \rp.
\mlabel{eq:dqshuf}
\end{equation}
Here the right hand side is defined in the same way as in Eq.~(\mref{eq:qsh1}).
\mlabel{co:qshuf}
\end{coro}

We next give an
explicit formula for the renormalized directional \mzvs.
\begin{mdefn} {\rm
Let $\Pi_k$ be the set of ordered partitions (compositions) of $k$,
consisting of ordered sequences $(i_1,\cdots,i_p)$ such that
$i_1+\cdots+i_p=k$. For $1\leq j\leq p$, define the partial sum
$I_j=i_1+\cdots+i_{j}$ with the convention that $I_0=0$.
The {\bf partition vectors} of $\vec{s}\in \RR^k$ from the ordered partition $(i_1,\cdots,i_p)$ are the vector $\vec{s}^{(j)}:=(s_{I_{j-1}+1},\cdots,s_{I_{j}})$, $1\leq j\leq p$.
}
\mlabel{de:part}
\end{mdefn}
\begin{theorem}
Let $P: \CC[T][[\vep,\vep^{-1}]\to \CC[T][\vep^{-1}]$ be the Rota--Baxter operator in Eq.~(\mref{eq:prb}). Denote $\check{P}=-P$ and $\tilde{P}=\id-P$. For $\vec{s}=(s_1,\cdots,s_k)\in \ZZ^k$ and
$\vec{r}=(r_1,\cdots,r_k)\in \NN_{>0}^k$,
%
\begin{eqnarray}
\lefteqn{{\ }\hspace{-.7cm} \uni{Z}_-\lp \wvec{\vec s}{\vec r};\vep\rp
= \hspace{-.5cm}
\sum_{(i_1,\cdots,i_p)\in \Pi_k}
\check{P}\llp
    \llp \cdots \check{P}\llp\check{P}\llp \uni{Z}\lp \wvec{\vec s ^{(1)}}{\vec r ^{(1)}};\vep\rp\rrp\,
    \uni{Z}\lp\wvec{\vec s^{(2)}}{\vec r ^{(2)}};\vep\rp \rrp \cdots\rrp
    \uni{Z}\lp\wvec{\vec s^{(p)}}{\vec r ^{(p)}};\vep\rp \rrp. } \notag \\
\lefteqn{{\ }\hspace{-.7cm}  \uni{Z}_+\lp \wvec{\vec s}{\vec r};\vep\rp
= \hspace{-.5cm}
\sum_{(i_1,\cdots,i_p)\in
\Pi_k}\tilde{P}\llp
    \llp \cdots \check{P}\llp\check{P}\llp \uni{Z}\lp \wvec{\vec s ^{(1)}}{\vec r ^{(1)}};\vep\rp\rrp\,
    \uni{Z}\lp\wvec{\vec s^{(2)}}{\vec r ^{(2)}};\vep\rp \rrp \cdots\rrp
    \uni{Z}\lp\wvec{\vec s^{(p)}}{\vec r ^{(p)}};\vep\rp \rrp } \notag \\
&=& \hspace{-.7cm}
\sum_{(i_1,\cdots,i_p)\in\Pi_k}\hspace{-.6cm}
\tilde{P}\llp \uni{Z}\lp\wvec{\vec s^{(p)}}{\vec r ^{(p)}};\vep\rp
    \check{P}\llp \uni{Z}\lp\wvec{\vec s^{(p-1)}}{\vec r ^{(p-1)}};\vep\rp
    \cdots \check{P}\llp\uni{Z}\lp\wvec{\vec s ^{(1)}}{\vec r ^{(1)}};\vep\rp\rrp\,
     \cdots \rrp
    \rrp.
\mlabel{eq:plus}
\end{eqnarray}
\label{thm:expform}
\end{theorem}
\proofbegin
This follows from Eq.~(\mref{eq:phi-}) and (\mref{eq:phi+}) by induction on $k$. There is nothing to prove when $k=1$. Assuming
the formulas for $\uni{Z}_-$ and $\uni{Z}_+$ are true for $k\leq n$. Then by Eq.~(\mref{eq:phi-}),
$$\uni{Z}_-\lp\wvec{s_1,\cdots,s_{n+1}}{r_1,\cdots,r_{n+1}}\rp
=\check{P}\Big (\uni{Z}\lp \wvec{s_1,\cdots,s_{n+1}}{r_1,\cdots,r_{n+1}}\rp
+\sum_{j=1}^k \uni{Z}_-\lp\wvec{s_1,\cdots,s_j}{r_1,\cdots,r_j}\rp
\uni{Z}\lp\wvec{s_{j+1},\cdots,s_{n+1}}{r_{j+1},\cdots,r_{n+1}}\rp\Big ).$$
Now the formula for $\uni{Z}_-$ follows by applying the induction
hypothesis to the $\uni{Z}_-$ factors in the sum and using the fact
that any ordered partition of $(1,\cdots,n+1)$ is either the one block partition $(n+1)$ or $(i_1,\cdots,i_p,n+1-j), 1\leq j\leq n,$ with $(i_1,\cdots,i_p)$ an ordered partition of $(1,\cdots,j)$. Then the first formula for $\uni{Z}_+$ follows from Eq.~(\mref{eq:phi+}). The second formula
for $\uni{Z}_+$ is just to put the $\uni{Z}$-factors to the front of $\check{P}(x)$ instead of after it.
\proofend

\section{Renormalized multiple zeta values}
\mlabel{sec:exp}
We now use the renormalized directional MZVs defined in Eq. (\mref{eq:dmzv}) to obtain {\bf renormalized MZVs}. Here we will focus on two cases, when the arguments are either all positive or all non-positive.

\subsection{The main definition and theorem}
\begin{mdefn}
For $\vec{s}\in \ZZ_{> 0}^k\cup \ZZ_{\leq 0}^k$, define
the {\bf renormalized MZV} at $\vec{s}$ to be
\begin{equation} \gzeta\lp \vec{s}\rp
= \lim_{\delta \to 0^+} \zeta\lp \wvec{\vec{s}}{|\vec{s}|+\delta}\rp ,
\mlabel{eq:gmzv}
\end{equation}
where, for $\vec s=(s_1,\cdots,s_k)$ and $\delta\in \RR_{>0}$, we denote $|\vec{s}|=(|s_1|,\cdots,|s_k|)$ and $|\vec s|+\delta=(|s_1|+\delta,\cdots,|s_k|+\delta).$

\mlabel{de:rmzv}
\end{mdefn}

\begin{remark} {\rm
Theorem~\mref{thm:main} below is our main theorem. It shows that our renormalized MZVs is well-defined and is compatible with known MZVs defined by either convergence, analytic continuation or the Ihara-Kaneko-Zagier regularization. It also proves that it satisfies the quasi-shuffle relation.
We are optimistic that this is in fact the only definition of $\gzeta(\vec{s})$ from $\zeta(\wvec{\vec s}{\vec r})$ with these properties and will elaborate on it in a subsequent work.
}
\mlabel{re:rmzv}
\end{remark}

\begin{theorem}
The limit in Eq.~(\mref{eq:gmzv}) exists for any $\vec{s}=(s_1,\cdots,s_k)\in \ZZ_{> 0}^k\cup \ZZ_{\leq 0}^k$.
More precisely,
\begin{enumerate}
\item
when $s_i$ are all positive with $s_1>1$, we have $\zeta\lp \wvec{\vec{s}}{\vec{r}}\rp =\zeta(\vec{s})$ independent of $\vec{r}\in \ZZ_{>0}^k$. In particular, $\gzeta(\vec{s})=\zeta(\vec{s})$;
\mlabel{it:g1}
\item
when $s_i$ are all positive, we have $\gzeta(\vec{s})=\zeta\lp \wvec{\vec{s}}{\vec{s}}\rp $. Further, $\gzeta(\vec{s})$ agrees with the regularized MZV $Z^*_{\vec{s}}(T)$ defined by Ihara-Kaneko-Zagier~\mcite{IKZ};
\mlabel{it:ge1}
\item
when $s_i$ are all negative, we have
$\gzeta(\vec{s})=\zeta\lp \wvec{\vec{s}}{-\vec{s}}\rp =\disp{\lim_{\vec{r}\to -\vec{s}} \zeta\lp \wvec{\vec{s}}{\vec{r}}\rp} ;$
\mlabel{it:neg}
\item
when $s_i$ are all non-positive, we have
$\gzeta(\vec{s})= \disp{\lim _{\vec r\to -\vec s}\zeta \lp \wvec{\vec s}{\vec
r}\rp ^{(\Sigma)}}$
where the right hand side is defined in Theorem~\mref{thm:gzeta}.
Further, $\gzeta(\vec{s})$ agrees with $\zeta(\vec{s})$ whenever the later is defined by analytic continuation.
\mlabel{it:nonp}\\
\noindent
{\ ${}$ } \hspace{-1.75cm} Furthermore,
\item
the set $\{\gzeta(\vec{s})\big| \vec{s}\in \ZZ^k_{>0}\}$ satisfies the quasi-shuffle relation;
\mlabel{it:pqs}
\item
the set $\{\gzeta(\vec{s})\big| \vec{s}\in \ZZ^k_{\leq 0}\}$ satisfies the quasi-shuffle relation.
\mlabel{it:nqs}
\end{enumerate}
\mlabel{thm:main}
\end{theorem}
\proofbegin
The items of this theorem will be proved in the rest of this paper.

(\mref{it:g1}) is a restatement of Theorem~\mref{thm:omzv}.
(\mref{it:ge1}) is Theorem~\mref{thm:gregmzv} combined with Proposition~\mref{pp:regzeta}.
(\mref{it:neg}) and the first statement of
(\mref{it:nonp}) are contained in Corollary~\mref{co:sigma}.
The second statement of (\mref{it:nonp}) is Proposition~\mref{pp:anamzv}.
(\mref{it:pqs}) is just Corollary~\mref{co:qsregmzv}.
(\mref{it:nqs}) is just Theorem~\mref{thm:gshuf}.
\proofend

\subsection{Renormalized multiple zeta values with positive arguments}
We first take care of the easy case when MZVs are define by the convergence of
the nested sums.

\begin{theorem}
Let $\vec{s}=(s_1,\cdots,s_k)$ with positive integers $s_1,\cdots,s_k$ and $s_1>1$. We have
$\zeta(\wvec{\vec{s}}{\vec{r}})=\zeta(\vec{s})$, independent of the choice of $\vec{r}\in \ZZ_{>0}^k$. In particular, $\gzeta(\vec{s})=\zeta(\vec{s})$.
\mlabel{thm:omzv}
\end{theorem}
\proofbegin For such an $\vec{s}$, $Z(\wvec{\vec s}{\vec r};\vep)$ is uniformly convergent in
$(-\infty, 0]$, and the summands are continuous functions.
So $Z(\wvec{\vec s}{\vec r};\vep)$ is continuous in $(-\infty, 0]$.
Therefore, the Laurent series of $Z(\wvec{\vec s}{\vec r};\vep)$ is a power series and,
by Theorem \mref {thm:expform},
$\zeta(\wvec{\vec{s}}{\vec{r}})=\lim_{\vep\to 0} \uni{Z}_+(\wvec{\vec{s}}{\vec{r}};\vep)
=\lim_{\vep\to 0} \uni{Z}(\wvec{\vec{s}}{\vec{r}};\vep)
= \uni{Z}(\wvec{\vec{s}}{\vec{r}};0)=\zeta(\vec{s}).$
\proofend

We now extend the last case to
include the possibility of $s_1=1$ and compare it with the regularized
MZVs of Ihara-Kaneko-Zagier~\cite{IKZ}.

Let $(\vec u,\vec v)$ denote the concatenation of two
vectors $\vec u$ and $\vec v$.

\begin {lemma}
For a log power series $f(\vep),g(\vep)\in \CC[T][[\vep]]$, denote
$f(\vep)=g(\vep)+ O(\vep)$ if $g(\vep)-f(\vep)\in \vep \CC[T][[\vep]]$.
Let $\vec s \in \ZZ ^k_{>0}$ be of the form
$\vec{s}=(\vec{1}_{m},\vec{s}\,')$ where $m\geq 1$, $\vec 1_m =(1,1,\cdots, 1)\in \ZZ ^m$ and either $k=m$ or $s_{m+1}>1$.
For $\ell=k-m$, let $\vec r\,' \in \ZZ ^\ell_{>0}$.
For $c>0$, denote $X=-\ln c +T$. Then
\begin{equation}
\uni{Z}(\wvec{\vec 1_m, \vec s\,'}{c \vec 1_m,\vec r\,'};\vep)=P_{m,\vec{s}\,'}(X)+ O(\vep),
\mlabel{eq:lead}
\end{equation}
where $P_{m,\vec{s}\,'}(X)$ is a degree $m$ polynomial in $X$
with leading coefficient $\zeta (\vec s\, ')/m!$ if $\ell>0$ and
$1/m!$ if $\ell=0$.
\mlabel{lem:lead}
\end{lemma}

\proofbegin
We prove by induction on $m\geq 1$.
First consider $m=1$. When $\ell=0$, note that
$$Z(\wvec{1}{c};\vep)=\sum_{n\geq 1} \frac
{e^{nc\vep}}{n}=-\ln(1-e^{c\vep})= -\ln c - \ln (-\vep)+\ln\big(\frac{-c\vep}{1-e^{c\vep}}\big).
$$
Since $ln \frac {-\vep}{1-e^{\vep }}$ is an analytic function at $\vep =0$ with
$\lim _ {\vep \to 0} \ln \frac {-\vep }{1-e^{\vep }}=0,$
we have
\begin{equation}
\uni{Z}(\wvec{1}{c};\vep)= -\ln c +T+O(\vep)=X+O(\vep).
\mlabel{eq:lead1}
\end{equation}
So Eq.~(\mref{eq:lead}) is proved for $m=1$ and $\ell=0$.
When $\ell\geq 1$, let $\vec{e}_j^{(\ell)}$ be the $j$-th unit vector of length $\ell$. Then by the
quasi-shuffle relation
$$
\uni{Z}(\wvec{\vec s\,'}{\vec r\,'}; \vep)\uni{Z}(\wvec{1}{c};\vep)
= \uni{Z}(\wvec{1,\vec s\,'}{c,\vec r\,'})+
\sum_{j=1}^\ell \uni{Z}(\wvec{s_1,\cdots,s_j,1,s_{j+1},\cdots,s_\ell}{r_1,\cdots,r_j,c,r_{j+1},\cdots,r_\ell})
+\sum _{j=1}^\ell \uni{Z}(\wvec{\vec s\,'+ \vec{e}^{(\ell)}_j}{\vec r\,' + c\vec{e}^{(\ell)}_j}; \vep).
$$
Since $s_1>1$, by the proof of Theorem~\mref{thm:omzv} and Eq.~(\mref{eq:lead1}), we have
$$ \uni{Z}(\wvec{1,\vec{s}\,'}{c,\vec{r}\,'})=\zeta(\vec{s}\,')X
-\sum_{j=1}^\ell \zeta(s_1,\cdots,s_j,1,s_{j+1},\cdots,s_\ell)
-\sum _{j=1}^\ell \zeta(\vec s\,'+ \vec{e}^{(\ell)}_j) + O(\vep).$$
This complete the proof for $m=1$.

Suppose the formula has been proved for $m\geq 1$ and consider
$\uni{Z}(\wvec{\vec 1_{m+1}, \vec s\,'}{c \vec 1_{m+1},\vec r\,'};\vep)$.
By the quasi-shuffle relation we have
$$
Z(\wvec{\vec 1_m,\vec s\,'}{c\vec 1_m,\vec r\,'}; \vep)Z(\wvec{1}{c};\vep)
= (m+1)Z(\wvec{\vec 1_{m+1},\vec s\,'}{c\vec 1_{m+1}, \vec r\,'}; \vep)
+\sum _{i=1}^m Z(\wvec{\vec 1_m + \vec{e}^{(m)}_i, \vec s\,'}{c \vec 1_m + c\vec{e}^{(m)}_i, \vec r\,'};
\vep)
+\sum _{j=1}^\ell Z(\wvec{\vec 1_m, \vec s\,'+\vec{e}^{(\ell)}_j}{c\vec 1_m, \vec r \,'+c \vec{e}^{(\ell)}_j}; \vep).
$$
By the induction hypothesis, all terms in the two sums on the right hand side are of the form
$f(X)+O(\vep)$ with $f$ polynomials in $X$ of degree $\leq m$.
Thus by Eq.~(\mref{eq:lead1}) and the induction hypothesis, we obtain
$$Z(\wvec{\vec 1_{m+1}, \vec s\,'}{c \vec 1_{m+1},\vec r\,'};\vep)
=P_{m+1,\vec{s}\,'}(X) + O(\vep),$$
where $P_{m+1,\vec{s}\,'}(X)$ has degree $\deg P_{m,\vec{s}\,'}+1=m+1$ and has leading coefficient the leading coefficient of
$P_{m,\vec{s}\,'}(X)X/(m+1)$, which is $\zeta (\vec s\,')/(m+1)!$
if $\ell>1$ and is $1/(m+1)!$ if $\ell=0$. This completes the induction.
\proofend

\begin {theorem} Let $\vec s \in \ZZ ^k_{>0}$ with
$\vec{s}=(\vec{1}_{m},\vec{s}\,')$ with $m\geq 1$ and $s_{m+1}>1$.
Then
$\gzeta(\vec s)=P_{m,\vec{s}\,'}(T)$, where $P_{m,\vec{s}\,'}$ is the polynomial in Lemma~\mref{lem:lead}. Further, $\gzeta(\vec s)=\zeta\lp\wvec{\vec s}{\vec s}\rp$.
\mlabel{thm:gregmzv}
\end{theorem}
\proofbegin  By Lemma~\mref{lem:lead} and Eq.~(\mref{eq:plus}),
$\gzeta\lp \wvec{\vec 1_m, \vec s\,'}{c\vec 1_m, \vec r\,'}\rp=P_{m,\vec{s}\,'}(X)$
independent of $\vec{r}\,'$. Since $\lim _{c\to 1} X=T=X\big|_{c=1},$
we obtain
$$\lim_{c\to 1, \vec r\,'\to \vec s\,'} \zeta\lp \wvec{\vec 1_m, \vec s\,'}{c\vec 1_m, \vec r\,'}\rp=P_{m,\vec{s}\,'}(T)=\zeta\lp\wvec{\vec s}{\vec s}\rp. \mbox{\proofend}
$$

\begin{coro}
$\gzeta(\vec{s}), \vec{s}\in \ZZ_{>0}^k$, satisfy the quasi-shuffle relation.
\mlabel{co:qsregmzv}
\end{coro}
\proofbegin
The subset $\{{\wvec{s}{s}} |\, s\in \ZZ_{>0}\}$ of the semigroup $\frakM$
in Eq.~(\mref{eq:mbase}) is a subsemigroup and the
$\CC$-space $\calh$ generated by it is a sub-algebra of
the quasi-shuffle algebra $\calh_\frakM$. Thus the algebra homomorphism $\uni{Z}_+:\calh_\frakM\to
\CC[T][[\vep]]$ in Corollary~\mref{co:renormz} restricts to an
algebra homomorphism $\uni{Z}_+:\calh\to \CC[T][[\vep]]$ with
$$\uni{Z}_+\big(\lp{\wvec{s_1}{s_1}}, \cdots, {\wvec{s_k}{s_k}}\rp;0\big)=\zeta\lp\wvec{\vec{s}}{\vec{s}}\rp=\gzeta(\vec s).$$
Hence the corollary.
\proofend

\begin{mprop}
Let $\vec{s}=(s_1,\cdots,s_k)$ with positive integers
$s_1,\cdots,s_k$ and $Z^\msh_{\vec{s}}(T)$ be the regularized MZVs of
Ihara-Kaneko-Zagier. Then $\gzeta(\vec s)=Z^\msh_{\vec{s}}(T)$.
\mlabel{pp:regzeta}
\end{mprop}
\proofbegin
We recall~\cite[Proposition 1]{IKZ} that $Z^\msh_{\vec{s}}(T)$ with
$s_i\geq 1, 1\leq i\leq k$, is obtained as the unique extension of
the MZVs $\zeta(\vec{s})$ with $s_i\geq 1$, $s_1>1$,
such that $Z^\msh_{(1)}(T)=T$ and such that the quasi-shuffle relation
still holds for $Z^*_{\vec{s}}(T)$. Since our definition of
$\gzeta(\vec{s})$ agrees with $\zeta(\vec{s})$ for $s_i\geq 1$ and
$s_1>1$, and our definition of $\gzeta(\vec{s})$ for $s_i \geq 1,
1\leq i\leq k$, also satisfies $\gzeta(1)=T$ (by Theorem~\mref{thm:gregmzv})
and the quasi-shuffle relation (by Corollary~\mref{co:qsregmzv}), the later $\gzeta(\vec{s})$ must
agree with the regularized MZVs $Z^\msh_{\vec{s}}(T)$.
\proofend

\subsection{Renormalized multiple zeta values with non-positive arguments}
We study $\gzeta(\vec{s})$ when $\vec{s}\in \ZZ_{\leq 0}^k$.
We then show that these values agrees with the special values of the multiple zeta functions with negative arguments defined by analytic continuation.

\subsubsection{The case of $\vec{s}=(0,\cdots,0)$}

\begin {mprop}
Let $\vec{s}=(s_1,\cdots,s_k)\in \ZZ_{\le 0}^k$ and $\vec{r}=(r_1,\cdots,r_k)\in \ZZ_{>0}^k$.
\begin{enumerate}
\item
For $k\ge 2$,  we have
\allowdisplaybreaks{
\begin{eqnarray}
\lefteqn{Z(\wvec{\vec s}{\vec r};\vep) =\sum _{j_1=0}^{-s_1}
\binc {-s_1}{j_1}Z(\wvec{-j_1}{r_1};\vep )Z(\wvec{s_1+s_2+j_1, s_3, \cdots
s_{k}}{r_1+r_2, r_3,\cdots, r_{k}};\vep)} \mlabel{eq:zformula}\\
&=& \sum _{j_1=0}^{-s_1} \sum _{j_2=0}^{-s_1-s_2-j_1}\cdots \hspace{-.3cm}\sum
_{j_{k-1}=0}^{-\sum _{\ell=1}^{k-1}s_\ell-\sum _{\ell=1}^{k-2}j_\ell}
\binc {-s_1}{j_1}\binc {-s_1-s_2-j_1}{j_2}\cdots \binc {-\sum
_{\ell=1}^{k-1}s_\ell-\sum _{\ell=1}^{k-2}j_\ell}{j_{k-1}}
\notag \\
&&\times Z(\wvec{-j_1}{r_1};\vep)Z(\wvec{-j_2}{r_1+r_2};\vep)\cdots Z(\wvec{-j_{k-1}}{\sum_{\ell=1}^{k-1}r_\ell};\vep)
Z(\wvec{\sum _{i=1}^{k}s_i+\sum _{i=1}^{k-1}j_i}
    {\sum _{i=1}^{k}r_i};\vep)
\notag
\end{eqnarray}
}
\mlabel{it:ration1}
\item
If $\ell<k$ and $s_1=\cdots=s_\ell=0$, then
$$Z\lp\wvec{\vec s}{\vec r}\rp= Z\lp \wvec{0,\cdots,0}{r_1,\cdots,r_\ell}\rp
    Z\lp \wvec{s_{\ell+1},s_{\ell+2},\cdots,s_k}{r_1+\cdots+r_{\ell+1},r_{\ell+2},\cdots,r_k}\rp.$$
\mlabel{it:ration2} \item Each coefficient in the Laurent  series
expansion of $Z (\wvec{\vec s}{\vec r};\vep )$ is a rational
function of the form $P(\vec r)/Q(\vec r)$, where $P$, $Q$ are in
$\CC[r_1,\cdots,r_k]$ with no common factors, and are of the form
$\Pi_{1\le j \leq k} (r_1+r_{2}+\cdots +r_j)^{c_{j}}$, $c_{j}\in
\ZZ _{\ge 0}$.  \mlabel{it:ration3}
\end{enumerate}
\mlabel{pp:zformula}
\end {mprop}
\proofbegin
(\mref{it:ration1}) Since $s_1\leq 0$, we have
\allowdisplaybreaks{
\begin{eqnarray*}
Z(\wvec{\vec s}{\vec r};\vep)&=&\sum _{n_2>n_3>\cdots n_{k}>0}\frac
{e^{n_2r_2\vep +\cdots +n_{k}r_{k}\vep} }{n_2^{s_2}\cdots
n_{k}^{s_{k}}}\sum _{m=1}^\infty
(n_2+m)^{-s_1}e^{n_2r_1\vep}e^{mr_1 \vep} \\
&=&\sum _{n_2>n_3>\cdots n_{k}>0}\frac {e^{n_2r_2\vep +\cdots
+n_{k}r_{k}\vep} }{n_2^{s_2}\cdots n_{k}^{s_{k}}}\sum _{m=1}^\infty
\sum _{j_1=0}^{-s_1} \binc
{-s_1}{j_1}m^{j_1}n_2^{-s_1-j_1}e^{n_2r_1\vep}e^{mr_1 \vep} \\
&=&\sum _{j_1=0}^{-s_1} \binc {-s_1}{j_1}Z(\wvec{-j_1}{r_1};\vep )
\sum_{n_2>n_3>\cdots n_{k}>0}\frac {e^{n_2(r_1+r_2)\vep +n_3r_3\vep
\cdots +n_{k}r_{k}\vep}
}{n_2^{s_1+s_2+j_1}n_3^{s_3}\cdots n_{k}^{s_{k}}} \\
&=&\sum _{j_1=0}^{-s_1} \binc {-s_1}{j_1}Z(\wvec{-j_1}{r_1};\vep
)
Z(\wvec{s_1+s_2+j_1, s_3, \cdots s_{k}}{r_1+r_2, r_3,\cdots, r_{k}};\vep).
\end{eqnarray*}
}
This gives the proposition when $k=2$. In general, applying the induction
hypothesis to the second $Z$-factor completes the proof.

\noindent
(\mref{it:ration2}) Applying the first equation of item (\mref{it:ration1}) repeatedly, we have
\begin{eqnarray*}
\lefteqn{ Z(\wvec{\vec s}{\vec r};\vep) = Z(\wvec{0}{r_1};\vep )
Z(\wvec{s_2, s_3, \cdots s_{k}}{r_1+r_2, r_3,\cdots, r_{k}};\vep)
=\cdots} \\
&=& Z(\wvec{0}{r_1};\vep)Z(\wvec{0}{r_1+r_2};\vep)\cdots
Z(\wvec{0}{\sum_{i=1}^{\ell}r_i};\vep)
Z(\wvec{s_{\ell+1},\cdots,s_k}
    {r_{\ell+1}+\sum_{i=1}^{\ell}r_i,r_{\ell+2},\cdots,r_k};\vep).
\end{eqnarray*}
Then applying the second equation of item (\mref{it:ration1}) to the product before the last factor
gives item (\mref{it:ration2}).

\noindent
(\mref{it:ration3}) From the generating series of the Bernoulli numbers
$\disp{\frac{\vep}{e^\vep-1}=\sum_{i\geq 0} B_i \frac{\vep^i}{i!}}$,
\begin{equation}
Z(\wvec{0}{1};\vep)=\sum_{n\geq 0} e^{n\vep}
= \frac{e^\vep}{1-e^\vep}= -\frac{1}{\vep}\frac {-\vep}{e^{-\vep}-1}
    = -\frac{1}{\vep}+ \sum_{i\geq 0}\zeta(-i) \frac{\vep^{i}}{i!}
    \mlabel{eq:zeta}
\end{equation}
since $B_0=1$ and $\zeta(-i)=(-1)^i\frac{B_{i+1}}{i+1}$ for $i\geq
0$.
For $s\in \ZZ_{<0}$, we have
$$ Z(\wvec{s}{1};\vep)=\sum_{n\geq 1}{n^{-s}} {e^{n\vep}}
= \frac{d^{-s}}{d\vep} \big(\frac{e^\vep}{1-e^\vep}\big) $$
which converges uniformly on any compact subset in ${\rm
Re}(\vep)<0$. So its Laurent series expansion at $\vep=0$ is obtained
by termwise differentiating Eq.~(\mref{eq:zeta}), yielding
\begin{equation}
Z(\wvec{s}{1};\vep)=(-1)^{s-1}(-s)!\,\vep^{s-1}+\sum _{j=0}^{\infty} \zeta
(s-j)\frac {\vep^j}{j!}. \mlabel{eq:zreg}
\end{equation}
Then for $r\in \ZZ _{>0}$, we have
\begin{equation}
Z(\wvec{s}{r};\vep)=(-1)^{s-1}(-s)!(r\vep)^{s-1}+\sum _{j=0}^{\infty} \zeta
(s-j)\frac {(r\vep)^j}{j!}
\mlabel{eq:zreg2}
\end{equation}
Then item~(\mref{it:ration3}) follows from item~(\mref{it:ration1}).
\proofend

Let $\Sigma_k$ denote the symmetric group on $k$ letters. For $\sigma\in \Sigma_k$ and $\vec{r}=(r_1,\cdots,r_k)$, denote
$\sigma(\vec{r})=(r_{\sigma(1)},\cdots,r_{\sigma(k)})$ and $f(\vec r)^{(\Sigma_k)}=\sum _{\sigma \in \Sigma_k} f(\sigma (\vec r))$.

\begin {mprop} \mlabel{pp:gzero}
Let $k\geq 1$ and $\vec{0}_k=(0,\cdots,0)\in \ZZ^k$.
Then $\zeta\lp\wvec{\vec 0_k}{\vec r}\rp^{(\Sigma_k)}$ is independent of the choice of  $\vec{r}\in \RR_{>0}^k$ and
$\gzeta (\vec 0_k)
=\frac 1{k!}\zeta\lp\wvec{\vec 0_k}{\vec r}\rp^{(\Sigma_k)}.$
\end{mprop}

\proofbegin This is proved by induction on $k\geq 1$. For $k=1$, by Eq.~(\mref{eq:zreg2}) we have
$$\zeta(\wvec{0}{r})= \tilde{P}\big(Z(\wvec{0}{r})\big)\big|_{\vep=0}=\zeta (0)$$
independent of $r>0$. Thus $\gzeta(0)$ is defined and the proposition holds.

In general, by Proposition~\mref{pp:prep} and Corollary~\mref{co:qshuf},
$$\zeta\lp\wvec{\vec 0_{k-1}}{r_1, \cdots, r_{k-1}}\rp^{(\Sigma_{k-1})}\zeta\lp\wvec{0}{r_k}\rp
=\zeta\lp\wvec{\vec 0_k}{\vec r}\rp^{(\Sigma_{k})}
+\sum _{i=1}^{k-1}\zeta\lp\wvec{\vec 0_{k-1}}{r_1, \cdots r_i',
\cdots, r_{k-1}}\rp^{(\Sigma_{k-1})}$$ where $r_i'=r_i+r_k$. So by the induction hypothesis,
$\zeta \lp\wvec{\vec 0_k}{\vec r}\rp^{(\Sigma_{k})}$ is independent of $\vec{r}$. In particular, taking $\vec{r}=(\delta,\cdots,\delta)\in \RR_{>0}^k$, we have $$ \frac{1}{k!}\, \zeta\lp\wvec{\vec{0_k}}{\vec{r}}\rp^{(\Sigma_k)} = \zeta\lp\wvec{\vec{0}_k}{\vec{r}}\rp= \gzeta(\vec{0}_k).
\mbox{\proofend} $$
It will be proved in Theorem~\mref{thm:gshuf} that $\gzeta(\vec{0}_k), k\geq 1,$ satisfy the quasi-shuffle relation. It is in fact the only way to define $\gzeta(\vec{0}_k), k\geq 1,$ with $\gzeta(0)=\zeta(0)$~\cite[Theorem 1.1]{Gu3}.
\subsubsection{The general case of $\vec s \in \ZZ _{\le 0}^k$}
\begin{mdefn} {\rm
Let $\vec{s}\in \ZZ_{\leq 0}^k$. Suppose $s_i=0$ exactly for $k_j'\leq i\leq k_j''$ with $k_j'\leq k_j''$ for $1\leq j\leq q$ and $k_j''<k_{j+1}'-1$ for $1\leq j\leq q-1$. Then $(s_{k_j'},\cdots,s_{k_j''})$ are the longest consecutive zero strings, called the {\bf zero clusters}. So (with the possibility of
$k'_1=1$ or $k''_q=k$)
$$
 \vec{s}=(s_1,\cdots,s_{k_1'-1},\!\!\!\!\underbrace{0,\cdots,0}_{k_1''-k_1'+1-{\rm terms}}\!\!\!\! , s_{k_1''+1},\cdots,s_{k_2'-1},
\cdots, \!\!\!\! \underbrace{0,\cdots,0}_{k_q''-k_q'+1-{\rm terms}}\!\!\!\! ,s_{k_q''+1},\cdots,s_k ).
$$
For each $1\leq i\leq q$, let $\vec{k}^{(i)}=(k_i',\cdots,k_i'')$ and
let $\Sigma_{\vec{k}^{(i)}}$ be the permutation group of $\vec{k}^{(i)}$, naturally a subgroup of $\Sigma_k$. Define
the subgroup
$$ \Sigma(\vec s)=\Sigma_{\vec{k}^{(1)}} \times \cdots \times \Sigma_{\vec{k}^{(q)}} \subseteq \Sigma_k.$$
}
\mlabel{de:cluster}
\end{mdefn}
%
So for $\sigma=(\sigma_1,\cdots,\sigma_q)\in  \Sigma(\vec{s})$ with each $\sigma_i\in \Sigma_{\vec{k}^{(i)}}$,
$1\leq i\leq q$, $\sigma(\vec{r})$ is obtained
by $\sigma_i$ permuting $r_{k_i'},\cdots,r_{k_i''}$ and leaving the other entries fixed.
Define
$$ \gzeta\lp\wvec{\vec{s}}{\vec{r}}\rp^{({\Sigma(\vec{s})})}=\sum_{\sigma\in {\Sigma(\vec{s})}}
\gzeta \lp\wvec{\vec{s}}{\sigma (\vec{r})}\rp,\quad
\uni{Z}_+(\wvec{\vec s}{\vec r}; \vep)^{({\Sigma(\vec{s})})} =\sum_{\sigma\in {\Sigma(\vec{s})}}
\uni{Z}_+(\wvec{\vec s}{\sigma (\vec{r})};\vep).$$

\begin {theorem} The limit
$\disp{\lim _{\vec r\to -\vec s}\zeta\lp \wvec{\vec s}{\vec r}\rp^{({\Sigma(\vec{s})})}}$ exists.
\mlabel{thm:gzeta}
\end {theorem}
We first give some applications of the theorem and defer its proof to Section~\mref{sec:gzeta}.

\begin{coro}
\mlabel{co:sigma}
For $\vec s \in \ZZ ^k_{\le 0}$, $\gzeta(\vec s)$ is well-defined and
$$\gzeta (\vec s)=\frac 1{|{\Sigma(\vec{s})}|}\lim _{\vec r\to -\vec s} \zeta\lp \wvec{\vec s}{\vec r}\rp^{({\Sigma(\vec{s})})}.
$$
If in addition $\vec s$ does not have consecutive zeros, then
$$ \gzeta(\vec{s})=\zeta \lp \wvec{\vec{s}}{-\vec{s}}\rp
=\lim_{\vec{r}\to -\vec{s}} \zeta\lp\wvec{\vec{s}}{\vec{r}}\rp.$$
\end{coro}
\proofbegin
Taking the limit in Theorem~\mref{thm:gzeta} when $\vec
r$ approaches $\vec s$ along the path $\vec r=\vec
s+\delta$, $\delta \to 0$, we have
$$\lim _{\vec r\to -\vec s}\zeta\lp\wvec{\vec s}{\vec
r}\rp^{({\Sigma(\vec{s})})}= \lim _{\delta\to 0} \zeta\lp\wvec{\vec s}{-\vec s+\delta}\rp^{({\Sigma(\vec{s})})}.$$
By the definition of ${\Sigma(\vec{s})}$ and our choice of $\vec r$, ${\Sigma(\vec{s})}$
permutes the components of $\vec r$ that equal $\delta$. Therefore,
$$ \frac{1}{|{\Sigma(\vec{s})}|}\zeta\lp\wvec{\vec s}{-\vec s+\delta}\rp^{({\Sigma(\vec{s})})}= \zeta\lp\wvec{\vec s}{-\vec s+\delta}\rp$$
giving the first limit in the corollary. The second part follows since then $\Sigma(\vec s)$ is trivial.
\proofend

We give an explicit formula when $k=2$. A similar formula holds for $k>2$, expressing
$\gzeta(\vec{s})$ as a polynomial in the Bernoulli numbers. As a consequence,
$\gzeta(\vec{s})$ is rational.
\begin{coro}
Let $s_1$, $s_2 \le 0$, but not both zero. Then $\gzeta(0,s_2)=\zeta(0)\zeta(s_2)-\zeta(s_2-1)$,
and for $s_1<0$,
\allowdisplaybreaks
\begin{eqnarray*}
\lefteqn{\gzeta (s_1, s_2)=\sum _{j=0}^{-s_1} \binc {-s_1}{j}
\zeta (-j)\zeta (s_1+s_2+j)-\frac{1}{1-s_1}\zeta(s_1+s_2-1)}\\
&&+\sum _{j=0}^{-s_1} \binc {-s_1}{j}\frac{(-1)^{s_1+s_2-j+1}}
{-s_1-s_2-j+1}\Big(\frac{s_1+s_2}{s_1}\Big)^{s_1+s_2+j-1}
\zeta(s_1+s_2-1).
\end{eqnarray*}
\mlabel{co:z2formula}
\end{coro}
\proofbegin
By Theorem~\mref{thm:expform}, Proposition~\mref{pp:zformula}
and Eq.~(\mref{eq:zreg2}), we have
\allowdisplaybreaks{
\begin{eqnarray*}
\lefteqn{\zeta \lp\wvec{s_1, s_2}{r_1,r_2}\rp=\sum _{j=0}^{-s_1} \binc {-s_1}{j} \zeta (-j)\zeta (s_1+s_2+j)} \\
&&\hspace{-.7cm} +\sum _{j=0}^{-s_1} \binc {-s_1}{j}\Big( \frac {(-1)^{j+1}}{j+1}\Big(\frac
{r_1+r_2}{r_1}\Big)^{j+1}
\hspace{-.3cm} +\frac{(-1)^{s_1+s_2-j+1}}
{-s_1-s_2-j+1}\Big(\frac
{r_1+r_2}{r_1}\Big)^{s_1+s_2+j-1}\Big)\zeta(s_1+s_2-1)\\
&&\hspace{-.7cm}+\frac{(-1)^{s_1}} {-s_1+1}\Big
(\frac{r_1}{r_2}\Big)^{s_1-1}\zeta(s_1+s_2-1).
\end{eqnarray*}
}
Since
\begin{eqnarray*} &\sum _{j=0}^{-s_1} \binc {-s_1}{j}
\frac {(-1)^{j+1}}{j+1}\Big(\frac {r_1+r_2}{r_1}\Big)^{j+1}=\frac
1{-s_1+1}\sum _{j=1}^{-s_1+1} \binc {-s_1+1}{j} (-1)^{j}\Big(\frac
{r_1+r_2}{r_1}\Big)^{j} \\
=&\frac 1{-s_1+1}\big ((1-\frac {r_1+r_2}{r_1})^{-s_1+1}-1\big
)=\frac 1{-s_1+1}(-\frac {r_2}{r_1})^{-s_1+1}-\frac 1{-s_1+1},
\end{eqnarray*}
the conclusion follows from
Corollary~\mref{co:sigma}. \proofend

In the following table, the element in row $s_1$ and column $s_2$ is $\gzeta(-s_1,-s_2)$,
$1\leq s_1\leq 7,\ 1\leq s_2\leq 8$. It can be seen that elements on each of the even numbered
subdiagonal lines are equal, and that for $s_1=s_2$ even, $\gzeta(s_1,s_2)=0$. They both follow from Eq.~(\mref{eq:negodd}). But the second one follows readily from the quasi-shuffle relation in
Theorem~\mref{thm:main}.(\mref{it:nqs}):
$ 2\gzeta(s_1,s_1)=\zeta(s_1)\zeta(s_1)-\zeta(2s_1)=0.$

{\small
$$
\hspace{-1.5cm}
\begin {array} {ccccccc}
\frac 1{288}& -\frac 1{240}& \frac {83}{64512}& \frac 1{504}& -\frac
{3925}{2239488}& -\frac 1{480}& \frac {342884347}{99656663040}
\medskip
\\
-\frac 1{240}& 0& \frac 1{504}& -\frac {319}{437400}&
-\frac 1{480}& \frac {2494519}{1362493440}& \frac 1{264}
\medskip
\\
-\frac {71}{35840}& \frac 1{504}&
\frac 1{28800}& -\frac 1{480}& \frac {114139507}{139519328256}&
\frac 1{264}& -\frac {313042283533}{93600000000000}
\medskip
\\
\frac 1{504}& \frac {319}{437400}& -\frac 1{480}& 0&
\frac 1{264}& -\frac {41796929201}{26873437500000}& -\frac
{691}{65520}
\medskip
\\
\frac {32659}{15676416}&
-\frac 1{480}& -\frac {21991341}{25836912640}&
\frac 1{264}& \frac 1{127008}& -\frac {691}{65520}& \frac {26194796926873}{5884626295848960}
\medskip\\
-\frac 1{480}& -\frac {2494519}{1362493440}& \frac 1{264}& \frac
{41796929201}{26873437500000}&
-\frac {691}{65520}& 0& \frac 1{24}
\medskip
\\
-\frac {75497471}{19931332608}& \frac 1{264}& \frac
{316292283533}{93600000000000}& -\frac {691}{65520}& -\frac
{36808933898915}{8238476814188544}& \frac 1{24}& \frac 1{115200}
\medskip
\\
\frac 1{264}& \frac
{16608667097}{2879296875000}& -\frac {691}{65520}& -\frac
{4607695}{491051484}& \frac 1{24}& \frac
{63967403428993199}{3561322226607185040}& -\frac {3617}{16320}
\end {array}
$$
}

\subsubsection{Compatibility with multiple zeta values defined by analytic continuation}
\mlabel{sss:neganal}

We recall that the multiple zeta function $\zeta(s_1,\cdots,s_k)$ has analytic continuation to $\CC^k$ with singularities on the subvarieties
in Eq.~(\mref{eq:pole}).
\begin{mprop}
For $(s_1,\cdots,s_k)\in \ZZ_{\leq 0}^k$, if $\zeta(s_1,\cdots,s_k)$ is well-defined by the analytic continuation, then it agrees with
$\gzeta(s_1,\cdots,s_k)$.
\mlabel{pp:anamzv}
\end{mprop}
\proofbegin
When $k=1$, by Eq.~(\mref{eq:zreg2})
for $s\leq 0$,
$$\gzeta(s)= \tilde{P}\big(Z(\wvec{s}{-s})\big)\big|_{\vep=0}=\zeta (s).$$
For $k=2$, by Eq.~(\mref{eq:pole}), exactly when $s_1+s_2$ is negative and
odd, the zeta values $\zeta(s_1,s_2)$ is defined by analytic continuation and thus agrees with the iterated limit
$\disp{ \lim_{z_2\to s_2} \lim_{z_1\to s_1} \zeta(z_1,z_2)}$
defined in~\mcite{A-T}, Eq. (3). Note that our order of arguments in the definition of multiple zeta functions is opposite to their order. So $\zeta(z_1,z_2)$ here is $\zeta(z_2,z_1)$ in their paper. Thus our order of limits here is also opposite to their order.

For $n\geq 0$ and $q\geq 1$, let $(n)_q=n(n+1)\cdots (n+q-1)$. Then by Eq.~(15) in~\mcite{A-T}:
\allowdisplaybreaks{
\begin{eqnarray*}
\zeta(s_1,s_2)&=& \lim_{z_2\to s_2}\lim_{z_1\to s_1}\zeta_2(z_1,z_2)\\
&=&-\frac{\zeta(s_1+s_2-1)}{1-s_1}-\frac{\zeta(s_1+s_2)}{2}
    +\sum_{q=1}^{-s_1} (s_1)_q \frac{(-1)^q}{q!} \zeta(-q) \zeta(s_1+s_2+q).
\end{eqnarray*}
}
Since $s_1+s_2$ is negative and odd, $s_1+s_2-1$ is negative and even. Hence the first term is zero. Further, for $1\leq q\leq -s_1$, either $-q$ or $s_1+s_2+q$ is negative and even. Thus the sum also vanishes, leaving
$\zeta(s_1,s_2)= -\zeta(s_1+s_2)/2$.
By the same argument, from Corollary~\mref{co:z2formula}, we have
\begin{equation} \gzeta(s_1,s_2)=\zeta(0)\zeta(s_1+s_2)
=-\frac{\zeta(s_1+s_2)}{2}=\zeta(s_1,s_2).
\mlabel{eq:negodd}
\end{equation}
By Eq.~(\mref{eq:pole}), for $k\geq 3$, $\zeta(s_1,\cdots,s_k)$ is not defined by analytic continuation for any non-positive integers $s_1,\cdots,s_k$. Thus we have completed the proof.
\proofend

\section{The proof of Theorem~\mref{thm:gzeta}}
\mlabel{sec:gzeta}

\subsection{Reduction to Proposition~\mref{pp:mgood}}
With the notations in Theorem~\mref{thm:expform}, define
$$
 \uni{Z}^m_+\lp \wvec{\vec s}{\vec r};\vep\rp
= \hspace{-.7cm}
\sum_{(i_1,\cdots,i_p)\in\Pi_k}\hspace{-.6cm}
\tilde{P}\llp \uni{Z}\lp\wvec{\vec s^{(p)}}{\vec r ^{(p)}};\vep\rp
    \check{P}\llp \uni{Z}\lp\wvec{\vec s^{(p-1)}}{\vec r ^{(p-1)}};\vep\rp
    \cdots \check{P}\llp \uni{Z}\lp\wvec{\vec s^{(2)}}{\vec r ^{(2)}};\vep\rp
    \check{P}\llp \vep^m\uni{Z}\lp\wvec{\vec s ^{(1)}}{\vec r ^{(1)}};\vep\rp\rrp\,
    \rrp \cdots \rrp
    \rrp.
$$
Then by Theorem~\mref{thm:expform}, we have
$\uni{Z}_+^0\lp \wvec{\vec s}{\vec r};\vep\rp
=\uni{Z}_+\lp \wvec{\vec s}{\vec r};\vep\rp.$
For $\sigma\in \Sigma_k$, let $\sigma(\vec{s})^{(j)}$, $1\leq j\leq p$, be the
partition vectors of $\sigma(\vec{s})=(s_{\sigma(1)},\cdots,s_{\sigma(k)})$ from the ordered partition $(i_1,\cdots,i_p)$
in Definition~\mref{de:part}.
Then
\begin{eqnarray}
\lefteqn{ \uni{Z}^m_+(\wvec{\vec{s}}{\vec{r}};\vep)^{(\sigma)}:=\uni{Z}^m_+(\wvec{\vec{s}}{\sigma(\vec{r})};\vep)}
\mlabel{eq:msum}
\\
&=& \hspace{-.7cm}
\sum_{(i_1,\cdots,i_p)\in\Pi_k}\hspace{-.6cm}
\tilde{P}\llp \uni{Z}\lp\wvec{\vec s^{(p)}}{\sigma(\vec r) ^{(p)}};\vep\rp
    \cdots \check{P}\llp \uni{Z}\lp\wvec{\vec s^{(2)}}{\sigma(\vec r) ^{(2)}};\vep\rp
    \check{P}\llp \vep^m\uni{Z}\lp\wvec{\vec s ^{(1)}}{\sigma(\vec r) ^{(1)}};\vep\rp\rrp\,
    \rrp \cdots \rrp.
\notag
\end{eqnarray}
Define
$$
\uni{Z}^m_+(\wvec{\vec s}{\vec r};\vep)^{(\Sigma(\vec{s}))}
=\sum_{\sigma\in \Sigma(\vec{s})} \uni{Z}^m_+(\wvec{\vec s}{\vec
r};\vep)^{(\sigma)}
$$
By Proposition \mref {pp:zformula}.(\mref{it:ration3}) and
Theorem~\mref{thm:expform},
each coefficient in the Laurent series expansion
of $\uni{Z}^m_+(\wvec{\vec s}{\vec r};\vep)$, and thus of
$\uni{Z}^m_+(\wvec{\vec s}{\vec r};\vep)^{(\Sigma(\vec{s}))}$, is a rational function
$P(\wvec{\vec s}{\vec r})/Q(\wvec{\vec s}{\vec r})\in \CC(\vec{s},\vec{r})$ with $P,Q\in \CC[\vec{s},\vec{r}]$.
We can assume that $P$ and $Q$ have no common factors.
We call this coefficient {\bf ordinary} at $\vec{r}=-\vec{s}$ if $Q(-\vec s)\not= 0$.
We
say that $ \uni {Z} _+^m\lp\wvec{\vec s}{\vec r}\rp^{(\Sigma(\vec{s}))}$ is {\bf
ordinary} if every coefficient of its Laurent series is ordinary.

\begin{lemma}
Let $\sigma\in \Sigma(\vec{s})$. Let $P(\vec{r})/Q(\vec{r})$ be a coefficient of the Laurent series of
$\uni{Z}_+^m(\wvec{\vec s}{\vec r})^{(\Sigma(\vec{s}))}$. The following statements are equivalent.
\begin{enumerate}
\item
$P(\vec{r})/Q(\vec{r})$ is ordinary at $\vec{r}=-\vec{s}$.
\mlabel{it:ord1}
\item
$\disp{\lim_{\vec r\to -\vec s} P(\vec{r})/Q(\vec{r})}$ exists.
\mlabel{it:ord2}
\item
$Q(\vec{r})$ does not have a linear factor
$r_{k_1}+\cdots+r_{k_t}$ such that
$\{k_1,\cdots,k_t\}$ is a subset of $\vec{k}^{(j)}$ for some $1\leq j\leq q$ in
Definition~\mref{de:cluster}.
\mlabel{it:ord3}
\end{enumerate}
\mlabel{lem:ord}
\end{lemma}
\proofbegin
(\mref{it:ord1}) $\Leftrightarrow$ (\mref{it:ord2}) holds for any rational functions, and
(\mref{it:ord1}) $\Rightarrow$ (\mref{it:ord3}) is clear since if $\{k_1,\cdots,k_t\}$ is
a subset of some $\vec{k}^{(j)}$, then $s_{k_1}+\cdots+s_{k_t}=0$.

\noindent
(\mref{it:ord3}) $\Rightarrow$ (\mref{it:ord1}): Suppose $P(\vec{r})/Q(\vec{r})$ is not ordinary
at $\vec{r}=-\vec{s}$.
Then by Proposition~\mref{pp:zformula}.(\mref{it:ration3}) and Theorem~\mref{thm:expform},
$Q(\vec{r})$ has a factor $r_{k_1}+\cdots+r_{k_t}$ with $s_{k_1}=\cdots s_{k_t}=0$.
Since the denominator of a sum of fractions is a factor of the product of the denominators
of the fractions, this factor is a factor in the denominator $Q_\sigma(\vec{r})$ of a coefficient
of the Laurent series of $\uni{Z}_+^m(\wvec{\vec s}{\sigma(\vec r)})$ for a $\sigma\in
\Sigma(\vec{s})$.
If this $\sigma$ is $\id$, then by Proposition~\mref{pp:zformula}.(\mref{it:ration3}) and Theorem~\mref{thm:expform}, $r_{k_1}+\cdots+r_{k_t}=r_i+\cdots+r_{i+t}$ for some $i$.
Thus from $s_{r_{k_1}}=\cdots=s_{r_{k_t}}=0$, $\{r_{k_1},\cdots,r_{k_t}\}=\{r_i,\cdots,r_i+t\}$
is a subset of $\vec{k}^{(j)}$ for some $1\leq j\leq q$ with the notation of Definition~\mref{de:cluster}. If $\sigma\neq \id$, then
$r_{k_1}+\cdots+r_{k_t}=r_{\sigma(i)}+\cdots+r_{\sigma(i+t)}$ for some $i$.
Since $\sigma\in \Sigma(\vec{s})$ permutes the components of $\vec{k}^{(j)}$ among themselves,
$\{r_{\sigma(i)},\cdots,r_{\sigma(i+t)}\}$ is still a subset of $\vec{k}^{(j)}$.
Thus in any case, $Q(\vec{r})$ has a linear factor
$r_{k_1}+\cdots+r_{k_t}$ such that
$\{k_1,\cdots,k_t\}$ is a subset of $\vec{k}^{(j)}$ for some $1\leq j\leq q$.
\proofend

Thus to prove Theorem~\mref{thm:gzeta}, we just apply the following
Proposition~\mref{pp:mgood} to the case when $m=0$ and then let $\vep$ go to $0$.
Note that even though we only need $m=0$ for Theorem~\mref{thm:gzeta}, we have
to consider other values of $m$ for the inductive proof.

\begin{mprop} Let $\vec s \in \ZZ _{\le
0}^k$. Take $m\in \ZZ _{\le 0}$ if $s_1=0$ and take $m\in \ZZ$ if
$s_1<0$. Then $ \uni {Z} _+^m\lp\wvec{\vec s}{\vec r}\rp^{(\Sigma(\vec{s}))}$ is
ordinary at $\vec r=-\vec s$.
\mlabel{pp:mgood}
\end{mprop}

\subsection{The proof of Proposition~\mref{pp:mgood}}
The following chart gives an outline of the proof.
$$
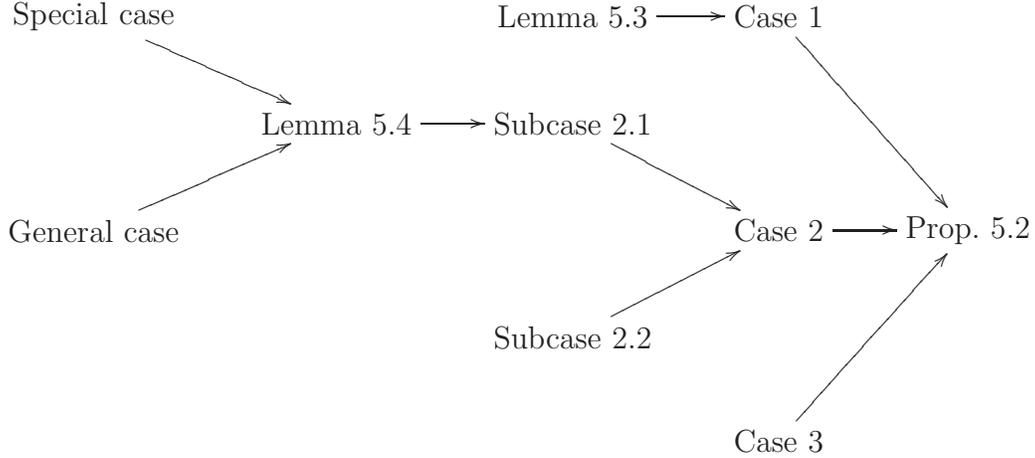
\xymatrix{
\mbox{Special case}\ar[rd] && \mbox{Lemma~\mref{lem:induc}}\ar[r]& \mbox{Case 1} \ar[rdd] & \\
& \mbox{Lemma~\mref{lem:cong}} \ar[r] & \mbox{Subcase 2.1} \ar[rd] && \\
\mbox{General case} \ar[ru] &&& \mbox{Case 2} \ar[r] & \mbox{Prop.~\mref{pp:mgood}}\\
&& \mbox{Subcase 2.2} \ar[ru] &&\\
&&& \mbox{Case 3} \ar[ruu] &
}
$$

We first introduce the following notations to
simplify our expressions. For any Laurent series $f$, we use $\lb f\rb$ to
denote $\check {P}\lp f\rp$. Also use
$\overline{\wvec{\vec{s}}{\vec{r}}}$ to
denote $\uni {Z}\lp \wvec{\vec s}{\vec r};\vep \rp$ and $\overline{\vec r}$ or $\overline{r_1\cdots r_k}$ to denote $\uni{Z}\lp \wvec{\vec 0}{\vec r};\vep\rp$.
With these abbreviations, we have
$$
 \uni{Z}_+^m\lp \wvec{\vec s}{\vec r};\vep\rp
= \sum_{ (i_1,\cdots,i_p)\in \Pi_k} \tilde{P}
    \left ( \overline {\wvec{\vec s^{(p)}}{\vec r ^{(p)}}}
    \lb\cdots \lb \overline {\wvec{\vec s^{( 2 )}}{\vec r ^{( 2)}}}\lb
    \vep ^m \overline {\wvec{\vec s ^{(1)}}{ \vec r ^{( 1 )}}}\rb\,
    \rb \cdots\rb
     \right ).
$$

We now prove Proposition~\mref{pp:mgood} by induction on $k$.
The case when $k=1$ is clear by Eq.~(\mref{eq:zreg2}). Assume that the proposition is true for vectors with length $\le k$
and let $\vec s\in \ZZ _{\le 0}^{k+1}$.
We will separately consider three cases with
$$ \mbox{\bf Case 1: when $s_1<0$; Case 2: $s_1=0$ but $\vec{s}\neq 0$ and
Case 3: $\vec{s}=0$}.$$
\noindent {\bf Case 1: assume $s_1<0$. } Then by our choice,
$m$ is in $\ZZ $. Let $\vec s=(s_1, \vec s\,')$, $\vec r=(r_1, \vec
r\,')$.
We clearly have the disjoint union
\begin{eqnarray}
\Pi_{k+1}&=&\{(1,i_1,\cdots,i_p)\ \big|\ (i_1,\cdots,i_p)\in \Pi_k\} \notag \\
&& \cup
    \{(i_1+1,i_2,\cdots,i_p)\ \big|\ (i_1,\cdots,i_p)\in \Pi_k\}.
\mlabel{eq:ab}
\end{eqnarray}
For $(i_1,\cdots,i_p)\in \Pi_k$, let $\vec{s}^{\,\prime (j)}, 1\leq j\leq p$ be the
partition vectors of $\vec{s}^{\,\prime}$ from $(i_1,\cdots,i_p)$ in Definition~\mref{de:part}.
Similarly define $\vec{r}^{\,\prime (j)}$, $1\leq j\leq p$.
Then by Eq.~(\mref{eq:ab}),
\allowdisplaybreaks{
\begin{eqnarray}
\lefteqn{ \uni {Z} _+^m(\wvec{\vec s}{\vec r})^{(\Sigma(\vec{s}))}
 = \sum_{(i_1,\cdots,i_p)\in \Pi_k}
    \tilde{P}\Llp \overline {\wvec{\vec s^{\,\prime(p)}}{\,\vec r^{\,\prime (p)}}}
    \lb\cdots \lb\overline {\wvec{\vec s^{\,\prime(2)}}{\vec r^{\,\prime (2)}}}
    \lb\overline {\wvec{\vec s^{\,\prime(1)}}{ \vec r^{\,\prime (1)}}}
    \lb\vep ^m \overline {\wvec{s_1}{r_1}}\rb   \rb
    \rb  \cdots\rb \Rrp^{(\Sigma(\vec{s}))}}
\notag \\
&&+\sum_{(i_1,\cdots,i_p)\in \Pi_k}
    \tilde{P}\Big(\overline {\wvec{\vec s^{\,\prime(p)}}{\vec r^{\,\prime (p)}}} \lb \cdots
    \lb \overline {\wvec{\vec s^{\,\prime(2)}}{\vec r^{\,\prime (2)}}}
    \lb \vep ^m \overline {\wvec{s_1,\vec s^{\,\prime (1)}}{r_1, \vec r ^{\,\prime (1)}}}\rb
    \rb  \cdots \rb
    \Big)^{(\Sigma(\vec{s}))}
\mlabel{eq:s19}
\end{eqnarray}
}
We will make use of the following lemma.
\begin{lemma}
Let $\vec{s}\in \ZZ_{\leq 0}^k$, $\vec{r}\in \ZZ_{>0}^k$.
Let $k'<k$ and $\vec{s}=(\vec{s}\,',\vec{s}\,''), \vec{r}=(\vec{r}\,',\vec{r}\,'')$
with $\vec{s}\,'$ and $\vec{r}\,'$ of length $k'$.
Suppose a Laurent series $f(\wvec{\vec{s}\,'}{\vec{r}\,'};\vep)=\sum_{i<0} c_i\vep^i\in \CC(\vec{s}\,',\vec{r}\,')[[\vep,\vep^{-1}]$
is ordinary at $\vec{r}\,'=-\vec{s}\,'$.
If $\uni{Z}_+^m(\wvec{\vec{s}\,''}{\vec{r}\,''};\vep)$ is ordinary at $\vec{r}\,''=-\vec{s}\,''$
for all $m\in \ZZ_{\leq 0}$ (resp. all $m\in \ZZ$), then
$$\sum_{(i_1,\cdots,i_p)\in \Pi_{k-k'}}
    \tilde{P}\Llp \overline {\wvec{\vec s^{\,\prime\prime(p)}}{\vec r^{\,\prime\prime (p)}}}
    \lb\cdots \lb\overline {\wvec{\vec s^{\,\prime\prime (2)}}{\vec r^{\,\prime\prime (2)}}}
    \lb \vep^{m} f(\wvec{\vec{s}\,'}{\vec{r}\,'};\vep) \overline {\wvec{\vec s^{\,\prime\prime(1)}}{ \vec r^{\,\prime\prime (1)}}}
       \rb
    \rb  \cdots\rb \Rrp^{(\Sigma(\vec{s}\,''))}
$$
is ordinary at $\vec{r}=-\vec{s}$
for all $m\in \ZZ_{\leq 0}$ (resp. all $m\in \ZZ$).
\mlabel{lem:induc}
\end{lemma}
\proofbegin
We have
\allowdisplaybreaks{
\begin{eqnarray*}
&&\sum_{(i_1,\cdots,i_p)\in \Pi_{k-k'}}
    \tilde{P}\Llp \overline {\wvec{\vec s^{\,\prime\prime(p)}}{\vec r^{\,\prime\prime (p)}}}
    \lb\cdots \lb\overline {\wvec{\vec s^{\,\prime\prime (2)}}{\vec r^{\,\prime\prime (2)}}}
    \lb \vep^{m} f(\wvec{\vec{s}\,'}{\vec{r}\,'};\vep) \overline {\wvec{\vec s^{\,\prime\prime(1)}}{ \vec r^{\,\prime\prime (1)}}}
    \rb \rb  \cdots\rb \Rrp^{(\Sigma(\vec{s}\,''))} \\
&=& \sum_{i<0} c_i \sum_{(i_1,\cdots,i_p)\in \Pi_{k-k'}}
    \tilde{P}\Llp \overline {\wvec{\vec s^{\,\prime\prime(p)}}{\vec r^{\,\prime\prime (p)}}}
    \lb\cdots \lb\overline {\wvec{\vec s^{\,\prime\prime (2)}}{\vec r^{\,\prime\prime (2)}}}
    \lb \vep^{m+i} \overline {\wvec{\vec s^{\,\prime\prime(1)}}{ \vec r^{\,\prime\prime (1)}}}
    \rb  \cdots\rb \Rrp^{(\Sigma(\vec{s}\,''))}
\end{eqnarray*}
}
Since $i<0$, we have $m+i<0$ if $m\in \ZZ_{\leq 0}$ and $m+i\in \ZZ$ if $m\in \ZZ$.
Hence each of the inner sum is ordinary at $\vec{r}\,''=-\vec{s}\,''$ and thus ordinary at
$\vec r=-\vec s$ since the inner sum does not involve $\vec s\,'$ and $\vec r\,'$.
By assumption each $c_i$ is ordinary at $\vec r\,'=-\vec s\,'$ and hence at $\vec r=-\vec s$ as $c_i$ does not involve $\vec s\,''$ and $\vec r\,''$. Thus the sum is ordinary at $\vec r=-\vec s$.
\proofend

Back to the proof of Proposition~\mref{pp:mgood} in Case 1,
since $s_1<0$, by Eq.~(\mref{eq:zreg2}), $\big\langle \vep ^m \overline {\wvec{s_1}{r_1}} \big\rangle =\check{P}(\vep^m \uni{Z}(\wvec{s_1}{r_1})$
is ordinary at $r_1=-s_1$. Therefore the first sum in Eq.~(\mref{eq:s19}) is ordinary at $\vec r=-\vec s$ by Lemma~\mref{lem:induc} and the induction hypothesis.

For the second term in Eq.~(\mref{eq:s19}), for fixed $1\leq i_1\leq k$, identify an ordered
partition $(i_2,\cdots,i_p)$ of
$k-i_1$ with the ordered partition $(i_1,i_2,\cdots,i_p)$ of $k$, we have
\allowdisplaybreaks{
\begin{eqnarray}
&&\sum_{(i_1,\cdots,i_p)\in \Pi_k}
    \tilde{P}\Big(\overline {\wvec{\vec s^{\,\prime(p)}}{\vec r^{\,\prime (p)}}} \lb\cdots
    \lb \overline {\wvec{\vec s^{\,\prime(2)}}{\vec r^{\,\prime (2)}}}
    \lb \vep ^m \overline {\wvec{s_1,\vec s^{\,\prime (1)}}{r_1, \vec r ^{\,\prime (1)}}}\rb
    \rb  \cdots \rb
    \Big)^{(\Sigma(\vec{s}))}
\notag \\
&=&\sum _{i_1=1}^{k} \sum_{(i_2,\cdots,i_p)\in \Pi_{k-i_1}}
    \tilde{P}\Big(\overline {\wvec{\vec s^{\,\prime(p)}}{\vec r^{\,\prime (p)}}} \lb \cdots
    \lb \overline {\wvec{\vec s^{\,\prime(2)}}{\vec r^{\,\prime (2)}}}
    \lb \vep ^m \overline {\wvec{s_1,\vec s^{\,\prime (1)}}{r_1, \vec r ^{\,\prime (1)}}}\rb
    \rb  \cdots \rb
    \Big)^{(\Sigma(\vec{s}))}
\mlabel{eq:s20}
\end{eqnarray}
}
Let $\Sigma'=\Sigma(s_1,\vec{s}\,'^{(1)})$
and $\Sigma''=\Sigma(\vec{s}\,'^{(2)},\cdots,\vec{s}\,'^{(p)})=\Sigma(s_{i_1+1},\cdots,s_k)$. Then $\Sigma'\Sigma''=\Sigma'\times \Sigma''$ and thus is a subgroup of $\Sigma(\vec{s})$. Let $S$ be a complete set of coset representatives
for the cosets $\{ \Sigma'\Sigma'' \sigma\ |\ \sigma\in \Sigma(\vec{s})\}$ of $\Sigma(\vec{s})$.
So
\begin{eqnarray}
&&\sum_{(i_2,\cdots,i_p)\in \Pi_{k-i_1}}
    \tilde{P}\Big(\overline {\wvec{\vec s^{\,\prime(p)}}{\vec r^{\,\prime (p)}}} \lb \cdots
    \lb \overline {\wvec{\vec s^{\,\prime(2)}}{\vec r^{\,\prime (2)}}}
    \lb \vep ^m \overline {\wvec{s_1,\vec s^{\,\prime (1)}}{r_1, \vec r ^{\,\prime (1)}}}\rb
    \rb  \cdots \rb\Big)^{(\Sigma(\vec{s}))}
\mlabel{eq:s21}
\\
&=&
\sum_{\sigma\in S} \sum_{(i_2,\cdots,i_p)\in \Pi_{k-i_1}}
    \tilde{P}\Big(\overline {\wvec{\vec s^{\,\prime(p)}}{\sigma(\vec r^{\,\prime})^{ (p)}}} \lb\cdots
    \lb \overline {\wvec{\vec s^{\,\prime(2)}}{\sigma(\vec r^{\,\prime})^{ (2)}}}
    \lb \vep ^m \overline {\wvec{s_1,\vec s^{\,\prime (1)}}{\sigma(r_1, \vec r ^{\,\prime (1)})}}\rb ^{(\Sigma')}
    \rb  \cdots \rb \Big)^{(\Sigma'')}.
\notag
\end{eqnarray}

Let $\sigma\in S$, by the definition of $\Sigma(\vec{s})$, we have
$\sigma(r_1,\vec{r}\,^{\prime (1)})=(r_1,r_{\sigma(2)},\cdots).$
So by Proposition~\mref{pp:zformula}, each coefficient in the
Laurent series expansion of $\lb \vep ^m \overline {\wvec{s_1,\vec
s^{\,\prime (1)}}{\sigma(r_1, \vec r ^{\,\prime} (1)})}\rb
^{(\Sigma')}$ has its denominator as a product of
$(r_1+r_{\sigma(j_1)}+\cdots +r_{\sigma(j_t)})$, $t\geq 0$, with
$1\leq j_1,\cdots,j_t\leq i_1$. Hence the expansion is ordinary at
 $\sigma(\vec r\,^{\prime})^{(1)}=-\sigma(\vec s\,^{\prime})^{(1)}$ since $s_1<0$ and $s_i\leq 0$.

Further by the induction hypothesis on $k$,
$$\sum_{(i_2,\cdots,i_p)\in \Pi_{k-i_1}}
    \tilde{P}\Big(\overline {\wvec{\vec s^{\,\prime(p)}}{\sigma(\vec r^{\,})^{\prime (p)}}} \lb\cdots
    \lb \vep^t \overline {\wvec{\vec s^{\,\prime(2)}}{\sigma(\vec r^{\,})^{\prime (2)}}}
    \rb  \cdots \rb \Big)^{(\Sigma'')}
$$
is ordinary at $(r'_{\sigma(i_1+1)},\cdots,r'_{\sigma(k+1)})=-(s_{i_1+1},\cdots, s_{k+1})$ for $t\in \ZZ_{\leq 0}$.
Thus by Lemma~\mref{lem:induc}, the inner sum on the
right hand side of Eq.~(\mref{eq:s21}) is ordinary at $\sigma(\vec r)=-\vec s$ for each $\sigma\in S$. Note that $\sigma(\vec{s})=\vec{s}$ for $\sigma\in \Sigma(\vec{s})$ and being ordinary at
$\sigma(\vec{r})=-\sigma(\vec{s})$ is equivalent to being ordinary at $\vec{r}=-\vec{s}$.
Hence the left hand sum is ordinary at $\vec r=-\vec s$ in
Eq.~(\mref{eq:s21}) and hence in Eq.~(\mref{eq:s20}) and hence in the second term of Eq.~(\mref{eq:s19}).

\medskip

\noindent {\bf Case 2:  assume $s_1=0$, but $\vec s\neq 0$.} Then
$m\le 0$. Assume $s_1=s_2=\cdots =s_\ell=0$, $s_{\ell+1}\not =0$,
$\ell \le k$.
Then $\vec s =(\vec s\,', \vec s\,'')$ with $\vec s\,'=\vec{0}_\ell$ and
$\vec{s}\,''=(s_{\ell+1},\cdots,s_{k+1})$. Similarly denote
$\vec r= (\vec r\,',\vec r\,'')$ with $\vec r\,'=(r_1,\cdots,r_\ell)$ and
$\vec r\,''=(r_{\ell+1},\cdots, r_{k+1})$.
In the notation of Definition~\mref{de:cluster}, $(1,\cdots,\ell)=\vec{k}^{(1)}$.

By Lemma~\mref{lem:ord}, in order to prove that $\uni{Z}_+^m(\wvec{\vec{s}}{\vec{r}};\vep)$ is
ordinary, we only need to show that no coefficient of its Laurent series expansion
has a denominator with either
\begin{enumerate}
\item a {\bf type (i) factor:} $r_{k_1}+\cdots+r_{k_t}$ where $k_1,\cdots,k_t\leq \ell$ or
\item a {\bf type (ii) factor:} $r_{k_1}+\cdots+r_{k_t}$ where $\{k_1,\cdots,k_t\}$ is a subset of $\vec{k}^{(j)}$ for some $2\leq j\leq q$ in
Definition~\mref{de:cluster}.

\end{enumerate}

\noindent
{\bf Subcase 2.1: there are no type (i) factors.}
Note that any ordered partition of $\Pi
_{k+1} $ is of the form $(i_1, i_2, \cdots, i_p, j_1, j_2, \cdots,
j_q)$ or $(i_1, i_2, \cdots, i_p+j_1, j_2, \cdots, j_q)$, with
$(i_1, i_2, \cdots, i_p)\in \Pi _\ell$, $(j_1, j_2, \cdots, j_q)\in
\Pi _{k+1-\ell}$.
Using the notations in Definition~\mref{de:part}, let $\vec{r}^{\,\prime(1)},\cdots,\vec{r}^{\,\prime(p)}$ be the partial vectors
of $\vec r\,'$ from the ordered partition $(i_1,\cdots,i_p)\in \Pi_\ell$.
Similarly let $\vec{s}^{\,\prime\prime(1)},\cdots, \vec{s}^{\,\prime\prime(q)}$ (resp.
$\vec{r}^{\,\prime\prime(1)},\cdots,\vec{r}^{\,\prime\prime(q)}$) be the partition vectors
of $\vec s\,''$ (resp. $\vec r\,''$) from the ordered partition $(j_1,\cdots,j_q)\in \Pi_{k+1-\ell}$.
Then we have
\begin{eqnarray*}
 \uni {Z} _+^m\lp\wvec{\vec s}{\vec r}\rp^{(\Sigma(\vec{s}))}&= & \hspace{-1cm}
 \sum_{\tiny{\begin{array}{c}(i_1,\cdots,i_p)\in \Pi_\ell\\
(j_1, j_2, \cdots, j_q)\in \Pi _{k+1-\ell}\end{array}}}
   \hspace{-1cm} \tilde{P}\Big(\overline {\wvec{\vec s^{\,\prime \prime(q)}}{\vec r^{\,\prime \prime (q)}}}\lb\cdots
   \lb \overline {\wvec{\vec s^{\,\prime \prime (1)}}{\vec r^{\,\prime \prime (1)}}}
   \lb \overline {\vec r^{\,\prime (p)}} \cdots
   \lb \overline {\vec r^{\,\prime (2)}}
   \lb \vep ^m \overline {\vec r^{\,\prime (1)}}\rb
    \rb  \cdots \rb
    \rb \cdots \rb
    \Big)^{(\Sigma(\vec{s}))}
\\
&& +  \hspace{-1.3cm}\sum_{\tiny{\begin{array}{c}(i_1,\cdots,i_p)\in \Pi_\ell\\
(j_1, j_2, \cdots, j_q)\in \Pi _{k+1-\ell}\end{array}}}
    \hspace{-1cm} \tilde{P}\Big(\overline {\wvec{\vec s^{\,\prime \prime(q)}}{\vec r^{\,\prime \prime (q)}}} \lb \cdots
    \lb \overline {\wvec{\vec 0, \vec s^{\,\prime \prime (1)}}{\vec r^{\,\prime (p)},\vec r^{\,\prime \prime (1)}}}\cdots
    \lb \overline {\vec r^{\,\prime (2)}}
    \lb \vep ^m \overline {\vec r^{\,\prime (1)}}\rb
    \rb  \cdots \rb \cdots \rb
    \Big)^{(\Sigma(\vec{s}))}
\end{eqnarray*}
Also define $\Sigma'=\Sigma_{\ell}=\Sigma(\vec{s}\,')$ and $\Sigma''=\Sigma(\vec{s}\,'')$ as in Definition~\mref{de:cluster}.
Then $\Sigma(\vec{s})=\Sigma_\ell \times \Sigma''$. So
\begin{eqnarray}
\lefteqn{\uni {Z} _+^m\lp\wvec{\vec s}{\vec r}\rp^{(\Sigma(\vec{s}))} \hspace{-.5cm}=\hspace{-1.3cm}\sum_{\tiny{\begin{array}{c}(i_1,\cdots,i_p)\in \Pi_\ell\\
(j_1, j_2, \cdots, j_q)\in \Pi _{k+1-\ell}\end{array}}}
    \hspace{-1cm} \tilde{P}\Big(\overline {\wvec{\vec s^{\,\prime \prime(q)}}{\vec r^{\,\prime \prime (q)}}} \lb \cdots
    \lb \overline {\wvec{\vec s^{\,\prime \prime (1)}}{\vec r^{\,\prime \prime (1)}}}
    \lb \overline {\vec r^{\,\prime (p)}}\cdots
    \lb \overline {\vec r^{\,\prime (2)}}
    \lb \vep ^m \overline {\vec r^{\,\prime (1)}}\rb
    \rb  \cdots \rb ^{(\Sigma _\ell)}
    \rb \cdots \rb
    \Big)^{(\Sigma'')}} \notag \\
&& + \hspace{-1.3cm} \sum_{\tiny{\begin{array}{c}(i_1,\cdots,i_p)\in \Pi_\ell\\
(j_1, j_2, \cdots, j_q)\in \Pi _{k+1-\ell}\end{array}}}
    \hspace{-1.3cm} \tilde{P}\Big(\overline {\wvec{\vec s^{\,\prime \prime(q)}}{\vec r^{\,\prime \prime (q)}}} \lb \cdots
    \lb \overline {\wvec{\vec 0, \vec s^{\,\prime \prime (1)}}{\vec r^{\,\prime (p)},\vec r^{\,\prime \prime (1)}}}\cdots
    \lb \overline {\vec r^{\,\prime (2)}}
    \lb \vep ^m \overline {\vec r^{\,\prime (1)}}\rb
    \rb  \cdots \rb ^{(\Sigma _\ell)} \cdots \rb
    \Big)^{(\Sigma'')}
\mlabel{eq:mzsum}
\end{eqnarray}

Now for fixed $(j_1,\cdots,j_q)\in\Pi _{k+1-\ell}$ and $\tau \in
\Sigma''$, the corresponding terms in the above summation are
\begin{eqnarray*}
&\sum_{\tiny{(i_1,\cdots,i_p)\in \Pi_\ell }}
    \tilde{P}\Big(\overline {\wvec{\vec s^{\,\prime \prime(q)}}{\tau(\vec r^{\,\prime \prime})^{(q)}}} \lb \cdots
    \lb \overline {\wvec{\vec s^{\,\prime \prime (1)}}{\tau(\vec r^{\,\prime \prime})^{ (1)}}}
    \lb \overline {\vec r^{\,\prime(p)}}\cdots
    \lb \overline {\vec r^{\,\prime (2)}}
    \lb \vep ^m \overline {\vec r^{\,\prime (1)}}\rb
    \rb  \cdots \rb ^{(\Sigma _\ell)}\rb \cdots \rb
    \Big) & \\
&+ \sum_{\tiny{(i_1,\cdots,i_p)\in \Pi_\ell}}
    \tilde{P}\Big(\overline {\wvec{\vec s^{\,\prime \prime(q)}}{\tau(\vec r^{\,\prime \prime})^{ (q)}}} \lb \cdots
    \lb \overline {\wvec{\vec 0, \vec s^{\,\prime \prime (1)}}{\vec r^{\,\prime (p)},\tau(\vec r^{\,\prime \prime})^{ (1)}}}\cdots
    \lb \overline {\vec r^{\,\prime (2)}}
    \lb \vep ^m \overline {\vec r^{\,\prime (1)}}\rb
    \rb  \cdots \rb ^{(\Sigma _\ell)} \cdots \rb
    \Big)&\\
&= \tilde{P}\Big(\overline {\wvec{\vec s^{\,\prime \prime(q)}}{\tau(\vec r^{\,\prime \prime})^{(q)}}} \lb \cdots
    \overline {\wvec{\vec s^{\,\prime \prime(2)}}{\tau(\vec r^{\,\prime \prime})^{(2)}}}
    \llb \sum_{\tiny{(i_1,\cdots,i_p)\in \Pi_\ell }}\Big (
  \overline {\wvec{\vec s^{\,\prime \prime (1)}}
  {\tau(\vec r^{\,\prime \prime})^{(1)}}}
  \lb \overline {\vec r^{\,\prime (p)}} \cdots
  \lb \overline {\vec r^{\,\prime (2)}}
  \lb \vep ^m \overline {\vec r^{\,\prime (1)}}
  \rb \rb  \cdots \rb
 & \\
  & \qquad \qquad \qquad
  +\overline {\wvec{\vec 0, \vec s^{\,\prime \prime (1)}}{\vec r^{\,\prime (p)},
  \tau(\vec r^{\,\prime \prime})^{ (1)}}}
    \cdots
  \lb \overline {\vec r^{\,\prime (2)}}
  \lb \vep ^m \overline {\vec r^{\,\prime (1)}}\rb
    \rb  \cdots
    \Big )^{(\Sigma _\ell)}
    \rrb  \cdots \rb
    \Big) &
\end{eqnarray*}

By Proposition~\mref {pp:zformula}.(\mref{it:ration2}),
$$\overline {\wvec{\vec 0, \vec s^{\,\prime \prime (1)}}{\vec r^{\,\prime (p)},\tau(\vec r^{\,\prime \prime})
^{(1)}}}=\overline { \vec r^{\,\prime (p)}}\overline {\wvec{\vec s^{\,\prime
\prime (1)}}{\tau(\vec r^{\,\prime \prime})^{ (1)}+(r'_{I_{p-1}+1}+\cdots +
r'_{I_p})\vec{e} _1^{(j_1)}}},
$$
where $\vec{e}^{(j_1)}_1$ is the first unit vector of length $j_1$ (which is the length of $\vec{r}\,^{\prime\prime (1)}$ and $\tau(\vec{r}\,'')^{(1)}$).
So the inner sum on the right hand side above becomes
\begin{eqnarray*}
&&\sum_{(i_1,\cdots,i_p)\in \Pi _\ell}\left (
  \overline {\wvec{\vec s^{\,\prime \prime (1)}}{\tau(\vec r^{\,\prime \prime})^{(1)}}}
  \lb \overline {\vec r^{\,\prime (p)}} \cdots
  \lb \overline {\vec r^{\,\prime (2)}}
  \lb \vep ^m \overline {\vec r^{\,\prime (1)}}\rb
  \rb  \cdots \rb \right .\\
&& \qquad \left . +\,
  \overline { \vec r^{\,\prime (p)}}
  \overline {\wvec{\vec s^{\,\prime\prime (1)}}
  {\tau(\vec r^{\,\prime \prime})^{ (1)}+(r'_{I_{p-1}+1}+\cdots +r'_{I_p})\vec{e} _1^{(j_1)}}} \lb \cdots
  \lb \overline {\vec r^{\,\prime (2)}}
  \lb \vep ^m \overline {\vec r^{\,\prime (1)}}\rb
    \rb  \cdots \rb
    \right )^{(\Sigma _\ell)}\\
&=&\sum_{(i_1,\cdots,i_p)\in \Pi _\ell}\sum_{\sigma\in \Sigma_{\ell}} \left (
  \overline {\wvec{\vec s^{\,\prime \prime (1)}}{\tau(\vec r^{\,\prime \prime})^{(1)}}}
  \lb \overline {\sigma(\vec r^{\,\prime})^{ (p)}} \cdots
  \lb \overline {\sigma(\vec r^{\,\prime})^{ (2)}}
  \lb \vep ^m \overline {\sigma(\vec r^{\,\prime})^{ (1)}}\rb
  \rb  \cdots \rb \right .\\
&& \qquad \left . +\,
  \overline { \sigma(\vec r^{\,\prime})^{ (p)}}
  \overline {\wvec{\vec s^{\,\prime\prime (1)}}
  {\tau(\vec r^{\,\prime \prime})^{ (1)}+(r'_{\sigma(I_{p-1}+1)}+\cdots +r'_{\sigma(I_p)})\vec{e} _1^{(j_1)}}} \lb \cdots
  \lb \overline {\sigma(\vec r^{\,\prime})^{ (2)}}
  \lb \vep ^m \overline {\sigma(\vec r^{\,\prime})^{ (1)}}\rb
    \rb  \cdots \rb
    \right )
\end{eqnarray*}
For a given pair $g:=(\sigma,\pi)\in \Sigma_\ell \times \Pi_\ell$,
let $(g)_1$ (resp. $(g)_2$) be the first term (resp. second term)
in the above sum. Thus the double sum can be simply denoted by
$$\sum _{g\in \Sigma_\ell\times \Pi_\ell}((g)_1+(g)_2).
$$

We denote
$f \equiv _\ell g$ if no coefficient of the Laurent series expansion of $f-g$
has a denominators with a factor $r'_\ell+r'_{k_1}+\cdots +r'_{k_t}$, $t\geq 0$, with $j_1,\cdots,j_t<\ell$. This is clearly an equivalence relation.
\begin{lemma}
$\disp{\sum _{g\in \Sigma_\ell \times \Pi_\ell}((g)_1+(g)_2) \equiv_\ell 0.}
$
\mlabel{lem:cong}
\end{lemma}
\proofbegin
\noindent
{\bf A special case.}
We first consider the special case when
$g=(\sigma,\pi)\in \Sigma_\ell\times\Pi_\ell$ is of the form
$(\cdots(\ell))$, that is, $\pi$ has $(\ell)$ as the last partition factor and $\sigma(\ell)=\ell$.
Denote $a= \lb \overline {\sigma(\vec r^{\,\prime})^{ (p-1)}}\cdots
  \lb \overline {\sigma(\vec r^{\,\prime})^{ (2)}}
  \lb \vep ^m \overline {\sigma(\vec r^{\,\prime})^{ (1)}}\rb
    \rb  \cdots \rb$.
Then we verify that
\begin{eqnarray}
(g)_1+(g)_2 &=&
  \overline {\wvec{\vec s^{\,\prime \prime (1)}}{\tau(\vec r^{\,\prime \prime})^{(1)}}}
  \lb\bar r' _\ell  \lb a\rb \rb
    +
   \bar r' _\ell\overline {\wvec{\vec s^{\,\prime
\prime (1)}}{\tau(\vec r^{\,\prime \prime})^{ (1)}+r'_{\ell} \vec{e} _1^{(j_1)}}} \lb a \rb
\notag \\
&& \equiv _\ell
  -\overline {\wvec{\vec s^{\,\prime \prime (1)}}{\tau(\vec r^{\,\prime \prime})^{(1)}}}
  \lb \lb \bar r' _\ell\rb \lb a \rb  \rb
    -  \lb \bar r' _\ell\rb
    \overline {\wvec{\vec s^{\,\prime\prime (1)}}
   {\tau(\vec r^{\,\prime \prime})^{ (1)}+r'_{\ell}\vec{e} _1^{(j_1)}}}
    \lb  a \rb
\mlabel{eq:lcong}\\
&&=\overline {\wvec{\vec s^{\,\prime \prime (1)}}{\tau(\vec r^{\,\prime \prime})^{ (1)}}}
    \lb \bar r' _\ell\rb \lb a\rb
    -\overline {\wvec{\vec s^{\,\prime \prime (1)}}{\tau(\vec r^{\,\prime \prime})^{ (1)}+r'_{\ell}\vec{e} _1^{(j_1)}}}
    \lb \bar r' _\ell\rb \lb  a\rb.
\notag
\end{eqnarray}
Here for the $\equiv_\ell$, note that
the series expansion of $\lb \bar r'_\ell \rb + \bar r'_\ell$ is
the power series part of the Laurent series of $\bar r'_\ell = \uni{Z}(\wvec{0}{r'_\ell};\vep)$, and hence by Eq.~(\mref{eq:zreg2}),
$r'_\ell$ does not occur in the denominators
the series expansion of $\lb \bar r'_\ell \rb + \bar r'_\ell$.
Then the $\equiv_\ell$ follows from the easily checked properties: if $f\equiv_\ell g$, then $ \lb f \rb \equiv \lb g\rb$, and if in addition $h\not\equiv_\ell 0$, then $fh \equiv_\ell gh$.
The last equation in Eq.~(\mref{eq:lcong}) holds since
$\check{P}$ is an idempotent Rota-Baxter operator and hence by
Eq.~(\mref{eq:rbe}),
$$ \lb \lb x \rb \lb y \rb \rb= \lb \lb \lp x \lb y\rb + \lb x\rb y +xy \rp \rb\rb
=\lb \lp x \lb y\rb + \lb x\rb y +xy \rp \rb = \lb x\rb \lb y \rb.$$

By Eq.~(\mref{eq:zreg2}), $\lb \bar {r'_\ell}\rb =-\frac 1{r'_\ell \vep}$.
Thus
\begin{eqnarray*}
&&\lim_{r'_\ell\to 0}\Big (\overline {\wvec{\vec s^{\,\prime \prime (1)}}{\tau(\vec r^{\,\prime \prime})^{(1)}}}
    \lb \bar {r'_\ell}\rb
    -\overline {\wvec{\vec s^{\,\prime\prime (1)}}
    {\tau(\vec r^{\,\prime \prime})^{ (1)}+r'_{\ell}\vec{e}_1^{(j_1)}}}
    \lb \bar {r'_\ell}\rb \Big) \\
&&    =\frac{1}{\vep} \lim_{r'_\ell\to 0} \Big(\overline {\wvec{\vec s^{\,\prime \prime (1)}}{\tau(\vec r^{\,\prime \prime})^{(1)}}}
    -\overline {\wvec{\vec s^{\,\prime\prime (1)}}
    {\tau(\vec r^{\,\prime \prime})^{ (1)}+r'_{\ell}\vec{e}_1^{(j_1)}}}\Big)/r'_\ell
    =\frac{1}{\vep} \frac{\partial}{\partial (r''_1)} \overline {\wvec{\vec s^{\,\prime \prime (1)}}{\tau(\vec r^{\,\prime \prime})^{(1)}}}
\end{eqnarray*}
exists. Here the differentiation is taken termwise in the Laurent series. Thus
$$\overline {\wvec{\vec s^{\,\prime \prime (1)}}{\tau(\vec r^{\,\prime \prime})^{(1)}}}
    \lb \bar {r'_\ell}\rb
    -\overline {\wvec{\vec s^{\,\prime\prime (1)}}
    {\tau(\vec r^{\,\prime \prime})^{ (1)}+r'_{\ell}\vec{e}_1^{(j_1)}}}
    \lb \bar {r'_\ell}\rb \equiv _\ell 0.
$$

\noindent
{\bf The general case. } We now prove the lemma in general by
induction on $\ell$. If $\ell=1$, then there can be only one
partition $(1)$. So the special case applies and we are done.
Assume the lemma is proved for $\ell-1$ and consider the case of
$\ell$.

For $g=(\sigma,(i_1,\cdots,i_p))\in \Sigma_\ell\times \Pi_\ell$, consider
the vector of partition vectors of $(\sigma(1),\cdots,\sigma(\ell))$ from $(i_1,\cdots,i_p)$, called a {\bf partitioned permutation},
$$((\sigma(1),\cdots,\sigma(I_1)),(\sigma(I_1+1),\cdots,\sigma(I_2)),\cdots,
    (\sigma(I_{p-1}+1),\cdots,\sigma(\ell))).$$
Here $I_j=i_1+\cdots+i_j$, $1\leq j\leq p$.
So $I_p=\ell$.
This gives a natural 1-1 correspondence between $\Sigma_\ell\times \Pi_\ell$ and
\begin{equation}
\Pi_\ell^\Sigma:=\{
((n_1, \cdots n_{I_1}),(n_{I_1+1}, \cdots
n_{I_2}),\cdots ,(n_{I_{p-1}+1}, \cdots, n_{\ell}))\}
\mlabel{eq:parper}
\end{equation}
where $(I_1,I_2-I_1,\cdots,\ell-I_{p-1})$ is in $\Pi_\ell$ and
$(n_1,\cdots,n_{\ell})$ is in $\Sigma_\ell$. We can thus identify
$\Sigma_\ell\times \Pi_\ell$ with $\Pi_\ell^\Sigma$ and
call $p=\length(g)$ the length of $g$.

For $1\leq p\leq \ell$, let
\begin{itemize}
\item
$\Sigma_{\ell,\,\leq p}\subseteq \Pi_\ell^\Sigma$ consisting of $g$ with $\length(g)\leq p$,
\item
$\Sigma_{\ell,\leq p}^{(1)}$ consisting of $g\in \Sigma_{\ell,\,\leq p}$ whose last partition factor is not $(\ell)$,
\item
$\Sigma_{\ell, \leq p}^{(2)}$ consisting of $g\in \Sigma_{\ell,\,\leq p}$ that do not contain $(\ell)$ as a partition factor,
\item
$\Sigma_{\ell, \leq p}^{(3)}=\Sigma_{\ell,\,\leq p}^{(1)}\backslash \Sigma_{\ell,\,\leq p}^{(2)}$, that is, consisting of $g\in \Sigma_{\ell,\,\leq p}$ that do contain $(\ell)$ as a partition factor, but not as the last factor.
\end{itemize}
Similarly define $\Sigma_{\ell,=p}$ and $\Sigma_{\ell, =p}^{(i)}$ for $i=1,2,3.$
Thus $\Sigma_{\ell,\,\leq \ell}=\Sigma_\ell^\Pi$ and by the special case, we have
$$ \sum_{g\in \Sigma_\ell^\Pi\backslash\Sigma_{\ell,\,\leq \ell}^{(1)}} ((g)_1+(g)_2)
\equiv_\ell 0.$$
So to prove Lemma~\mref{lem:cong} we only need to prove
\begin{equation}
 \sum_{g\in \Sigma_{\ell,\,\leq \ell}^{(1)}} ((g)_1+(g)_2)
\equiv_\ell 0.
\mlabel{eq:cong1}
\end{equation}
For this we first use the induction on $p=\length(g)$ to prove that for $(g)=(g)_1$ or $(g)_2$,

\begin{equation}
\sum_{g\in \Sigma_{\ell,\le p}^{(1)}}(g)
\equiv _\ell \sum _{g\in \Sigma_{\ell,=p}^{(2)}}(g)
-\sum _{i=1}^{\ell-1}\sum _{g\in \Sigma_{\ell-1,\leq p-1}}(g)^{(i)}
 \mlabel {eqn:length}
\end{equation}
In the last term $(g)^{(i)}$ means replacing $(r'_1,\cdots,r'_{\ell-1})$ by $(r'_1, \cdots, r'_i+r'_\ell ,\cdots, r'_{\ell -1})$ in $(g)$.

The case of $\ell=1$ and thus $p=1$ is covered by the special
case. For $\ell \ge 2$, we use induction on $p$. When $p=1$, there is only one partition.
So Eq.~(\mref{eqn:length}) is an identity.
If for $p$,
formula (\mref {eqn:length}) is true, then we have
\begin{eqnarray}
\sum_{g\in \Sigma_{\ell,\le (p+1)}^{(1)}}(g)
&=& \sum_{g\in \Sigma_{\ell,\le p}^{(1)}}(g)
+ \sum_{g\in \Sigma_{\ell,= (p+1)}^{(1)}}(g) \mlabel{eq:gind}
\\
&\equiv _\ell& \sum _{g\in \Sigma_{\ell,=p}^{(2)}}(g)
-\sum _{i=1}^{\ell-1}\sum _{g\in \Sigma_{\ell-1, \leq p-1}}(g)^i +\sum _{g\in \Sigma_{\ell,=p+1}^{(2)}}(g)
+\sum _{g\in \Sigma_{\ell,=p+1}^{(3)}}(g).
\notag
\end{eqnarray}

It is easily verified that the following relations are equivalence
relations.
\begin{itemize}
\item An element in $\Sigma_{\ell, =p+1}^{(3)}$ is of the form
$(\cdots (\ell)(a_1,\cdots,a_j)\cdots )$ with $\{a_1,\cdots,a_j\}
\subseteq [\ell]$. Define $g_1=(\cdots (\ell)(a_1, \cdots,
a_j)\cdots )\sim_3 g_2$ if $g_2$ can be obtained from $g_1$ by a
permutation of $(a_1,\cdots,a_j)$. Thus an equivalence class for
$\sim_3$ is of the form
$(\cdots(\ell)(a_1,\cdots,a_j)^{\Sigma_j}\cdots)$. \item An
element in $\Sigma_{\ell, =p}^{(2)}$ is of the form $(\cdots (a_1,
\cdots, a_j)\cdots )$ with $\ell\in \{a_1,\cdots,a_j\}$. Define
$g_1=(\cdots (a_1, \cdots, a_j)\cdots )\sim_2 g_2$, where $\ell\in
\{a_1,\cdots,a_j\}$, if $g_2$ can be obtained from $g_1$ by by a
permutation of $(a_1, \cdots a_j)$. Thus an equivalence class for
$\sim_2$ is of the form
$(\cdots(a_1,\cdots,a_j)^{\Sigma_{j}}\cdots)$.

\item An element in $[\ell -1]\times \Sigma_{\ell-1,=p}$ is of the
form $\big(i, (\cdots (\sigma(I_{j-1}+1), \cdots,
\sigma(I_j))\cdots)\big)$ where $\sigma\in \Sigma_{\ell-1}$,
$(I_{j-1}+1,\cdots,I_j)$ is a block of an ordered partition of
$\ell-1$ of length $p$ and $I_{j-1}+1\leq i\leq I_j$. Define $(i,
(\cdots (\sigma(I_{j-1}+1),\cdots,\sigma(I_j))\cdots) \sim
(i',g')$ if $I_{j-1}+1\leq i' \leq I_j$ and $g'$ can be obtained
from $g$ by a permutation of $(\sigma(I_{j-1}+1),\cdots,
\sigma(I_j))$. An equivalence class for $\sim$ is of the form
$$\cup_{i=I_{j-1}+1}^{I_j} (i, (\cdots(\sigma(I_{j-1}+1),\cdots,\sigma(I_j))
^{\Sigma_{I_j-I_{j-1}+1}}, \cdots)).$$
\end{itemize}

There are obvious one-to-one correspondences between these equivalence classes
\allowdisplaybreaks{
\begin{eqnarray}
&& \psi:  \Sigma_{\ell,=p+1}^{(3)}/\!\!\sim_3\,
\longrightarrow \Sigma_{\ell, =p}^{(2)}/\!\!\sim_2
\mlabel{eq:psi}
\\
&& (\cdots (\ell)(a_1,\cdots,a_j)^{\Sigma_j}\cdots) \mapsto
(\cdots (\ell,a_1,\cdots,a_j)^{\Sigma_{j+1}}\cdots),
\notag \\
&& \rho:  \Sigma_{\ell,=p+1}^{(3)}/\sim_3
 \longrightarrow   [\ell-1]\times \Sigma_{\ell-1,=p}/\sim
\mlabel{eq:rho}\\
&& (\cdots (\ell)(a_1,\cdots,a_j)^{\Sigma_j}\cdots)
\mapsto \cup_{i=I_{j-1}+1}^{I_j} (i,(\cdots(a_1,\cdots,a_j)
^{\Sigma_{j}}\cdots)). \notag
\end{eqnarray}
}
Here $(a_1,\cdots,a_j)=(\sigma(I_{j-1}+1),\cdots,\sigma(I_j))$.

For $g=(\cdots (\ell)(a_1, \cdots, a_j)\cdots )
\in \Sigma_{\ell,=p+1}^{(3)}$, as in Eq.~(\mref{eq:lcong}) we have
\begin{eqnarray*}
 (g)&=&\lb \lb \cdots \rb \bar r'_\ell \rb \overline {r'_{a_1} \cdots r'
_{a_j}}\cdots \equiv _\ell -\lb \lb \cdots \rb \lb \bar r'_\ell\rb
\rb \overline {r'_{a_1} \cdots r' _{a_j}}\cdots \\
&=&\lb \cdots \rb \lb
\bar r'_\ell \rb \overline {r'_{a_1} \cdots r' _{a_j}}\cdots \equiv
_\ell-\lb \cdots \rb \bar r'_\ell \overline {r'_{a_1} \cdots r'
_{a_j}}\cdots.
\end{eqnarray*}
By Proposition~\mref{pp:prep} and Eq.~(\mref{eq:zmap}),  we have
$$\bar {r'_\ell}\ \overline{r'_{a_1}r'_{a_2}\cdots r'_{a_{j}}}^{(\Sigma_{j})}
=\overline{r'_{\ell} r'_{a_1}r'_{a_2}\cdots r'_{a_{j}}}^{(\Sigma_{j+1})}+\sum
_{i=1}^{j}\overline{r'_{a_1}\cdots \tilde{r}_{a_i}'\cdots
r'_{a_j}}^{(\Sigma_{j})}
$$
where $\tilde{r}_{a_i}'=r'_{a_i}+r'_\ell$.
Using Eq.~(\mref{eq:psi}) and (\mref{eq:rho}), we obtain
$$\sum _{h\sim_3 g}(h)=-\sum _{h\sim_2 \psi (g)}(h)-\sum _{h\sim
\rho (g)}(h).
$$
Here, for $h=(i,g)$, $(h)=(g)^{(i)}$.
Summing over all the equivalence classes, we have
$$
\sum _{h\in \sigma_{\ell,=p+1}^{(3)}}(h)
=-\sum _{h\in \sigma_{\ell,=p}^{(2)}}(h)
-\sum _{h\in [\ell-1]\times \Sigma_{\ell-1,=p}} (h)
=-\sum _{h\in \sigma_{\ell,=p}^{(2)}}(h)
-\sum_{i=1}^{\ell-1} \sum _{h\in \Sigma_{\ell-1,=p}} (h)^{(i)}.
$$
Combining this with Eq.~(\mref{eq:gind}) gives
$$ \sum_{g\in \Sigma_{\ell,\le p+1}^{(1)}}(g)
\equiv _\ell  \sum _{g\in \Sigma_{\ell,=p+1}^{(2)}}(g)
-\sum _{i=1}^{\ell-1}\sum _{g\in \Sigma_{\ell-1, \leq p}}(g)^{(i)},
$$
completing the inductive proof of Eq.~(\mref {eqn:length}).

Take $p=\ell+1$ in Eq.~(\mref {eqn:length}). Since the maximal length of
an ordered partition of $\ell$ is $\ell$, we have
$$ \sum_{g\in \Sigma_{\ell,\le \ell+1}^{(1)}}((g)_1+(g)_2)
\equiv _\ell
-\sum _{i=1}^{\ell-1}\sum _{g\in \Sigma_{\ell-1, \leq \ell}}((g)_1^{(i)}+(g)_2^{(i)})
=-\sum _{i=1}^{\ell-1}\sum _{g\in \Pi_{\ell-1}^\Sigma}((g)_1^{(i)}+(g)_2^{(i)}).
$$
Now by the induction hypothesis on $\ell$, the right hand side is  $\equiv _\ell 0$. On the other hand, by its definition,
$\Sigma_{\ell,\,\leq \ell+1}^{(1)}=\Sigma_{\ell,\,\leq \ell}^{(1)}.$
Therefore Eq.~(\mref{eq:cong1}), and hence Lemma~\mref{lem:cong},
is proved.
\proofend

Thus we have proved
\begin{eqnarray*}
&&\sum_{\Pi _\ell}\left (
  \overline {\wvec{\vec s^{\,\prime \prime (1)}}{\vec r^{\,\prime \prime(1)}}}
  \lb \overline {\vec r^{\,\prime (p)}} \lb \cdots
  \lb \overline {\vec r^{\,\prime (2)}}
  \lb \vep ^m \overline {\vec r^{\,\prime (1)}}\rb
    \rb  \cdots \rb \rb
     \right .\\
&&  + \left .
    \overline {\wvec{\vec s^{\,\prime\prime (1)}}
    {\vec r^{\,\prime \prime (1)}+(r'_{I_{p-1}+1}+\cdots +r'_{I_p})\vec{e} _1^{(1)}}}
    \overline { \vec r^{\,\prime (p)}} \lb \cdots
    \lb\overline {\vec r^{\,\prime (2)}}
    \lb \vep ^m \overline {\vec r^{\,\prime (1)}}\rb
    \rb  \cdots \rb \right )^{(\Sigma _\ell)}\equiv _\ell0
\end{eqnarray*}
Because the action of $\Sigma_\ell$, the role of $\ell$ in the above expression is symmetric to any
$1\leq t\leq \ell-1$. Thus no coefficient of its Laurent series has a denominator with
a homogeneous linear factor $r'_{a_1}+\cdots +r'_{a_m}$ with $\{a_1,
\cdots , a_m \}\subset \{1,\cdots , \ell \}$. This completes Subcase 2.1.

\noindent
{\bf Subcase 2.2: There are no type (ii) factors.}
Now for a fixed $\pi=(i_1,\cdots,i_p)\in \Pi_{\ell}$ and $\sigma\in \Sigma_\ell$, let
$$\lb \overline {\sigma(\vec r^{\,\prime})^{(p)}}\cdots
    \lb \overline {\sigma(\vec r^{\,\prime})^{ (2)}}
    \lb \vep ^m \overline {\sigma(\vec r^{\,\prime})^{(1)}}
    \rb\rb  \cdots \rb
    =\sum_{i<0} c_i^{\pi,\sigma} \vep^i \in \CC(\vec{s}\,',\vec{r}\,')[[\vep,\vep^{-1}]
$$
be the Laurent series expansion.
Then the first sum for $\uni{Z}^m_+(\wvec{\vec{s}}{\vec{r}};\vep)^{(\Sigma(\vec{s}))}$ in Eq.~(\mref{eq:mzsum}) becomes
\begin{eqnarray*}
&& \hspace{-1cm}\sum_{\tiny{\begin{array}{c}(i_1,\cdots,i_p)\in \Pi_\ell \\ \tiny{(j_1,\cdots,j_q)\in \Pi_{k+1-\ell}}\end{array}}}
    \hspace{-.3cm}\sum_{\sigma\in \Sigma_\ell}
     \tilde{P}\Big(
    \overline {\wvec{\vec s^{\,\prime \prime(q)}}{\vec r^{\,\prime \prime (q)}}}\cdots
    \lb \overline {\wvec{\vec s^{\,\prime \prime (1)}}{\vec r^{\,\prime \prime (1)}}}
    \lb \overline {\sigma(\vec r^{\,\prime})^{(p)}}\cdots
    \lb \overline {\sigma(\vec r^{\,\prime})^{ (2)}}
    \lb \vep ^m \overline {\sigma(\vec r^{\,\prime})^{(1)}}\rb
    \rb  \cdots \rb \rb \cdots
    \Big)^{(\Sigma'')} \\
&=& \sum_{\tiny{\begin{array}{c}(i_1,\cdots,i_p)\in \Pi_\ell \\ \tiny{(j_1,\cdots,j_q)\in \Pi_{k+1-\ell}}\end{array}}}
    \hspace{-.3cm}\sum_{\sigma\in \Sigma_\ell}
     \tilde{P}\Big(
    \overline {\wvec{\vec s^{\,\prime \prime(q)}}{\vec r^{\,\prime \prime (q)}}}\cdots
    \lb \sum_{i<0} c_i^{\pi,\sigma} \vep^i
    \overline {\wvec{\vec s^{\,\prime \prime (1)}}{\vec r^{\,\prime \prime (1)}}}\rb \cdots
    \Big)^{(\Sigma'')} \\
&=& \sum_{\tiny{\begin{array}{c}(i_1,\cdots,i_p)\in \Pi_\ell \\ \sigma\in \Sigma_\ell\end{array}}}
     \sum_{i<0} c_i^{\pi,\sigma}  \sum_{\tiny{(j_1,\cdots,j_q)\in \Pi_{k+1-\ell}}}
    \tilde{P}\Big(
    \overline {\wvec{\vec s^{\,\prime \prime(q)}}{\vec r^{\,\prime \prime (q)}}}\cdots \lb \vep^i
    \overline {\wvec{\vec s^{\,\prime \prime (1)}}{\vec r^{\,\prime \prime (1)}}}\rb \cdots
    \Big)^{(\Sigma'')}
\end{eqnarray*}
By the induction hypothesis on $k$, the inner most sum has a Laurent series expansion
whose coefficients have denominators with no type (ii) factors. Since $c_i^{\pi,\sigma}$ is
a rational function in $r_1,\cdots,r_\ell$, the same can be said of the whole sum.

Similarly, let
$$ \overline {\sigma(\vec r^{\,\prime})^{(p)}} \lb \cdots
    \lb \overline {\sigma(\vec r^{\,\prime})^{ (2)}}
    \lb \vep ^m \overline {\sigma(\vec r^{\,\prime})^{(1)}}
    \rb\rb  \cdots \rb
  =\sum_{-\infty<i} d_i^{\pi,\sigma} \vep^i \in \CC(\vec{s}\,',\vec{r}\,')[[\vep,\vep^{-1}]
$$
be the Laurent series expansion. Then the second sum for $\uni{Z}^m_+(\wvec{\vec{s}}{\vec{r}};\vep)^{(\Sigma(\vec{s}))}$ in Eq.~(\mref{eq:mzsum}) becomes
\begin{eqnarray*}
&& \hspace{-1.5cm}\sum_{\tiny{\begin{array}{c}(i_1,\cdots,i_p)\in \Pi_\ell\\
(j_1, j_2, \cdots, j_q)\in \Pi _{k+1-\ell}\end{array}}}
    \hspace{-1.3cm} \tilde{P}\Big(
    \overline {\wvec{\vec s^{\,\prime \prime(q)}}{\vec r^{\,\prime \prime (q)}}}\cdots
    \lb \overline {\wvec{\vec s^{\,\prime\prime (1)}}
    {\vec r^{\,\prime \prime (1)}+(r'_{I_{p-1}+1}+\cdots +r'_{I_p})\vec{e} _1^{(1)}}}
    \Big( \overline { \vec r^{\,\prime (p)}}\cdots
    \lb \overline {\vec r^{\,\prime (2)}}
    \lb \vep ^m \overline {\vec r^{\,\prime (1)}}\rb
    \rb  \cdots \Big)^{(\Sigma_\ell)}
    \rb  \cdots\Big)^{(\Sigma'')}
\\
&=& \hspace{-.8cm}\sum_{\tiny{\begin{array}{c}(i_1,\cdots,i_p)\in \Pi_\ell\\ \sigma\in \Sigma_\ell\end{array}}}
     \sum_{-\infty < i} d_i^{\pi,\sigma}\hspace{-.4cm}  \sum_{\tiny{(j_1,\cdots,j_q)\in \Pi_{k+1-\ell}}}
    \hspace{-.7cm} \tilde{P}\Big(
    \overline {\wvec{\vec s^{\,\prime \prime(q)}}{\vec r^{\,\prime \prime (q)}}}\cdots
    \lb \vep^i \overline {\wvec{\vec s^{\,\prime\prime (1)}}
    {\vec r^{\,\prime \prime (1)}+(r'_{\sigma(I_{p-1}+1)}+\cdots +r'_{\sigma(I_p)})\vec{e} _1^{(1)}}}
    \rb \cdots\Big)^{(\Sigma'')}
\end{eqnarray*}
Since $I_{p-1}+1,\cdots,I_{p}\leq \ell$ and $\sigma\in \Sigma_\ell$, we have $\sigma(I_{p-1}+1),\cdots,\sigma(I_p)\leq \ell$. Thus when $\vec r\to -\vec s$, we have
$r'_{\sigma(I_{p-1}+1)},\cdots, r'_{\sigma(I_p)} \to 0$. Therefore, for the inner sum above,
$$ \sum_{\tiny{(j_1,\cdots,j_q)\in \Pi_{k+1-\ell}}}
    \hspace{-.7cm} \tilde{P}\Big(
    \overline {\wvec{\vec s^{\,\prime \prime(q)}}{\vec r^{\,\prime \prime (q)}}}\cdots
    \lb \vep^i \overline {\wvec{\vec s^{\,\prime\prime (1)}}
    {\vec r^{\,\prime \prime (1)}+(r'_{\sigma(I_{p-1}+1)}+\cdots +r'_{\sigma(I_p)})\vec{e} _1^{(1)}}}
    \rb \cdots\Big)^{(\Sigma'')}\Big |_{\vec r=-\vec s}
$$
exists since it equals to
$${
\sum_{\tiny{(j_1,\cdots,j_q)\in \Pi_{k+1-\ell}}}
    \hspace{-.7cm} \tilde{P}\Big(
    \overline {\wvec{\vec s^{\,\prime \prime(q)}}{\vec r^{\,\prime \prime (q)}}}\cdots
    \lb \vep^i \overline {\wvec{\vec s^{\,\prime\prime (1)}}
    {\vec r^{\,\prime \prime (1)}}} \rb \cdots\Big)^{(\Sigma'')}_{\Big |_{\vec r=-\vec s}}
= \hspace{-.3cm}\sum_{\tiny{(j_1,\cdots,j_q)\in \Pi_{k+1-\ell}}}
    \hspace{-.7cm} \tilde{P}\Big(
    \overline {\wvec{\vec s^{\,\prime \prime(q)}}{\vec r^{\,\prime \prime (q)}}}\cdots
    \lb \vep^i \overline {\wvec{\vec s^{\,\prime\prime (1)}}
    {\vec r^{\,\prime \prime (1)}}} \rb \cdots\Big)^{(\Sigma'')}_{\Big |_{\vec r\,''=-\vec s\,''}}
}
$$
which exists by the induction hypothesis on $k$ and Lemma~\mref{lem:ord}. Then by Lemma~\mref{lem:ord} again, the above inner sum is ordinary at $\vec r=-\vec s$ and hence at $\vec r\,'=-\vec s\,'$
and so is free of type (ii) factors in the denominators of it Laurent series coefficients.
Therefore the whole expression has a Laurent series expansion
whose coefficients do not have any denominators with a type (ii) factor since $d_i^{\pi,\sigma}$
does not involve $\vec s\,''$ and $\vec r\,''$.

Thus $\uni{Z}^m_+(\wvec{\vec{s}}{\vec{r}})^{(\Sigma(\vec{s}))}$ has no type (ii) factors, completing
the proof of Case 2.

\noindent {\bf Case 3: assume $\vec s=\vec 0$}. The proof is basically the same as for Case 2 except that there is not type (ii) factors to exclude.
For $\sigma\in \Sigma_{k+1}$ and $\pi=(i_1,\cdots,i_p)\in \Pi_{k+1}$, as in Eq.~(\mref{eq:parper}),
let
$$g:=((n_1,\cdots,n_{I_1}),(n_{{I_1}+1},\cdots,n_{I_2}),\cdots,(n_{I_{p-1}+1},\cdots,n_{I_p}))\in \Sigma_{k+1}^\Pi$$
be the corresponding partitioned permutation. Let
$$(g):= \tilde{P}\Big(\overline{\sigma(\vec{r})^{(p)}} \lb \cdots \lb\overline{\sigma(\vec{r})^{(2)}}\lb \vep^m \overline{\sigma(\vec{r})^{(1)}}\rb \rb \cdots \rb  \Big).
$$
Then
$ \uni{Z}^m_+(\wvec{\vec{s}}{\vec{r}};\vep)= \sum_{g\in \Sigma_{k+1}^\Pi} (g).$
We first prove the following analog of Lemma~\mref{lem:cong} by adapting its proof.
\begin{equation}
\sum_{g\in \Sigma_{k+1}^\Pi} (g) \equiv _{k+1} 0.
\mlabel{eq:cong2}
\end{equation}
Suppose $g$ has the last partition factor as $(k+1)$. Then $\sigma(k+1)=k+1$.
Let $\lb a \rb = \lb \cdots \lb\overline{\sigma(\vec{r})^{(2)}}\lb \vep^m \overline{\sigma(\vec{r})^{(1)}}\rb \rb \cdots \rb
=\sum_{i<0} a_i \vep^i$. So $a_i$ are rational functions in $\{r_{\sigma(1)},\cdots,r_{\sigma(k)}\}
=\{r_1,\cdots,r_k\}$.
Then by Eq.~(\mref{eq:zreg2}), we have
\begin{eqnarray*}
(g)&=&\tilde{P}( \overline{r_{k+1}}\lb a\rb)\\
&=& \sum_{i<0} a_i \tilde{P}\llp\vep^i ( -(r_{k+1}\vep)^{-1}+\sum_{j=0}^\infty \zeta(-j)\frac{(r_{k+1}\vep)^j}{j!} )\rrp\\
&=& \sum_{i<0} a_i \tilde{P}\llp\vep^i (\sum_{j=0}^\infty \zeta(-j)\frac{(r_{k+1}\vep)^j}{j!} )\rrp
\end{eqnarray*}
since $\tilde{P}=\id-{P}$ is the projection of a Laurent series to its power series part.
Thus no linear factor involving $r_{k+1}$ appears in the denominator of any coefficient
of the Laurent series expansion of $(g)$.

Thus using the notation of the General Case in the proof of Lemma~\mref{lem:cong}, we only
need to prove
$
\sum_{g\in \Sigma_{k+1, \leq k+1}^{(1)}} (g) \equiv _{k+1} 0.
$
For this we prove the following analog of Eq.~(\mref{eqn:length}) by using the same proof.
\begin{eqnarray}
\sum_{g\in \Sigma_{k+1,\le p}^{(1)}}(g)
\equiv _{k+1} \sum _{g\in \Sigma_{{k+1},=p}^{(2)}}(g)
-\sum _{i=1}^{k}\sum _{g\in \Sigma_{k,\leq p-1}}(g)^{(i)}.
 \mlabel {eqn:length1}
\notag
\end{eqnarray}
Then the rest of the proof of Lemma~\mref{lem:cong} carries through and gives Eq.~(\mref{eq:cong2}).
Then again similar to Case 2, the symmetry of $r_1,\cdots, r_{k+1}$ in $\uni{Z}_+^m(\wvec{\vec{s}}{\vec{r}};\vep)$
shows that it is ordinary.
\medskip

We have completed our inductive proof of Proposition~\mref{pp:mgood} in all three cases.

\section{The quasi-shuffle relation for non-positive MZVs}
\mlabel{sss:nonpoqs}
The purpose of this section is to prove the following theorem and thus
to complete the proof of Theorem~\mref{thm:main}.

\begin {theorem} $\gzeta(\vec{s})$, $\vec{s}\in
\ZZ_{\leq 0}^k$, satisfy the quasi-shuffle relation.
\mlabel{thm:gshuf}
\end{theorem}
%
\subsection{A lemma on stuffles}
For the proof of Thoerem~\mref{thm:gshuf}, we use the stuffle interpretation of the quasi-shuffle product.
The mathematics formulation of stuffles already appeared in Cartier's construction of free commutative Rota-Baxter algebras~\mcite{Ca2} in 1972, even though stuffle was defined using the same formulation in the study of MZVs 20 years later~\mcite{3BL,Br}.  It is well-known in the literatures of MZVs that quasi-shuffle product is the same as the stuffle product~\mcite{Ho3}. To see it in another way, it was proved in~\mcite{G-K1} that the stuffle product, in the variation of Cartier, is equivalent to the mixable shuffle product and the mixable shuffle product is shown in~\mcite{E-G1} to to be the same as the quasi-shuffle product.

For an integer $n\geq 1$, denote $[n]=\{1,\cdots,n\}$ which is also
identified with the vector $(1,\cdots,n)$.
For integers $k,\ell\geq 1$, a {\bf $(k,\ell)$-stuffle triple} is a triple $(r,\alpha,\beta)$ in
$$\frakS:= \left \{ (r,\alpha,\beta)\ \left |  \begin{array}{l}
\max(k,\ell)\leq r\leq k+\ell, \\
\alpha:[k]\to [r], \beta: [\ell]\to [r] \mbox{\ order preserving injections}, \\
\im(\alpha)\cup \im(\beta)=[r] \end{array} \right . \right \}$$
Thus for each $1\leq u\leq r$, at least one of $\alpha^{-1}(u)$ and $\beta^{-1}(u)$ is a singleton $\{w\}$ which we just write $w$. Similarly denote
$$\tilde{\frakS}:= \left \{ (r,\alpha,\beta)\ \left |  \begin{array}{l}
\max(k,\ell)\leq r\leq k+\ell, \\
\alpha:[k]\to [r], \beta: [\ell]\to [r] \mbox{\ injective (might not preserve order)}, \\
\im(\alpha)\cup \im(\beta)=[r] \end{array} \right . \right \}$$

Let $\vec x=(x_1, \cdots, x_k)$ and $\vec y=(y_1,\cdots, y_\ell)$ be vectors of symbols.
The {\bf stuffle of $\vec x$ and $\vec{y}$} corresponding to $(r,\alpha,\beta)$ is defined by, in the notations of~\mcite{Br,Ca2},
{
\begin{equation}
\hspace{-.5cm}\vec z:= \Phi_{r,\alpha,\beta}(\vec{x},\vec{y})=(z_1,\cdots,z_r), \
z_u=\left \{ \begin{array}{ll}
x_{\alpha^{-1}(u)}, & \alpha^{-1}(u)\neq \emptyset, \beta^{-1}(u)=\emptyset,\\
y_{\beta^{-1}(u)}, & \alpha^{-1}(u)= \emptyset, \beta^{-1}(u) \neq \emptyset,\\
x_{\alpha^{-1}(u)}y_{\beta^{-1}(u)}, & \alpha^{-1}(u)\neq \emptyset, \beta^{-1}(u) \neq \emptyset
\end{array} \right .
\mlabel{eq:Phi}
\end{equation}
}
With the convention that $x_{\emptyset}=y_{\emptyset}=\bfone$ and thus
$x_{\alpha^{-1}(u)}y_{\beta^{-1}(u)}=x_{\alpha^{-1}(u)}$ if $\beta^{-1}(u)=\emptyset$ and
$x_{\alpha^{-1}(u)}y_{\beta^{-1}(u)}=y_{\beta^{-1}(u)}$ if $\alpha^{-1}(u)=\emptyset$,
we simply have
\begin{equation}
\Phi_{r,\alpha,\beta}(\vec{x},\vec{y})=(x_{\alpha^{-1}(1)}y_{\beta^{-1}(1)},\cdots,
    x_{\alpha^{-1}(r)}y_{\beta^{-1}(r)}).
\mlabel{eq:Phi2}
\end{equation}
Note that this definition makes sense even for $(r,\alpha,\beta)\in \tilde{\frakS}$.
More generally, for subvectors (that is, subsequences) $\vec{k}\subv,\vec{\ell}\subv,
\vec{r}\subv=(r_{i_1},\cdots,r_{i_{p\subv}})$ of $[k],[\ell],[r]$ respectively with
$\alpha(\vec{k}\subv)\cup \beta(\vec{\ell}\subv)=\vec{r}\subv$, consider the corresponding
subvectors $\vec{x}\subv,\vec{y}\subv,\vec{z}\subv$ of
$\vec{x},\vec{y},\vec{z}$ respectively.
Define
\begin{equation}
\Phi_{r\subv,\alpha|_{\vec{k}\subv},\beta|_{\vec{\ell}\subv}}
(\vec{x}\subv,\vec{y}\subv)
=(x_{\alpha^{-1}(r_{i_1})}y_{\beta^{-1}(r_{i_1})}, \cdots,
x_{\alpha^{-1}(r_{i_{p\subv}})}y_{\beta^{-1}(r_{i_{p\subv}})}).
\mlabel{eq:phi1}
\end{equation}
\begin{lemma}
\begin{enumerate}
\item
For $\sigma\in \Sigma_k, \tau\in \Sigma_\ell$,
$ \Phi_{r,\alpha,\beta}(\sigma(\vec{x}),\tau(\vec{y}))= \Phi_{r,\alpha\circ \sigma^{-1}, \beta \circ \tau^{-1}}(\vec{x},\vec{y}).$
\mlabel{it:phi-equal}
\item
Distinct triples $(r,\alpha,\beta)$ in $\tilde{\frakS}$
give distinct vectors $\Phi_{r,\alpha,\beta}(\vec{x},\vec{y})$.
\mlabel{it:phi-dist}
\item
We have
$ \tilde{\frakS}=\{ (r,\alpha\circ \sigma, \beta\circ \tau) |
    (r,\alpha,\beta)\in \frakS, \sigma\in \Sigma_k, \tau\in \Sigma_\ell\}$,
giving natural actions of $\Sigma_k\times \Sigma_\ell$ on $\tilde{\frakS}$ and on
$\{\Phi_{r,\alpha,\beta}(\vec{x},\vec{y}) \ |\ (r,\alpha,\beta)\in \tilde{\frakS}\}$.
Furthermore the later action is free.
\mlabel{it:phi-act}
\mlabel{it:phi-free}
\item
Fix an $r_0$ with $\max(k,\ell)\leq r_0 \leq k+\ell$. Let $$\tilde{\frakS}_{r_0}=\{ (r,\alpha,\beta)\in\tilde{\frakS}\ |\ r=r_0\}.$$
Then for any $\pi\in \Sigma_{r_0}$,
$ \pi(\Phi_{r_0,\alpha,\beta}(\vec{x},\vec{y}))
    =\Phi_{r_0, \pi^{-1}\circ\alpha, \pi^{-1}\circ\beta} (\vec{x},\vec{y}),$
giving a natural action of $\Sigma_{r_0}$ on $\tilde{\frakS}_{r_0}$.
Furthermore, this action is free.
\mlabel{it:phi-r}
\mlabel{it:phi-rfree}
\item
For $(r,\alpha,\beta)\in \tilde{\frakS}$, denote
$\Phi_{r,\alpha,\beta}(\vec{x},\vec{y})=(z_1,\cdots,z_r)$.
For $1\leq r'\leq r''\leq r$, denote $\vec{r}\subv=(r',r'+1,\cdots,r'')$ and $\vec{z}\,^\sharp=(z_{r'},z_{r'+1},\cdots,z_{r''})$.
Further denote
$$ \vec{k}\subv= \alpha^{-1}(\vec{r}\subv)=(k_{i_1},\cdots,k_{i_p}), \quad
 \vec{\ell}\subv= \beta^{-1}(\vec{r}\subv)=(\ell_{j_1},\cdots,\ell_{j_q}).$$
Then
$$ \vec{z}\,^\sharp=\Phi_{r\subv,\alpha|_{\vec{k}\subv}, \beta|_{\vec{\ell}\subv}}(\vec{x}\subv,\vec{y}\subv)$$
where the right hand side is defined by Eq.~(\mref{eq:phi1}).
\mlabel{it:phi-sub}
\end{enumerate}
\mlabel{lem:sigact}
\end{lemma}
In words, part~(\mref{it:phi-sub}) says that a part of a stuffle of
two vectors is a stuffle of parts of the two vectors.

\proofbegin
(\mref{it:phi-equal}) follows from the definition of $\Phi$ and the bijectivity of $\sigma$ and $\tau$:
\begin{eqnarray*}
\Phi_{r,\alpha,\beta}(\sigma(\vec{x}),\tau(\vec{y}))
&=&(x_{\sigma(\alpha^{-1}(1))}y_{\tau(\beta^{-1}(1))}, \cdots,
x_{\sigma(\alpha^{-1}(r))}y_{\tau(\beta^{-1}(r))})\\
&=& (x_{(\alpha\circ\sigma^{-1})^{-1}(1)}y_{(\beta\circ\tau^{-1})^{-1}(1)}, \cdots,
x_{(\alpha\circ\sigma^{-1})^{-1}(r)}y_{(\beta\circ\tau^{-1})^{-1}(r)}).
\end{eqnarray*}

\noindent
(\mref{it:phi-dist}) Suppose $\Phi_{r,\alpha,\beta}(\vec{x},\vec{y})
    =\Phi_{r',\alpha',\beta'}(\vec{x},\vec{y})$
for $(r,\alpha,\beta)$ and $(r',\alpha',\beta')$ in $\tilde{\frakS}$.
Then $r=r'$ by comparing the length of the vectors.
Further by the equation, and the fact that the $x_i$s and $y_j$s are distinct variables, we have, for $1\leq u\leq r$,
$ \alpha^{-1}(u)\neq \emptyset$ if and only if
$ \alpha'{}^{-1}(u) \neq \emptyset$. This implies $\alpha=\alpha'$
since both maps are injective with the same domain and codomain.
Similarly $\beta=\beta'$.

\noindent
(\mref{it:phi-act}) As is well-known, for an injective $\alpha:[k]\to [r]$, let $\im (\alpha) = \{n_1<\cdots<n_k\}$ and let
$\sigma(i)=\alpha^{-1}(n_i)$, $1\leq i\leq k$. Then $\sigma\in \Sigma_k$ and $\alpha\circ \sigma(i)=n_i$, $1\leq i\leq k$. Thus
$\alpha\circ \sigma$ is order preserving and injective. Similarly
for injective $\beta:[\ell]\to [r]$ we have $\tau\in \Sigma_{\ell}$ with $\beta\circ \tau: [\ell]\to [r]$ order preserving and injective. Then for $(r,\alpha,\beta)\in \tilde{\frakS}$ we have
$(r,\alpha,\beta)=(r,(\alpha\circ \sigma)\circ \sigma^{-1},
(\beta\circ \tau)\circ \tau^{-1})$, as needed.

The natural action of $(\sigma,\tau)\in \Sigma_k\times \Sigma_\ell$ on $\{\Phi_{r,\alpha,\beta}\ |\ (r,\alpha,\beta)\in \tilde{\frakS}\}$ is given by
$$\Phi_{r,\alpha,\beta}(\vec{x},\vec{y})^{(\sigma,\tau)}:
=\Phi_{r,\alpha,\beta}(\sigma(\vec{x}),\tau(\vec{y}))
=\Phi_{r,\alpha\circ \sigma^{-1},\beta\circ \tau^{-1}}(\vec{x},\vec{y})$$
by part (\mref{it:phi-equal}). So by part (\mref{it:phi-dist}),
$\Phi_{r,\alpha,\beta}(\vec{x},\vec{y})^{(\sigma,\tau)}
=\Phi_{r,\alpha,\beta}(\vec{x},\vec{y})$ implies
$(r,\alpha,\beta)=(r,\alpha\circ \sigma^{-1},\beta\circ \tau^{-1}).$
Thus $\alpha=\alpha\circ \sigma^{-1}$ and then $\sigma^{-1}=\id$ on $[k]$ since $\alpha$ is injective. Similarly, $\tau^{-1}=\id$ on $[\ell]$.

\noindent
(\mref{it:phi-r}) Fix a $(r_0,\alpha,\beta)\in \tilde{\frakS}_{r_0}$ and let $\Phi_{r_0,\alpha,\beta}(\vec{x},\vec{y})=(z_1,\cdots,z_{r_0})$.
Let $\pi\in \Sigma_{r_0}$. Then
$ \pi(\Phi_{r_0,\alpha,\beta}(\vec{x},\vec{y}))
    =(z_{\pi(1)},\cdots,z_{\pi(r_0)}).$
By the definition of $\Phi_{r_0,\alpha,\beta}$ in Eq.~(\mref{eq:Phi2}) and the bijectivity of $\pi$, we have
$$
z_{\pi(u)}= x_{\alpha^{-1}(\pi(u))}y_{\beta^{-1}(\pi(u))}
= x_{(\pi^{-1}\circ\alpha)^{-1}(u)}
y_{(\pi^{-1}\circ\beta)^{-1}(u)}.$$
So
$\pi(\Phi_{r_0,\alpha,\beta}(\vec{x},\vec{y}))$ is
$\Phi_{r_0, \pi^{-1}\circ\alpha, \pi^{-1}\circ\beta} (\vec{x},\vec{y})$.

Now let $\pi, \pi'\in \Sigma_{r_0}$. By item (\mref{it:phi-dist}),
$ \Phi_{r_0,\pi^{-1}\circ \alpha, \pi^{-1}\circ \beta}(\vec{x},\vec{y})
=\Phi_{r_0,\pi'{}^{-1}\circ \alpha, \pi'{}^{-1}\circ \beta}(\vec{x},\vec{y})$
implies $(r_0,\pi^{-1}\circ \alpha, \pi^{-1}\circ \beta)
    =(r_0,\pi'{}^{-1}\circ \alpha, \pi'{}^{-1}\circ \beta)$.
Thus $\pi^{-1}\circ \alpha =\pi'{}^{-1}\circ \alpha$ and
$\pi^{-1}\circ \beta =\pi'{}^{-1}\circ \beta$.
Since $\alpha$ and $\beta$ are injective and $\im\,\alpha\cup \im\, \beta = [r_0]$. We have $\pi=\pi'$.

\noindent
(\mref{it:phi-sub}) This follows directly from Eq.~(\mref{eq:Phi}) and Eq.~(\mref{eq:phi1}).
\proofend

\subsection{The proof of Theorem~\mref{thm:gshuf}}

As remarked at the beginning of this section, the quasi-shuffle product is also given by
$$ \vec{x} \msh \vec{y} = \sum_{(r,\alpha,\beta)\in \frakS } \Phi_{r,\alpha,\beta} (\vec{x},\vec{y}).$$
%
%
For
$\vec{a}\in \ZZ^k_{\leq 0}$ and $\vec{b}\in \ZZ^\ell_{\leq 0}$, only some of the stuffles $\Phi_{r,\alpha,\beta}(\vec{a},\vec{b})$ are distinct. Denote them by $\vec{d}_j, j\in J$. So we have
$ \vec{a}\msh \vec{b} =\sum_{j\in J} n_j \vec{d}_j.$
Denote
$$I_j=\{ (r,\alpha,\beta) |\ \Phi_{r,\alpha,\beta}(\vec{a},\vec{b})=\vec d_j\}.$$
Define the evaluation map
\begin{equation}
\hspace{-.45cm}f: \{ \Phi_{r,\alpha,\beta}(\vec{x},\vec{y}) | (r,\alpha,\beta)\in \frakS\}
\to \{ \Phi_{r,\alpha,\beta}(\vec{a},\vec{b}) |
(r,\alpha,\beta) \in \frakS\} = \{ \vec{d}_j | j\in J\}
\mlabel{eq:eval}
\end{equation}
by sending $x_i$ to $a_i$ and $y_j$ to $b_j$. Then
$ I_j =\{(r,\alpha,\beta)\ |\ \Phi_{r,\alpha,\beta}(\vec x,\vec y)\in f^{-1}(\vec{d}_j)\}$. So $|I_j|=n_j$.

Considering the actions of $\Sigma(\vec{a})\subset \Sigma_k$ and
$\Sigma(\vec{b})\subset \Sigma_\ell$, we have
\begin{eqnarray*}
\lefteqn{\hspace{-1cm}\sum _{\sigma \in \Sigma (\vec a), \tau \in \Sigma (\vec b)}(\vec
a *\vec b, \sigma (\vec x)*\tau (\vec y))
=\sum _{\sigma\in \Sigma(\vec{a}), \tau \in \Sigma(\vec{b})}
\sum_{(r,\alpha,\beta)\in \frakS}
\big (\Phi_{r,\alpha,\beta}(\vec a, \vec b), \Phi_{r,\alpha,\beta} (\sigma(\vec x), \tau (\vec y))\big)}
\\
&=& \sum _{\sigma,\tau} \sum_{j\in J}\sum _{(r,\alpha,\beta)\in I_j}
(\vec d_j, \Phi_{r,\alpha,\beta}(\sigma(\vec x), \tau (\vec y)))\\
&=& \sum _{j\in J} \sum
_{\sigma\in \Sigma(\vec{a}), \tau\in \Sigma(\vec{b}), (r,\alpha,\beta)\in I_j} (\vec d_j, \Phi_{r,\alpha,\beta}(\sigma(\vec x), \tau (\vec y)))
\end{eqnarray*}
For a fixed $\vec d_j$, the directions in the inner sum are $$S_j=\{\Phi_{r,\alpha,\beta}(\sigma(\vec{x}), \tau (\vec{y}))\ | \ (r,\alpha,\beta)\in I_j, \sigma\in \Sigma(\vec{a}), \tau\in \Sigma(\vec{b}) \}. $$
By its definition,
$S_j$ carries a $\Sigma (\vec a)\times \Sigma (\vec b)$ action.
Also by Lemma~\mref{lem:sigact}.(\mref{it:phi-free}),
$$|S_j|=|I_j| |\Sigma (\vec a)||\Sigma (\vec b)|=n_j|\Sigma (\vec a)||\Sigma (\vec b)|. $$
We next consider the action of $\Sigma (\vec d_j)$ and prove the
following lemma.

\begin{lemma} Let $r_0$ be the length of the vector $\vec{d}_j$.
The free action of $\Sigma_{r_0}$ on $$\{\Phi_{r_0,\alpha,\beta}(\vec{x},\vec{y})\ |\ (r_0,\alpha,\beta)\in \tilde{\frakS}_{r_0}\}$$
defined in Lemma~\mref{lem:sigact}.(\mref{it:phi-r}) restricts to a free action of $\Sigma(\vec d_j)$ on $S_j$
\mlabel{lem:sigd}
\end{lemma}

\proofbegin Once it is shown that $S_j$ is closed under the action of $\Sigma(\vec{d}_j)$, its freeness is automatic since it is the restriction of a free action.

To prove the closeness of the action,
let $\Phi_{r_0,\alpha,\beta}(\sigma(\vec{x}),\tau(\vec{y}))\in S_j$ and $\pi\in \Sigma(\vec{d}_j)$.
We only need to show that there are
$(r_0,\tilde{\alpha},\tilde{\beta})\in I_j$ and $(\tilde{\sigma},\tilde{\tau})\in \Sigma(\vec{a})\times
\Sigma(\vec{b})$ such that
$\pi\big(\Phi_{r_0,\alpha,\beta}(\sigma(\vec{x}),\tau(\vec{y}))\big)
=\Phi_{r_0,\tilde{\alpha},\tilde{\beta}}(\tilde{\sigma}(\vec{x}), \tilde{\tau}(\vec{y})).$
By Lemma~\mref{lem:sigact}.(\mref{it:phi-equal}) and (\mref{it:phi-r}), this means
$$
\Phi_{r_0,\pi^{-1}\circ\alpha\circ\sigma^{-1}, \pi^{-1}\circ\beta\circ\tau^{-1}}(\vec{x},\vec{y})
=\Phi_{r_0,\tilde{\alpha}\circ \tilde{\sigma}^{-1},\tilde{\beta}\circ\tilde{\tau}^{-1}}
(\vec{x}, \vec{y}).
$$
So we only need to prove
$
\pi^{-1}\circ\alpha\circ\sigma^{-1} = \tilde{\alpha}\circ \tilde{\sigma}^{-1}, \quad
 \pi^{-1}\circ\beta\circ\tau^{-1}
=\tilde{\beta}\circ\tilde{\tau}^{-1}.
$
That is, to show that the following diagram commutes.
\begin{equation}
\xymatrix{
[k] \ar[rr]^{\sigma^{-1}} \ar[d]^= && [k] \ar[rr]^\alpha &&
[r_0] \ar[d]^{\pi^{-1}} && [\ell] \ar[ll]_{\beta} &&
[\ell] \ar[ll]_{\tau^{-1}} \ar[d]^=\\
[k] \ar[rr]^{\tilde{\sigma}^{-1}} && [k] \ar[rr]^{\tilde{\alpha}} &&
[r_0] && [\ell] \ar[ll]_{\tilde{\beta}} && [\ell] \ar[ll]_{\tilde{\tau}^{-1}}
}
\mlabel{eq:sigdiag}
\end{equation}
Note that the existence of such $\tilde{\sigma}$ and $\tilde{\tau}$ with $(\tilde{\sigma},\tilde{\tau})\in \Sigma_k\times \Sigma_\ell$ is already given
in Lemma~\mref{lem:sigact}.(\mref{it:phi-act}). We want to show that $(\tilde{\sigma},\tilde{\tau})$ is in $\Sigma(\vec{a})\times \Sigma(\vec{b})$ for $\pi\in \Sigma(\vec{d}_j)$.

First let us look at the action of $\pi$ ``locally".
Let
$$\vec d_j=\Phi_{r_0,\alpha,\beta}(\vec{a},\vec{b})
=(d_1,\cdots,d_{r_0})$$
and let the sub-vector $\vec{d}\subv=(d_{r'},\cdots,d_{r''})$,
$1\leq r'\leq r'' \leq r_0$, be a 0-cluster of $\vec d_j$
as defined in
Definitin~\mref{de:cluster}. Denote $\vec{r}\subv=(r',\cdots,r'')$.
Let $\vec{k}\subv=\alpha^{-1}(\vec{r}\subv)$ and $\vec{\ell}\subv=\beta^{-1}(\vec{r}\subv)$. Since $\alpha$ and $\beta$ are
order preserving maps between vectors, we see that
$\vec{k}\subv=(k',k'+1,\cdots,k'')$ and $\vec{\ell}\subv=(\ell', \ell'+1,\cdots,\ell'')$ are sub-vectors of $[k]=(1,\cdots,k)$ and $[\ell]=(1,\cdots,\ell)$ respectively and
$$ \alpha\subv:=\alpha|_{\vec{k}\subv}, \quad \beta\subv:=\beta|_{\vec{\ell}\subv}$$
are order preserving maps with $\im(\alpha\subv)\cup \im(\beta\subv)=\vec{r}\subv$.

With the given $\sigma\in \Sigma(\vec a)$, $\tau\in \Sigma(\vec{b})$ and $\pi\in \Sigma_{\vec{r}\subv}\subset \Sigma(\vec{d}_j)$,
the same proof for Lemma~\mref{lem:sigact}.(\ref{it:phi-act}) shows that there is a bijection $\sigma\subv:\sigma(\vec{k}\subv)\to \vec{k}\subv$ and an order preserving injection $\alpha\subv: \vec{k}\subv \to \vec{r}\subv$ such that
$ \pi^{-1}\circ\alpha\circ\sigma^{-1}|_{\sigma(\vec{k}\subv)}
    = \alpha\subv \circ \sigma\subv,$
that is, the left square of the following diagram is commutative
\begin{equation}
\xymatrix{
\sigma(\vec{k}\subv) \ar[rr]^{\sigma^{-1}} \ar[d]^= && \vec{k}\subv \ar[rr]^\alpha &&
\vec{r}\subv \ar[d]^{\pi^{-1}} && \vec{\ell}\subv \ar[ll]_{\beta} &&
\tau(\vec{\ell}\subv) \ar[ll]_{\tau^{-1}} \ar[d]^=\\
\sigma(\vec{k}\subv) \ar[rr]^{\sigma\subv} && \vec{k}\subv \ar[rr]^{\alpha\subv} &&
\vec{r}\subv && \vec{\ell}\subv \ar[ll]_{\beta\subv} && \tau(\vec{\ell}\subv) \ar[ll]_{\tau\subv}
}
\mlabel{eq:sigdiag1}
\end{equation}
Similarly, there is a bijection $\tau\subv:\tau(\vec{\ell}\subv)\to \vec{\ell}\subv$ and an order preserving injection $\beta\subv: \vec{\ell}\subv \to \vec{r}\subv$ such that
the right square of the diagram is commutative.

Since our choice of $\pi$ is the identity when restricted to $[r_0]\backslash \vec{r}\subv$, we have the trivial commutative diagram of bijections and order preserving maps
\begin{equation}
\xymatrix{
[k]\backslash \sigma(\vec{k}\subv) \ar[r]^{\sigma^{-1}} \ar[d]^= & [k]\backslash \vec{k}\subv \ar[r]^\alpha &
[r_0]\backslash \vec{r}\subv \ar[d]^{\pi^{-1}} & [\ell]\backslash \vec{\ell}\subv \ar[l]_{\beta} &
[\ell]\backslash \tau(\vec{\ell}\subv) \ar[l]_{\tau^{-1}} \ar[d]^=\\
[k]\backslash \sigma(\vec{k}\subv) \ar[r]^{\sigma^{-1}} & [k]\backslash \vec{k}\subv \ar[r]^{\alpha} &
[r_0]\backslash \vec{r}\subv & [\ell]\backslash \vec{\ell}\subv \ar[l]_{\beta} & [\ell]\backslash \tau(\vec{\ell}\subv) \ar[l]_{\tau^{-1}}
}
\mlabel{eq:sigdiag2}
\end{equation}
Taking the union of these two diagrams, we obtain a commutative diagram~(\mref{eq:sigdiag}) where the bijections and order preserving maps are defined by
$$ \tilde{\sigma}(i)=\left \{ \begin{array}{ll}
    \sigma\subv^{-1}(i), & i\in \vec{k}\subv, \\
    \sigma(i), & i\not\in \vec{k}\subv
    \end{array} \right . \quad
     \tilde{\alpha}(i)=\left \{ \begin{array}{ll}
    \alpha\subv(i), & i\in \vec{k}\subv, \\
    \alpha(i), & i \not\in \vec{k}\subv
    \end{array} \right . \quad 1 \leq i\leq k.
    $$
$$ \tilde{\tau}(i)=\left \{ \begin{array}{ll}
    \tau\subv^{-1}(i), & i\in \vec{\ell}\subv, \\
    \tau(i), & i\not\in \vec{\ell}\subv
    \end{array} \right . \quad
     \tilde{\beta}(i)=\left \{ \begin{array}{ll}
    \beta\subv(i), & i\in \vec{\ell}\subv, \\
    \beta(i), & i \not\in \vec{\ell}\subv
    \end{array} \right . \quad 1 \leq i\leq \ell.
    $$

We next show that $\tilde{\sigma}\in \Sigma(\vec{a})$ and $\tilde{\tau}\in \Sigma(\vec{b})$.
Let $\vec{z}\subv=(z_{r'},\cdots,z_{r''})$ be the the sub-vector of
$\Phi_{r_0,\alpha,\beta}(\vec{x},\vec{y})
=(z_1,\cdots,z_{r_0})$ corresponding to $\vec{r}\subv$.
Similarly, denote
$$\vec{x}\subv=(x_{k'},\cdots,x_{k''}),
\vec{y}\subv=(y_{\ell'},\cdots,y_{\ell''}),
\vec{a}\subv=(a_{k'},\cdots,a_{k''}),
\vec{b}\subv=(b_{\ell'},\cdots,b_{\ell''}).$$
By Lemma~\mref{lem:sigact}.(\mref{it:phi-sub}), we have
$$\vec{z}\subv = \Phi_{r\subv,\alpha\subv,\beta\subv}(\vec{x}\subv,\vec{y}\subv),
\quad
\vec{d}\subv = \Phi_{r\subv,\alpha\subv,\beta\subv}(\vec{a}\subv,\vec{b}\subv).$$
Here $r\subv=r''-r'+1$.

Under the evaluation map~(\mref{eq:eval}),
 $$f(\vec{z}\subv)
=f(\Phi_{r\subv,\alpha\subv,\beta\subv}(\vec{x}\subv, \vec{y}\subv))= \Phi_{r\subv,\alpha\subv,\beta\subv}(\vec{a}\subv, \vec{b}\subv)= \vec{d}\subv=\vec{0},$$ and
$f(x_i)=a_i, f(y_j)=b_j, f(x_iy_j)=a_i+b_j$ are non-positive numbers or their sums. So we must have
$\vec{a}\subv=\vec{0}$ and $\vec{b}\subv=\vec{0}.$
Thus $\vec{a}\subv$ (resp. $\vec{b}\subv$) is a part of a zero cluster of $\vec{a}$ (resp. $\vec{b}$).
Thus $\vec{k}\subv\subseteq \vec{k}^{(i)}\subseteq [k]$ where $\vec{k}^{(i)}$ is the index set of a zero cluster of $\vec{a}$
(see Definition~(\mref{de:cluster})).
Since $\Sigma(\vec{a})|_{\vec{k}^{(i)}}=\Sigma_{\vec{k}^{(i)}}$,
we have
$\sigma(\vec{k}\subv)\subseteq \sigma(\vec{k}^{(i)})
=\vec{k}^{(i)}.$
Further, there is $\tilde{\sigma}'\in \Sigma_{\vec{k}^{(i)}}$ such that
$\tilde{\sigma}'|_{\vec{k}\subv}=\sigma\subv{}^{-1}$ and
$\tilde{\sigma}'|_{\vec{k}^{(i)}\backslash\vec{k}\subv}=\sigma.$
Since $\Sigma(\vec{a})=\Sigma_{\vec{k}^{(i)}}\times \Sigma''$,
we further have $\tilde{\sigma}=(\tilde{\sigma}',\sigma'')\in \Sigma(\vec{a}).$ Similarly, $\tilde{\tau}\in \Sigma(\vec{b})$.

Thus $S_j$ is closed under the action of $\Sigma(\vec{r}\subv)$ and hence of $\Sigma(\vec{d}_j)$
since $\Sigma (\vec d_j)$ is the direct product of such $\Sigma(\vec{r}\subv)$ from all the zero clusters of $\vec{d}_j$.
\proofend

By Lemma~\mref{lem:sigd}, there are $n_j|\Sigma (\vec a)||\Sigma (\vec b)|/|\Sigma (d_j)|$ orbits. Let
$$
\vec r_h, \ h=1, \cdots, n_j|\Sigma (\vec a)||\Sigma (\vec b)|/|\Sigma (d_j)|
$$
be a complete set of representatives of $S_j$, then
$$\sum _{\sigma \in \Sigma (\vec a), \tau \in \Sigma (\vec b)}(\vec
a *\vec b, \sigma (\vec x)*\tau (\vec y))=\sum _{j} \sum _h(\vec
d_j, \vec r_h)^{(\Sigma (\vec d_j))}
$$
Thus
\allowdisplaybreaks{
\begin{eqnarray*}
\gzeta(\vec a)\gzeta(\vec b)
&=&\frac 1{|\Sigma (\vec a)||\Sigma
(\vec b)|}\lim _{\vec x \to -\vec a, \vec y\to -\vec b}
\zeta\lp\wvec{\vec a}{\vec x}\rp^{(\Sigma (\vec a))}
\zeta\lp\wvec{\vec b}{\vec y}\rp^{(\Sigma (\vec b))}
\\
&=& \frac 1{|\Sigma (\vec a)||\Sigma
(\vec b)|}\lim _{\vec x \to -\vec a, \vec y\to -\vec b}\sum _{\sigma
\in \Sigma (\vec a), \tau \in \Sigma (\vec b)}
\zeta\lp\wvec{\vec a *\vec b}{\sigma (\vec x)*\tau (\vec y)}\rp
\\
&=& \frac 1{|\Sigma (\vec a)||\Sigma
(\vec b)|}\lim _{\vec x \to -\vec a, \vec y\to -\vec b}\sum _{j}
\sum _h\zeta\lp\wvec{\vec d_j}{\vec r_h}\rp^{(\Sigma (\vec d_j))}
\\
&=& \frac 1{|\Sigma (\vec a)||\Sigma
(\vec b)|}\sum _{j}\sum _h\lim _{\vec x \to -\vec a, \vec y\to -\vec
b} \zeta\lp\wvec{\vec d_j}{\vec r_h}\rp^{(\Sigma (\vec d_j))}.
\end{eqnarray*}
}
Note that, for a fixed $h$, $\vec x \to -\vec a, \vec y\to -\vec b$
means $\vec r_h \to -\vec d_j$. Then
$$\lim _{\vec x \to -\vec a, \vec y\to -\vec b}
\zeta\lp\wvec{\vec d_j}{\vec r_h}\rp^{(\Sigma (\vec d_j))}
= \lim _{r_h\to -\vec d_j}
\zeta\lp\wvec{\vec d_j}{\vec r_h}\rp^{(\Sigma (\vec d_j))}.$$
So by Theorem~\mref{thm:gzeta}, we have
$$\gzeta(\vec a)\gzeta(\vec b)=\frac 1{|\Sigma (\vec a)||\Sigma
(\vec b)|}\sum _{j} n_j|\Sigma (\vec a)||\Sigma (\vec b)|\gzeta
(\vec d_j)=\sum _j n_j\gzeta (\vec d_j)=\gzeta (\vec a *\vec b).
$$
This completes the proof of Theorem~\mref{thm:gshuf}.



\begin{thebibliography}{abcdsfgh}

\mbibitem{A-H} M. Aguiar and S. Hsiao, Canonical characters on
quasi-symmetric functions and bivariate Catalan numbers,
{\em Electron. J. Combin.} {\bf 11(2)} (2005), \#R15, 34pp.

\mbibitem{AET} S. Akiyama, S. Egami and Y. Tanigawa,
    Analytic continuation of multiple zeta-functions
    and their values at non-positive integers,
    {\it Acta Arith.} {\bf 98} (2001), 107--116.

\mbibitem{A-T} S Akiyama and Y. Tanigawa,
    Multiple zeta values at non-positive integers,
    {\it Ramanujan J.} {\bf 5} (2001), 327--351.

\mbibitem{A-K} T. Arakawa and M. Kaneko,
    Multiple zeta values, poly-Bernoulli numbers, and related zeta
functions, {\it Nagoya Math. J.} {\bf 153} (1999), 189--209.


\mbibitem{3BL} J.~M.~Borwein, D.~J.~Broadhurst, D.~M.~Bradley,
and P.~Lison\v ek,
  {Special values of multiple polylogarithms},
  \textit{Trans.\ Amer.\ Math.\ Soc.},  \textbf{353}, (2001), no.~3, 907--941.
       {arXiv: {math.CA/9910045}}

\mbibitem{Br1} D.~M. Bradley, Partition identities for the multiple zeta function,
in ``Zeta Functions, Topology, and Quantum Physics" T. Aoki et. al. (eds.),
Springer, 2005, 19--29.

\mbibitem{Br} D.~M.~Bradley, {Multiple $q$-zeta values},
   {\it{J. Algebra}}, {\bf{283}}, (2005), no. 2, 752--798
                                  {arXiv: {math.QA/0402093}}

\mbibitem{B-K} D.~J.~Broadhurst and D.~Kreimer,
  {Association of multiple zeta values with positive knots via Feynman diagrams up to $9$ loops},
                        {\it Phys. Lett. B}, {\bf{393}}, (1997), no. 3-4, 403--412.

\mbibitem{B-Z} J.-L. Bryliski, B. Zhang, {Equivariant Todd classes
for toric varieties,} arXiv:math.AT/0311318.

\mbibitem{Ca2} P.~Cartier, {On the structure of free Baxter
algebras},
   {\it Adv. in Math.}, {\bf 9}, (1972), 253--265.

\mbibitem{Ca1} P.~Cartier, {Fonctions polylogarithmes, nombres
polyzêtas et groupes pro-unipotents},
   {\it Astérisque}, {\bf{282}}, (2002), 137--173, (Sém. Bourbaki no. 885).

\mbibitem{C-K1} A.~Connes and D.~Kreimer, {Renormalization in
quantum field theory and the Riemann-Hilbert problem. I. The Hopf algebra structure of graphs and the main theorem},
{\it Comm. Math. Phys.}, {\bf 210}, (2000), no. 1, 249--273.

\mbibitem{C-K2} A.~Connes and D.~Kreimer, {Renormalization in
quantum field theory and the Riemann-Hilbert problem. II. The $\beta$-function, diffeomorphisms and the renormalization group},
 {\it Comm. Math. Phys.}, {\bf 216}, (2001), no. 1, 215--241.

\mbibitem{C-M}
A. Connes and M. Marcolli, From physics to
number theory via noncommutative geometry,
  part II: renormalization, the Riemann-Hilbert correspondence,
  and motivic Galois theory, {arXiv:hep-th/0411114}.

\mbibitem{Dr} V. G.~Drinfel'd, {On quasitriangular quasi-Hopf algebras and a group closely related to ${\rm Gal}(\bar{\QQ}/\QQ)$,}
    {\em Leningrad Math. J.} {\bf 2} (1991), 829--860.


\mbibitem{E-G1} K.~Ebrahimi-Fard and L. Guo {Mixable Shuffles,
Quasi-shuffles and Hopf Algebras}, {\em J. Alg. Combinatorics}, {\bf 24}
(2006), 83--101, {arXiv: {math.RA/0506418}}.


\mbibitem{E-G5} K.~Ebrahimi-Fard and L.~Guo, Multiple zeta values
and Rota-Baxter algebras, {arXiv:math.NT/0601558}.



\mbibitem{EGK2} K.~Ebrahimi-Fard, L.~Guo and D.~Kreimer,
{Integrable Renormalization II: the General case,}
 {Annales Henri Poincar\'e}, {\bf{6}}, (2005), 369--395.

\mbibitem{EGK3} K.~Ebrahimi-Fard, L.~Guo and D.~Kreimer,
  {Spitzer's Identity and the Algebraic Birkhoff Decomposition in pQFT},
  {\it{J. Phys. A: Math. Gen.}}, {\bf{37}}, (2004) 11037--11052.
     {arXiv: {hep-th/0407082}}

\mbibitem{Eh} R. Ehrenborg, On postes and Hopf algebras, {\em Adv. Math.} {\bf 119} (1996), 1-25.

\mbibitem{Fa} F. Fares, Quelques constructions d'alg\`ebres et de coalg\`ebres, Universit\'e du Qu\'ebec \`a
Mont\'eal (1999).

\bibitem{F-G} H.~Figueroa, J.~M.~Gracia-Bond\'ia,
 {{Combinatorial Hopf algebras in quantum field theory I}},
 {\it Reviews of Mathematical Physics}, {\bf{17}} (2005), 881--976.
 arXiv: {arXiv:hep-th/0408145}

\mbibitem{Go3} A.~G.~Goncharov, {Periods and mixed motives},
     {arXiv: {math.AG/0202154}}.

\mbibitem{G-M} A. Goncharov and Y. Manin,
    {Multiple $\zeta$-motives and moduli spaces $\overline{\mathcal M}_{0,n}$},
    {\em Comp. Math.} {\bf 140} (2004), 1 - 14.



\mbibitem{Gu3} L.~Guo, {Baxter algebras, Stirling numbers and
partitions},
  {\it{J. Algebra Appl.}}, {\bf{4}}, (2005), no. 2, 153--164.
        {arXiv: {math.AC/0402348}}

\mbibitem{G-K1} L.~Guo and W.~Keigher, {Free Baxter algebras and
shuffle products},
   {\it Adv. in Math.,} {\bf 150}, (2000), 117--149.

\mbibitem{G-K2} L.~Guo and W.~Keigher, {On Baxter algebras:
completions and the internal construction},
  {\it Adv. in Math.}, {\bf{151}}, (2000), 101--127.

\mbibitem{Ha} M. Hazewinkel, {Generalized overlapping shuffle algebras
{\em J. Math. Sci. (New York)}, {\bf 106} (2001), 3168-3186.

\mbibitem{Ho0} M.~E.~Hoffman, {Multiple harmonic series},
  {\em Pacific J. Math.}, {\bf 152} (1992), no. 2, 275--290.

\mbibitem{Ho1} M.~E.~Hoffman, {The algebra of multiple harmonic
series},
     {\it J. Algebra}, {\bf 194}, no. 2, (1997), 477--495.

\mbibitem{Ho2} M.~E.~Hoffman, {Quasi-shuffle products},
   {\it J. Algebraic Combin.}, {\bf 11}, no. 1, (2000), 49--68.

\mbibitem{Ho3} M.~E.~Hoffman, {Algebraic aspects of multiple zeta
values},  in ''Zeta Functions, Topology and Quantum Physics"
    T.~Aoki et. al. (Eds.),
    Springer, (2005), 51--73, arXiv:{{math.QA/0309425}}

\mbibitem{IKZ} K.~Ihara, M.~Kaneko and D.~Zagier,
  {Derivation and double shuffle relations for multiple zeta values}, Compos. Math. {\bf 142} (2006), 307--338.


\mbibitem{Kr} D.~Kreimer, {Knots and Feynman Diagrams},
  Cambridge Lecture Notes in Physics, {\bf{13}}.
  Cambridge University Press, Cambridge, (2000).

\mbibitem{Lo} J.-L. Loday, On the algebra of quais-shuffles, arXiv: math.QA/0506498.

\mbibitem{Ma} D. Manchon,
  {{Hopf algebras, from basics to applications to renormalization}},
  Comptes-rendus des Rencontres mathm\'ematiques de Glanon 2001.
  ArXiv:math.QA/0408405

\mbibitem{M-P1} D. Manchon and S. Paycha, Shuffle relations for regularized integrals of symbols, arXiv:math-ph/0510067.

\mbibitem{M-P2} D. Manchon and S. Paycha, Renormalized Chen integrals for symbols on $\RR^n$ and renormlized polyzeta functions, arXiv:math.NT/0604562.


\mbibitem{Ma2} K. Matsumoto,
    The analytic continuation and the asymptotic behaviour of certain multiple zeta-functions I, {\it J. Number Theory}, {\bf 101} (2003), 223--243.


\mbibitem{Ra} G. Racinet, {Doubles m\'elanges des polylogarithmes multiples aux racines de l'unit\'e}, {\em Pub. Math. IHES}, {\bf 95} (2002), 185-231.


\mbibitem{Ro3} G.-C.~Rota and D.~A.~Smith, {Fluctuation theory and
Baxter algebras},
  Symposia Mathematica, Vol. IX
 (Convegno di Calcolo delle Probabilità, INDAM, Rome, 1971),
  pp. 179--201. Academic Press, London, (1972).



\mbibitem{Te} T.~Terasoma, {Mixed Tate motives and multiple zeta values,}
                        {\it Invent.~Math.}, {\bf{149}}, (2002), no. 2, 339--369.
                        {arXiv: {math.AG/0104231}}

\mbibitem{Za} D.~Zagier, {Values of zeta functions and their
applications},
  First European Congress of Mathematics, Vol. II (Paris, 1992), 497--512,
{\it Progr. Math.}, {\bf{120}}, Birkhäuser, Basel, 1994

\mbibitem{Zh2} J. Zhao, Analytic continuation of multiple zeta
    functions. {\it Proc. Amer. Math. Soc.} {\bf 128} (2000),
        1275-1283.

}


\end{thebibliography}
\end{document}